\documentclass[11pt]{article}

\usepackage{amsfonts}
\usepackage{amscd}
\usepackage{amssymb}
\usepackage{amsthm}
\usepackage{amsmath, xspace}
\usepackage{enumitem}
\usepackage{fancyhdr}
\usepackage{mathrsfs}
\usepackage[plainpages,backref]{hyperref}
\usepackage{lscape}

\theoremstyle{plain}
\newtheorem{thm}{Theorem}[section]

\newtheorem{prop}[thm]{Proposition}
\newtheorem{cor}[thm]{Corollary}
\theoremstyle{definition}
	
\newtheorem{remark}[thm]{Remark}
\newtheorem{example}{Example}[section]

\theoremstyle{remark}

\numberwithin{equation}{section}

\providecommand{\keywords}[1]{\textbf{\textit{Key words---}} #1}

\def\CC{\mathbb{C}}

\def\KK{\mathbb{K}}

\def\ZZ{\mathbb{Z}}

\def\fa{\mathfrak{a}}

\def\Card{\mathrm{Card}}

\def\ev{\mathrm{ev}}

\def\Ind{\mathrm{Ind}}

\usepackage{etex} %fixes the fight that pictex has with every other drawing package
\usepackage{pictexwd}
\usepackage{tikz}
	\usepgflibrary[patterns] % ConTEXt and pure pgf 
	\usetikzlibrary{patterns} % LATEX and plain TEX when using TikZ 
	\usepgflibrary{shapes.geometric}
\usetikzlibrary{arrows}

\newcounter{r}
\newcounter{s}

\newcommand\Part[1]{
        \setcounter{r}{1}
	 \foreach \x in {#1}{
 	{\ifnum\value{r}=1
		\draw (0,\value{r}-1)--(\x,\value{r}-1); 
		\fi}
	\draw (0,\value{r}) to (\x,\value{r});
   	\foreach \y in {0, ..., \x} {\draw (\y,\value{r})--(\y,\value{r}-1);}
	\addtocounter{r}{1}
 }}
 
\newcommand\Tableau[1]{
        \foreach \x [count = \c from 1] in {#1} {
		\foreach \y [count = \d from 1] in \x{
			\node at (\d-.5,\c-.5) {\scriptsize$\y$}; %Want to change the color of the numbers? insert \color{???} just after \scriptsize
			\draw (\d,\c) to (\d,\c-1);
			{\ifnum\d=1
				\draw (0,\c) to (0,\c-1);
				\fi}
			\setcounter{r}{\d}
		}
		{\ifnum\c=1
			\draw (0,0)--(\value{r},0);
			\fi}
		\draw(0,\c) to (\value{r},\c);
		\setcounter{s}{\c}}}
		
\newcommand\sTableau[1]{
        \foreach \x [count = \c from 1] in {#1} {
		\foreach \y [count = \d from 1] in \x{
			\node at (\d-.5,\c-.5) {\tiny$\y$}; %Want to change the color of the numbers? insert \color{???} just after \scriptsize
			\draw (\d,\c) to (\d,\c-1);
			{\ifnum\d=1
				\draw (0,\c) to (0,\c-1);
				\fi}
			\setcounter{r}{\d}
		}
		{\ifnum\c=1
			\draw (0,0)--(\value{r},0);
			\fi}
		\draw(0,\c) to (\value{r},\c);
		\setcounter{s}{\c}}}
\tikzstyle{V}=[draw, fill =black, circle, inner sep=0pt, minimum size=1.5pt]
\tikzstyle{wV}=[draw, fill =white, circle, inner sep=0pt, minimum size=4.5pt]
\tikzstyle{bV}=[draw, fill =black, circle, inner sep=0pt, minimum size=4.5pt]
\tikzstyle{over}=[draw=white,double=black,line width=2pt, double distance=.5pt]

\def\Over[#1,#2][#3,#4]{ %1,2=start position; 3,4=end position
	\draw[style=over]   (#2,#1) .. controls ++(#4*.5-#2*.5,0) and ++(-#4*.5+#2*.5,0) .. (#4,#3);}
\def\Under[#1,#2][#3,#4]{ %1,2=start position; 3,4=end position
	\draw  (#2,#1) .. controls ++(#4*.5-#2*.5,0) and ++(-#4*.5+#2*.5,0) .. (#4,#3);}
\def\Cross[#1,#2][#3,#4]{%Mimic over, under follows
	\Under[#3,#2][#1,#4]\Over[#1,#2][#3,#4]}

\def\Tops[#1][#2][#3]{%1=pole locations, 2=top, 3=k 
	\foreach\x in {#1}{
		\draw (#2,\x+.15) -- (#2+.1, \x+.15) (#2, \x-.15) -- (#2+.1, \x-.15) ;
		\draw (#2+.1,\x) arc (0:360:.75mm and 1.5mm);}
	\foreach \x in {1,...,#3} {\draw (#2,\x)  to (#2+.05,\x); \node[V] at (#2+.05,\x){};}
	}
\def\Bottoms[#1][#2][#3]{%1=pole locations, 2 = bottom, 3=top 
	\foreach\x in {#1}{
		\draw (#2, \x+.15) -- (#2-.1, \x+.15) (#2, \x-.15) -- (#2-.1, \x-.15) ;
		\draw (#2-.1, \x+.15) arc (90:270:.75mm and 1.5mm);}
	\foreach \x in {1,...,#3} {\draw (#2, \x)  to (#2-.05, \x); \node[V] at (#2-.05, \x){};}
	}
\def\Caps[#1][#2,#3][#4]{%1=pole locations, 2 = bottom, 3=top, 4=k 
	\Tops[#1][#3][#4]
	\Bottoms[#1][#2][#4]
	}
\def\Pole[#1][#2,#3]{%1=horizontal location, 2 = bottom, 3=top
	\shade[left color=white,right color=white] (#2,#1+.15) rectangle (#3,#1-.15);
	\draw[over] (#2,#1+.15) to (#3,#1+.15) (#2,#1-.15) to (#3,#1-.15) ;}
\def\Label[#1,#2][#3][#4]{%1,2 = top/bot position, 3=i, 4=label
	\node[right] at (#2+.1,#3) {#4};
	\node[left] at (#1-.1,#3) {#4};		}
\def\Nodes[#1][#2]{
	 \foreach \x in {1,...,#2} {\node[V] at (#1,\x){};	}
	}
\def\PoleCaps[#1][#2,#3]{%1=pole location, 2 = bottom, 3=top, 4=k 
	\foreach\x in {#1}{
		\draw (#2,\x+.15) -- (#2-.1,\x+.15) (#2,\x-.15) -- (#2-.1,\x-.15) ;
		\draw (#2-.1,\x+.15) arc (0:-180:1.5mm and .75mm);}
	\foreach\x in {#1}{
		\draw (#3,\x+.15) -- (#3+.1,\x+.15) (#3,\x-.15) -- (#3+.1,\x-.15) ;
		\draw (#3+.1,\x+.15) arc (0:360:1.5mm and .75mm);}
	}
\def\PoleTwist[#1,#2]{%1 = bottom, 2=top
	\foreach \x/\y in {-1/1L, -.7/1R, 0/2L, .3/2R}{\coordinate(T\y) at (#2,\x); \coordinate(B\y) at (#1,\x);}
	\draw[thin] (B1R) .. controls ++(#2*.5-#1*.5-.1,0) and ++(-#2*.5+#1*.5-.1,0) ..  (T2R)
			(B1L)   .. controls ++(#2*.5-#1*.5+.1,0) and ++(-#2*.5+#1*.5+.1,0) ..    (T2L) ;
	\draw[line width=2pt, white]
			(#1,.15)  .. controls +(#2*.5-#1*.5,0) and +(-#2*.5+#1*.5,0) ..   (#2,-.85) ;
	\draw[thin,over] 
		(B2R) .. controls ++(#2*.5-#1*.5+.1,0) and ++(-#2*.5+#1*.5+.1,0) ..  (T1R) 
			(B2L)  .. controls +(#2*.5-#1*.5-.1,0) and +(-#2*.5+#1*.5-.1,0) ..   (T1L) ;
			}

\def\SymPolesCaps[#1,#2][#3]{%1 = Vertical position, 2=k
	\draw (#1,.3) -- (#1-.1,.3) (#1,.15) -- (#1-.1, .15) ;
	\draw (#1-.1, .3) arc (0:-180:2pt and 1.5pt);
	\draw (#1,#3+.7) -- (#1-.1,#3+.7) (#1,#3+.85) -- (#1-.1,#3+.85) ;
	\draw (#1-.1,#3+.85)  arc (0:-180:2pt and 1.5pt);
	\draw (#2,.3) -- (#2+.1, .3) (#2, .15) -- (#2+.1, .15) ;
	\draw (#2+.1, .3) arc (0:360:2pt and 1.5pt);
	\draw (#2, #3+.7) -- (#2+.1, #3+.7) (#2, #3+.85) -- (#2+.1, #3+.85) ;
	\draw (#2+.1, #3+.85) arc (0:360:2pt and 1.5pt);}

\newcommand{\posleq}[1]{
	\hspace{0.1cm}
	\begin{tikzpicture}
	\draw (-0.8ex, -0.5ex) -- (0.8ex, -0.5ex);
	\draw (-0.8ex, 0.4ex) -- (0.7ex, -0.2ex);
	\draw (-0.8ex, 0.4ex) -- (0.7ex, 1ex);
	\draw (0.4ex,0.4ex) --(1.1ex, 0.4ex);
	\draw (0.75ex,0.75ex) --(0.75ex, 0.05ex);
	\end{tikzpicture}
	\hspace{0.1cm}
	}
\newcommand{\negleq}[1]{
	\hspace{0.1cm}
	\begin{tikzpicture}
	\draw (-0.8ex, -0.5ex) -- (0.8ex, -0.5ex);
	\draw (-0.8ex, 0.4ex) -- (0.7ex, -0.2ex);
	\draw (-0.8ex, 0.4ex) -- (0.7ex, 1ex);
	\draw (0.4ex,0.4ex) --(1.1ex, 0.4ex);
	\end{tikzpicture}
	\hspace{0.1cm}
	}
	
\newcommand{\zeroleq}[1]{
	\hspace{0.1cm}
	\begin{tikzpicture}
	\draw (-0.8ex, -0.5ex) -- (0.8ex, -0.5ex);
	\draw (-0.8ex, 0.4ex) -- (0.7ex, -0.2ex);
	\draw (-0.8ex, 0.4ex) -- (0.7ex, 1ex);
	\draw  (0.75ex,0.4ex) ellipse (0.2ex and 0.35ex);
	\end{tikzpicture}
	\hspace{0.1cm}
	}
	
\newcommand{\posgeq}[1]{
	\hspace{0.1cm}
	\begin{tikzpicture}
	\draw (-0.8ex, -0.5ex) -- (0.8ex, -0.5ex);
	\draw (0.8ex, 0.4ex) -- (-0.7ex, -0.2ex);
	\draw (0.8ex, 0.4ex) -- (-0.7ex, 1ex);
	\draw (-0.4ex,0.4ex) --(-1.1ex, 0.4ex);
	\draw (-0.75ex,0.75ex) --(-0.75ex, 0.05ex);
	\end{tikzpicture}
	\hspace{0.1cm}
	}
\newcommand{\neggeq}[1]{
	\hspace{0.1cm}
	\begin{tikzpicture}
	\draw (-0.8ex, -0.5ex) -- (0.8ex, -0.5ex);
	\draw (0.8ex, 0.4ex) -- (-0.7ex, -0.2ex);
	\draw (0.8ex, 0.4ex) -- (-0.7ex, 1ex);
	\draw (-0.4ex,0.4ex) --(-1.1ex, 0.4ex);
	\end{tikzpicture}
	\hspace{0.1cm}
	}
	
\newcommand{\zerogeq}[1]{
	\hspace{0.1cm}
	\begin{tikzpicture}
	\draw (-0.8ex, -0.5ex) -- (0.8ex, -0.5ex);
	\draw (0.8ex, 0.4ex) -- (-0.7ex, -0.2ex);
	\draw (0.8ex, 0.4ex) -- (-0.7ex, 1ex);
	\draw  (-0.75ex,0.4ex) ellipse (0.2ex and 0.35ex);
	\end{tikzpicture}
	\hspace{0.1cm}
	}

\newcommand{\posl}[1]{
	\hspace{0.1cm}
	\begin{tikzpicture}
	\draw (-0.8ex, 0.4ex) -- (0.7ex, -0.2ex);
	\draw (-0.8ex, 0.4ex) -- (0.7ex, 1ex);
	\draw (0.4ex,0.4ex) --(1.1ex, 0.4ex);
	\draw (0.75ex,0.75ex) --(0.75ex, 0.05ex);
	\end{tikzpicture}
	\hspace{0.1cm}
	}
\newcommand{\negl}[1]{
	\hspace{0.1cm}
	\begin{tikzpicture}
	\draw (-0.8ex, 0.4ex) -- (0.7ex, -0.2ex);
	\draw (-0.8ex, 0.4ex) -- (0.7ex, 1ex);
	\draw (0.4ex,0.4ex) --(1.1ex, 0.4ex);
	\end{tikzpicture}
	\hspace{0.1cm}
	}
	
\newcommand{\zerol}[1]{
	\hspace{0.1cm}
	\begin{tikzpicture}
	\draw (-0.8ex, 0.4ex) -- (0.7ex, -0.2ex);
	\draw (-0.8ex, 0.4ex) -- (0.7ex, 1ex);
	\draw  (0.75ex,0.4ex) ellipse (0.2ex and 0.35ex);
	\end{tikzpicture}
	\hspace{0.1cm}
	}
	
\newcommand{\posg}[1]{
	\hspace{0.1cm}
	\begin{tikzpicture}
	\draw (0.8ex, 0.4ex) -- (-0.7ex, 1ex);
	\draw (0.8ex, 0.4ex) -- (-0.7ex, -0.2ex);
	\draw (-0.4ex,0.4ex) --(-1.1ex, 0.4ex);
	\draw (-0.75ex,0.75ex) --(-0.75ex, 0.05ex);
	\end{tikzpicture}
	\hspace{0.1cm}
	}
\newcommand{\negg}[1]{
	\hspace{0.1cm}
	\begin{tikzpicture}
	\draw (0.8ex, 0.4ex) -- (-0.7ex, -0.2ex);
	\draw (0.8ex, 0.4ex) -- (-0.7ex, 1ex);
	\draw (-0.4ex,0.4ex) --(-1.1ex, 0.4ex);
	\end{tikzpicture}
	\hspace{0.1cm}
	}
	
\newcommand{\zerog}[1]{
	\hspace{0.1cm}
	\begin{tikzpicture}
	\draw (0.8ex, 0.4ex) -- (-0.7ex, -0.2ex);
	\draw (0.8ex, 0.4ex) -- (-0.7ex, 1ex);
	\draw  (-0.75ex,0.4ex) ellipse (0.2ex and 0.35ex);
	\end{tikzpicture}
	\hspace{0.1cm}
	}

\makeatletter
\renewcommand{\@makefnmark}{\mbox{\textsuperscript{}}}
\makeatother

\title{$c$-functions and Koornwinder polynomials}
\author{
Laura Colmenarejo\quad\ \ email:\ lcomen@ncsu.edu \\
Arun Ram\quad\quad\ \,\ email:\ aram@unimelb.edu.au \\
\\
}
\date{}

\usetikzlibrary{arrows.meta}

\begin{document}

\maketitle

\vspace{-3em}
\begin{center}
{\sl In memory of Ian G. Macdonald}
\end{center}

%\smallskip

\begin{abstract}
\noindent
This paper develops the theory of Macdonald-Koornwinder polynomials in parallel analogy with the work
done for the $GL_n$ case in~\cite{CR22}.   In the context of the type $CC_n$ affine root system the
Macdonald polynomials of other root systems of classical type are specializations of the Koornwinder polynomials.
We derive $c$-function formulas for symmetrizers and use them to give $E$-expansions, principal specializations 
and norm formulas for bosonic, mesonic and fermionic Koornwinder polynomials.  Finally, we explain the proof of the
norm conjectures and constant term conjectures for the Koornwinder case.
\end{abstract}

\keywords{Macdonald polynomials, symmetric functions, Hecke algebras}
\footnote{AMS Subject Classifications: Primary 05E05; Secondary  33D52.}

\setcounter{section}{-1}

\section{Introduction}

The Koornwinder polynomials are multivariate generalizations of the classical orthogonal polynomials that appear in the Askey scheme
\cite{K92}.  
At the top of this hierarchy we find the Askey-Wilson polynomials and the other families of orthogonal polynomials
in the Askey scheme are obtained from the Askey-Wilson polynomials by specializing parameters.

Macdonald's 1987 paper~\cite{Mac87} provides a very general framework for associating orthogonal polynomials 
$P_\lambda$ to any affine root system.
It turns out that the Koornwinder polynomials are the Macdonald polynomials for the affine root system of type $CC_n$ (in
Macdonald's notation $(C^\vee_n, C_n)$) and the Askey-Wilson polynomials are the Macdonald polynomials for the
affine root system of type $CC_1$.  One of the key features of Macdonald's picture is that the norms
$\langle P_\lambda, P_\lambda\rangle_+$ are generalizations of Macdonald's ``constant terms'' $\langle P_0,P_0 \rangle_+$.
In this way, Macdonald stated conjectures for the values of $\langle P_\lambda, P_\lambda\rangle_+$
which vastly expanded his earlier constant-term conjectures.

Cherednik introduced the double affine Hecke algebra as a tool for extending Opdam's ideas to prove the norm conjectures~\cite{C03}.
This perspective pointed to a larger family of orthogonal polynomials $E_\mu$, from which the $P_\lambda$
are obtained by a process of symmetrization.  All of these tools, including the proof of the norm conjectures, 
were wonderfully exposited in the full generality of a possibly non-reduced affine root system in Macdonald's book~\cite{Mac03}.
For a wonderful history of the exciting trajectory of these amazing developments see~\cite{HKO24}.

Particularly in the type $GL_n$ case, the Macdonald polynomials have been of interest to the combinatorial community because of
the wealth of wonderful $q,t$-generalizations of classical combinatorial formulas in symmetric function theory.  It is also
stimulating that there are many fascinating connections to adjacent fields (representations of $p$-adic groups and affine Lie algebras, geometry of Hilbert schemes and affine Springer fibers, torus knot invariants, vertex models in statistical mechanics, particle process in probability, etc).  For this reason it is desirable to provide expositions of the tools that bridge the language gaps between
the generality of affine root systems and the standard conventions in classical symmetric function theory.
In~\cite{CR22} we explained how, in the type $GL_n$ case, many of the combinatorial formulas can be understood from the theory of $c$-functions, which are the analogs of the Harish-Chandra $c$-functions surveyed by Helgason~\cite{Hel94} and that appear everywhere in  Macdonald's monograph~\cite{Mac03}.

In the type $GL_n$ case, the Macdonald polynomials depend on two parameters $q$ and $t$.  In the Koornwinder case
(type $CC_n$), the polynomials depend on 6 parameters, $q,t, t_0^{\frac12}, u_0^{\frac12}, t_n^{\frac12}$ and $u_n^{\frac12}$.  
Because of this plethora of parameters, sometimes it is not so easy to see how the combinatorial formulas familiar 
in the type $GL_n$ case generalize to the Koornwinder case.
This paper follows the same pattern as our earlier paper~\cite{CR22}, generalizing from the type $GL_n$ case
to the Koornwinder (type $CC_n$) case.

There is a constantly increasing literature on Koornwinder polynomials.  There are interesting
technical advances and also fascinating connections to other fields (see, for example,~\cite{CGdGW16, CMW23, Ra17, RW15, Yam20,YY21}).
The foundational work in~\cite{Nou95, Ra01, Sah99, Sto00, vD95, vD96}, among others, continues to be extremely useful
for clarifying the role and position of the double affine Hecke algebra as a tool for the Koornwinder case.
There are also important and very useful surveys of the theory of Koornwinder polynomials (see, for example,~\cite{Sto04, Sto21}).
There is a significant intersection between the content of this paper and the content presented in these surveys.
We hope that our combinatorial and $c$-function point of view can be useful in continuing research on Koornwinder polynomials.

\bigskip
The plan of this paper is as follows:

\bigskip

At the end of this introduction we include some remarks on the different sets of parameters used in the literature 
and establish the ones we will use for the paper.  Since there are 6 different parameters to keep track of
and lots of literature to navigate, perhaps this dictionary will be useful to readers (as it was for us).
Following these remarks we present a diagram of the affine root systems of classical types 
together with the specializations of the parameters
for obtaining the Macdonald polynomials of the corresponding type from Koornwinder polynomials, which are the Macdonald polynomials for the affine root system of
type $(C^\vee_n, C_n)$. A thorough study can be found in~\cite{YY21}.

Section~\ref{sec: affineWG} is dedicated to the affine Weyl group and root system for type $CC_n$, and we include the affine coroots and affine roots.  Our new contribution in Section 1 is the diagram giving the relative positions of the 
affine root systems of classical type and the specializations that give the Macdonald polynomials of 
the other classical type from Koornwinder polynomials.
Although these relationships are, in principle, well-known (from Macdonald~\cite{Mac72} and Bruhat-Tits~\cite{BT72}), we have not
seen this way of presenting this information, which we find very useful, and not broadly known.

In Section~\ref{sec: c-funct and DAHA} we introduce the main tools for working with Koornwinder polynomials, 
including the $c$-functions and the double affine Hecke algebra. 
Our contribution here is to provide a framework for the DAHA
in terms of $c$-functions, which makes the, sometimes daunting, formulas for the operators on the polynomial representation seem obvious and natural.

In the second half of Section~\ref{sec: c-funct and DAHA} we introduce four families of symmetrizers together 
with the relations between them, their $c$-function formulas and the case when the stabilizers are nontrivial. 
Our contribution here is to treat the four types of symmetrizers in tandem so that their role in the theory
(and the symmetry between them) becomes clearer.

Sections~\ref{sec: electronic MK} and~\ref{sec:BFM} examine the main objects of study:
\begin{enumerate}[itemsep=-0.2em, topsep=0.3em]
\item[(a)] electronic Macdonald-Koornwinder polynomials (Section~\ref{sec: electronic MK});
\item[(b)] bosonic Macdonald-Koornwinder polynomials;
\item[(c)] fermionic Macdonald-Koornwinder polynomials;
\item[(d)] mesonic Macdonald-Koornwinder polynomials.
\end{enumerate}
We introduce the electronic Macdonald-Koornwinder polynomials $E_\mu$ as eigenfunctions of Cherednik-Dunkl operators 
and then give a recursive formula and a creation formula for the $E_\mu$. For the other variants, our study includes:
\begin{enumerate}[itemsep=-0.2em, topsep=0.3em]
\item[(a)] definition of the  Weyl denominators;
\item[(b)] study of the bosonic, fermionic and mesonic spaces;
\item[(c)] formulas for the Poincar\'e polynomial; 
\item[(d)] expansions in terms of $E_\mu$, and
\item[(e)] principal specializations.
\end{enumerate}
Our contribution here is to put the focus on the fermionic and mesonic Koornwinder polynomials so that the four-fold structure
is clearly visible.
This four-fold structure eventually leads to powerful recursions for computing norms.

Finally, Section~\ref{sec:orth} is dedicated to the study of the Macdonald-Koornwinder polynomials as a family of orthogonal polynomials. In particular, we 
\begin{enumerate}[itemsep=-0.2em, topsep=0.3em]
\item[(a)] define the Macdonald-Koornwinder inner product via multiplication by a kernel and taking the constant term
and characterize the electronic and bosonic Macdonald-Koornwinder polynomials in terms of the inner product;
\item[(b)] compute adjoints of the operators from the double affine Hecke algebra;
\item[(c)] prove the going up a level and Weyl character formulas to provide recursions for norms; and
\item[(d)] use the recursions for norms to compute the norms $\langle P_\lambda, P_\lambda\rangle_+$ 
and the Macdonald constant term for type $CC_n$.
\end{enumerate}
This section follows the same trajectory as that taken in~\cite[Ch.\ 5]{Mac03}.
Our contribution here is to use the fermionic and mesonic framework to organize the recursions for norms and make
the proof of norm conjectures easy and natural.

\subsection*{Acknowledgments}

We thank Zajj Daugherty, Stephen Griffeth, Bogdan Ion, Alejandro H. Morales and Martha Yip for discussions which greatly helped to
solidify our understanding of the topics that we treat in this paper.
We thank O.\ Warnaar and J.\ Stokman for helpful and insightful comments on a preliminary version of this paper.

L.C.\  was partially supported by a grant from the Simons Foundation (MPS-TSM-00007362).

 \begin{landscape}

\subsubsection{The poset of affine root systems of classical type}
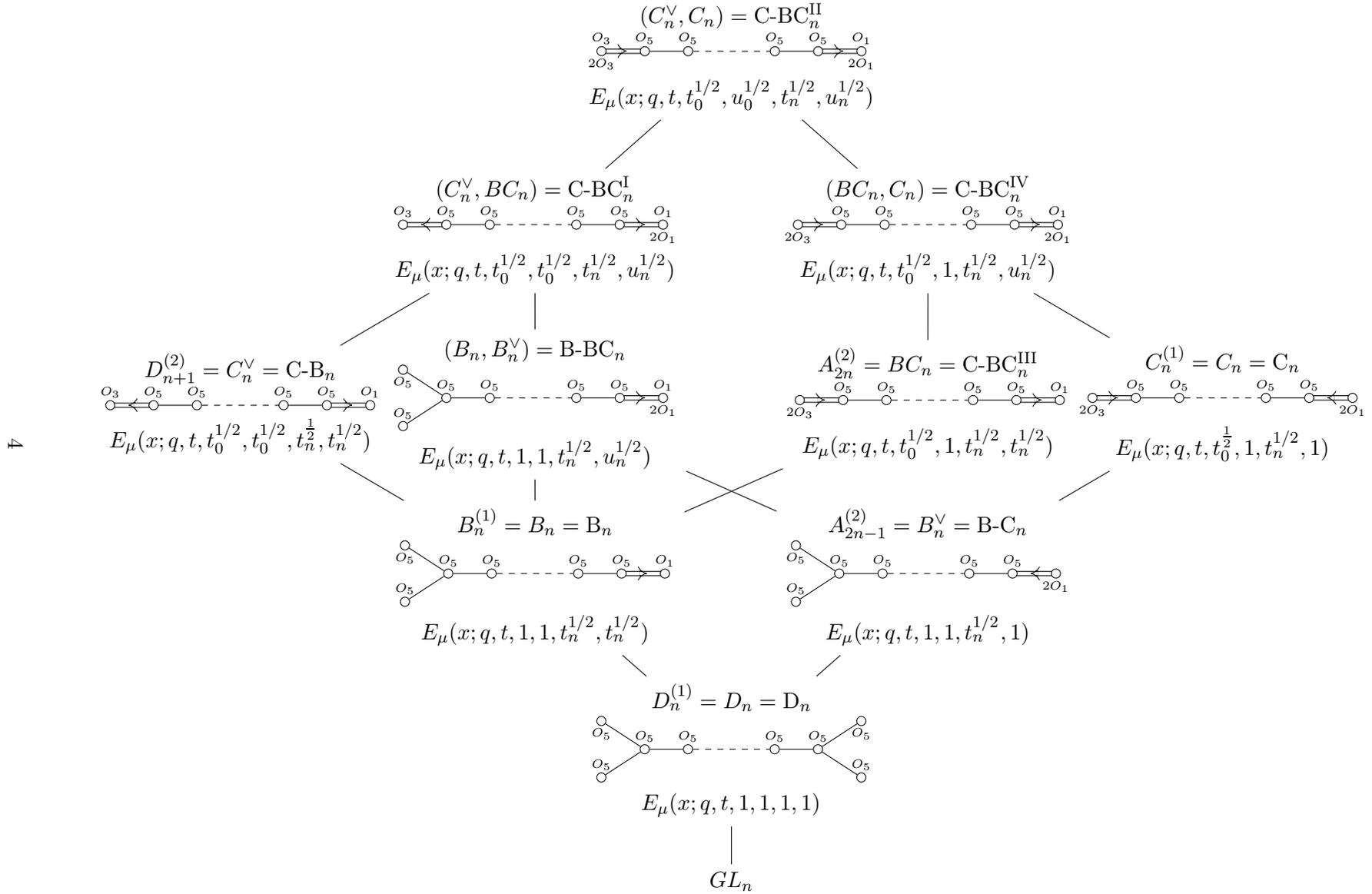
\label{affrootsystemsfig}

$$\begin{tikzpicture}[xscale=1.7,yscale=-1.5]
\node (CvC) at (0,0) 
{$\begin{matrix}
(C^\vee_n,C_n)=\hbox{C-BC}_n^{\mathrm{II}}
\\
\begin{matrix}
\begin{tikzpicture}[every node/.style={inner sep=1}, scale=.75]
	\node at (.5,0) {\tikz \draw[black, ->, >=implies, double distance=2] (0,0) -- +(.127,0);};
	\node at (5.5,0) {\tikz \draw[black, ->, >=implies, double distance=2] (0,0) -- +(.127,0);};
	\foreach \x in {0,5}{\foreach \y in {-1,1}{\draw (\x,\y*.05) to (\x+1,\y*.05);}}
	\node[wV, label=below:{\tiny$2O_3$}, label=above:{\tiny$O_3$}] (0) at (0,0) {};
	\node[wV, label=below:{\tiny$2O_1$}, label=above:{\tiny$O_1$}] (6) at (6,0) {};
	\foreach \x/\y in {1/1,2/2,4/n\!-\!2,5/n\!-\!1}{
		\node[wV,label=above:{\tiny$O_5$}] (\x) at (\x,0) {};}
	\draw (1) to (2) (4) to (5);
	\draw[dashed] (2) to (4);
	\end{tikzpicture}\end{matrix}
\\
E_\mu(x;q,t, t_0^{1/2},u_0^{1/2}, t_n^{1/2},u_n^{1/2})
\end{matrix}
$};
\node (CvBC) at (-2,2) {
$\begin{matrix}
(C^\vee_n,BC_n)=\hbox{C-BC}_n^{\mathrm{I}}
\\
\begin{matrix}
\begin{tikzpicture}[every node/.style={inner sep=1}, scale=.75]
	\node at (.5,0) {\tikz \draw[black, ->, >=implies, double distance=2] (0,0) -- +(-.127,0);};
	\node at (5.5,0) {\tikz \draw[black, ->, >=implies, double distance=2] (0,0) -- +(.127,0);};
	\foreach \x in {0,5}{\foreach \y in {-1,1}{\draw (\x,\y*.05) to (\x+1,\y*.05);}}
	\node[wV, label=above:{\tiny$O_3$}] (0) at (0,0) {};
	\node[wV, label=below:{\tiny$2O_1$}, label=above:{\tiny$O_1$}] (6) at (6,0) {};
	\foreach \x/\y in {1/1,2/2,4/n\!-\!2,5/n\!-\!1}{
		\node[wV,label=above:{\tiny$O_5$}] (\x) at (\x,0) {};}
	\draw (1) to (2) (4) to (5);
	\draw[dashed] (2) to (4);
	\end{tikzpicture}\end{matrix}
\\
E_\mu(x;q,t, t_0^{1/2},t_0^{1/2},t_n^{1/2},u_n^{1/2})
\end{matrix}
$};
\node (BCC) at (2,2) {
$\begin{matrix}
(BC_n,C_n)=\hbox{C-BC}_n^{\mathrm{IV}}
\\
\begin{matrix}
\begin{tikzpicture}[every node/.style={inner sep=1}, scale=.75]
	\node at (.5,0) {\tikz \draw[black, ->, >=implies, double distance=2] (0,0) -- +(.127,0);};
	\node at (5.5,0) {\tikz \draw[black, ->, >=implies, double distance=2] (0,0) -- +(.127,0);};
	\foreach \x in {0,5}{\foreach \y in {-1,1}{\draw (\x,\y*.05) to (\x+1,\y*.05);}}
	\node[wV, label=below:{\tiny$2O_3$}] (0) at (0,0) {};
	\node[wV, label=below:{\tiny$2O_1$}, label=above:{\tiny$O_1$}] (6) at (6,0) {};
	\foreach \x/\y in {1/1,2/2,4/n\!-\!2,5/n\!-\!1}{
		\node[wV,label=above:{\tiny$O_5$}] (\x) at (\x,0) {};}
	\draw (1) to (2) (4) to (5);
	\draw[dashed] (2) to (4);
	\end{tikzpicture}\end{matrix}
\\
E_\mu(x;q,t, t_0^{1/2},1,t_n^{1/2},u_n^{1/2})
\end{matrix}
$};
\node (Cv) at (-5,4) {
$\begin{matrix}
D_{n+1}^{(2)} = C_n^\vee = \hbox{C-B}_n
\\
\begin{matrix}
\begin{tikzpicture}[every node/.style={inner sep=1}, scale=.75]
	\node at (.5,0) {\tikz \draw[black, ->, >=implies, double distance=2] (0,0) -- +(-.127,0);};
	\node at (5.5,0) {\tikz \draw[black, ->, >=implies, double distance=2] (0,0) -- +(.127,0);};
	\foreach \x in {0,5}{\foreach \y in {-1,1}{\draw (\x,\y*.05) to (\x+1,\y*.05);}}
	\node[wV, label=above:{\tiny$O_3$}] (0) at (0,0) {};
	\node[wV, label=above:{\tiny$O_1$}] (6) at (6,0) {};
	\foreach \x/\y in {1/1,2/2,4/n\!-\!2,5/n\!-\!1}{
		\node[wV,label=above:{\tiny$O_5$}] (\x) at (\x,0) {};}
	\draw (1) to (2) (4) to (5);
	\draw[dashed] (2) to (4);
	\end{tikzpicture}\end{matrix}
\\
E_\mu(x;q,t,t_0^{1/2},t_0^{1/2},t_n^{\frac12},t_n^{1/2})
\end{matrix}
$};
\node (BBv) at (-2,4) {
$\begin{matrix}
(B_n,B_n^\vee)= \hbox{B-BC}_n
\\
	\begin{matrix}\begin{tikzpicture}[every node/.style={inner sep=1}, scale=.75]
	\node at (6.5,0) {\tikz \draw[black, ->, >=implies, double distance=2] (0,0) -- +(.127,0);};
		\foreach \y in {-1,1}{\draw (6,\y*.05) to (6+1,\y*.05);}
		\node[wV, label=below:{\tiny$2O_1$}, label=above:{\tiny$O_1$}] (7) at (7,0) {};
	\foreach \x/\y in {2/2,3/3,5/n\!-\!2,6/n\!-\!1}{
		\node[wV, label=above:{\tiny$O_5$}] (\x) at (\x,0) {};}
	\node[wV, label=above:{\tiny$O_5$}] (0) at (1,-.65) {};
	\node[wV, label=below:{\tiny$O_5$}] (1) at (1,.65) {};
	\draw (1) to (2) (0) to (2) (2) to (3) (5) to (6);
	\draw[dashed] (3) to (5);
	\end{tikzpicture}\end{matrix}
\\
E_\mu(x;q,t,1,1,t_n^{1/2},u_n^{1/2})
\end{matrix}
$};
\node (BC) at (2,4) {
$\begin{matrix}
A_{2n}^{(2)}=BC_n=\hbox{C-BC}_n^{\mathrm{III}}
\\
\begin{matrix}
\begin{tikzpicture}[every node/.style={inner sep=1}, scale=.75]
	\node at (.5,0) {\tikz \draw[black, ->, >=implies, double distance=2] (0,0) -- +(.127,0);};
	\node at (5.5,0) {\tikz \draw[black, ->, >=implies, double distance=2] (0,0) -- +(.127,0);};
	\foreach \x in {0,5}{\foreach \y in {-1,1}{\draw (\x,\y*.05) to (\x+1,\y*.05);}}
	\node[wV, label=below:{\tiny$2O_3$}] (0) at (0,0) {};
	\node[wV, label=above:{\tiny$O_1$}] (6) at (6,0) {};
	\foreach \x/\y in {1/1,2/2,4/n\!-\!2,5/n\!-\!1}{
		\node[wV,label=above:{\tiny$O_5$}] (\x) at (\x,0) {};}
	\draw (1) to (2) (4) to (5);
	\draw[dashed] (2) to (4);
	\end{tikzpicture}\end{matrix}
\\
E_\mu(x;q,t,t_0^{1/2},1,t_n^{1/2},t_n^{1/2})
\end{matrix}
$};
\node (C) at (5,4) {
$\begin{matrix}
C_n^{(1)}=C_n=\hbox{C}_n
\\
\begin{matrix}
\begin{tikzpicture}[every node/.style={inner sep=1}, scale=.75]
	\node at (.5,0) {\tikz \draw[black, ->, >=implies, double distance=2] (0,0) -- +(.127,0);};
	\node at (5.5,0) {\tikz \draw[black, ->, >=implies, double distance=2] (0,0) -- +(-.127,0);};
	\foreach \x in {0,5}{\foreach \y in {-1,1}{\draw (\x,\y*.05) to (\x+1,\y*.05);}}
	\node[wV, label=below:{\tiny$2O_3$}] (0) at (0,0) {};
	\node[wV, label=below:{\tiny$2O_1$}] (6) at (6,0) {};
	\foreach \x/\y in {1/1,2/2,4/n\!-\!2,5/n\!-\!1}{
		\node[wV,label=above:{\tiny$O_5$}] (\x) at (\x,0) {};}
	\draw (1) to (2) (4) to (5);
	\draw[dashed] (2) to (4);
	\end{tikzpicture}\end{matrix}
\\
E_\mu(x;q,t,t_0^{\frac12},1,t_n^{1/2},1)
\end{matrix}
$};
\node (B) at (-2,6) {
$\begin{matrix}
B_n^{(1)} = B_n = \hbox{B}_n
\\
	\begin{matrix}\begin{tikzpicture}[every node/.style={inner sep=1}, scale=.75]
	\node at (6.5,0) {\tikz \draw[black, ->, >=implies, double distance=2] (0,0) -- +(.127,0);};
		\foreach \y in {-1,1}{\draw (6,\y*.05) to (6+1,\y*.05);}
		\node[wV, label=above:{\tiny$O_1$}] (7) at (7,0) {};
	\foreach \x/\y in {2/2,3/3,5/n\!-\!2,6/n\!-\!1}{
		\node[wV, label=above:{\tiny$O_5$}] (\x) at (\x,0) {};}
	\node[wV, label=above:{\tiny$O_5$}] (0) at (1,-.65) {};
	\node[wV, label=below:{\tiny$O_5$}] (1) at (1,.65) {};
	\draw (1) to (2) (0) to (2) (2) to (3) (5) to (6);
	\draw[dashed] (3) to (5);
	\end{tikzpicture}\end{matrix}
\\
E_\mu(x;q,t,1,1,t_n^{1/2},t_n^{1/2})
\end{matrix}
$};
\node (Bv) at (2,6) {
$\begin{matrix}
A_{2n-1}^{(2)} = B_n^\vee =  \hbox{B-C}_n
\\
	\begin{matrix}\begin{tikzpicture}[every node/.style={inner sep=1}, scale=.75]
	\node at (6.5,0) {\tikz \draw[black, ->, >=implies, double distance=2] (0,0) -- +(-.127,0);};
		\foreach \y in {-1,1}{\draw (6,\y*.05) to (6+1,\y*.05);}
		\node[wV, label=below:{\tiny$2O_1$}] (7) at (7,0) {};
	\foreach \x/\y in {2/2,3/3,5/n\!-\!2,6/n\!-\!1}{
		\node[wV, label=above:{\tiny$O_5$}] (\x) at (\x,0) {};}
	\node[wV, label=above:{\tiny$O_5$}] (0) at (1,-.65) {};
	\node[wV, label=below:{\tiny$O_5$}] (1) at (1,.65) {};
	\draw (1) to (2) (0) to (2) (2) to (3) (5) to (6);
	\draw[dashed] (3) to (5);
	\end{tikzpicture}\end{matrix}
\\
E_\mu(x;q,t,1,1,t_n^{1/2},1)
\end{matrix}
$};
\node (D) at (0,8) {
$\begin{matrix}
D_n^{(1)}=D_n=\hbox{D}_n
\\
	\begin{matrix}\begin{tikzpicture}[every node/.style={inner sep=1}, scale=.75]
	\foreach \x/\y in {2/2,3/3,5/n\!-\!3,6/n\!-\!2\phantom{-} }{
		\node[wV, label=above:{\tiny$O_5$}] (\x) at (\x,0) {};}
	\node[wV, label=above:{\tiny$O_5$}] (0) at (1,-.65) {};
	\node[wV, label=below:{\tiny$O_5$}] (1) at (1,.65) {};
	\node[wV, label=above:{\tiny$O_5$}] (7) at (7,-.65) {};
	\node[wV, label=below:{\tiny$O_5$}] (8) at (7,.65) {};
	\draw (1) to (2) (0) to (2) (2) to (3) (5) to (6) (6) to (7) (6) to (8);
	\draw[dashed] (3) to (5);
	\end{tikzpicture}\end{matrix}
\\
E_\mu(x;q,t,1,1,1,1)
\end{matrix}
$};

\node (A) at (0,9.5) {$GL_n$};
\draw (CvC) -- (CvBC);
\draw (CvC) -- (BCC);
\draw (CvBC) -- (Cv);
\draw (CvBC) -- (BBv);
\draw (BCC) -- (BC);
\draw (BCC) -- (C);
\draw (Cv) -- (B);
\draw (BBv) -- (B);
\draw (BBv) -- (Bv);
\draw (BC) -- (B);
\draw (C) -- (Bv);
\draw (B) -- (D);
\draw (Bv) -- (D);
\draw (D) -- (A);
\end{tikzpicture}
$$
See Section~\ref{ssec:paramspec} for the explanation of this diagram.
\end{landscape}

\subsection{Remarks on parameters}

Depending on the reference, the notation for the parameters varies. 
In this article, we follow~\cite{Nou95} and~\cite{Sah99} and mostly use the parameters
$$
q, t^{\frac12}, t_0^{\frac12}, u_0^{\frac12}, t_n^{\frac12}, u_n^{\frac12}. 
$$

\begin{remark} \label{rem:paramconv}
In an attempt to relate the parameter notations in~\cite[\S 4.7]{Mac03},~\cite[\S 3]{Sah99} and~\cite[Def.\ 2.1]{C03}
let
$$\tau_0' = u_0^{\frac12}, \quad
\tau_n' = u_n^{\frac12}, \quad
\tau_0  = t_0^{\frac12},\quad
\tau_n = t_n^{\frac12}, \quad\hbox{and}\quad
\tau_i = \tau_i' = t^{\frac12}\ \ \hbox{for $i\in \{1,\ldots, n-1\}$.}
$$
The summary of (1.5.1), (4.4.1), (4.4.2), (4.4.3), and (5.1.4) in~\cite{Mac03} is that, for an affine root $a$,
\begin{equation*}
\tau_a = (t_at_{2a})^{\frac12} = q^{\frac12\kappa_a} = q^{\frac12(k(a)+k(2a))},
\qquad\hbox{and} \qquad
\tau_a' = t_a^\frac12 = q^{\frac12\kappa_a'} = q^{\frac12(k(a)-k(2a))}.
\end{equation*}
In our situation
\begin{align*}
t_n^{\frac12} &= \tau_n = q^{\frac12\kappa_n}  
= t_{\varepsilon_n}^{\frac12}t_{2\varepsilon_n}^{\frac12}
= q^{\frac12 k(\varepsilon_n) + \frac12 k(2\varepsilon_n)} = q^{\frac12 k_1 + \frac12 k_2},  \\
u_n^{\frac12} &= \tau_n' =  q^{\frac12\kappa_n'} = t_{\varepsilon_n}^{\frac12}
= q^{\frac12 k(\varepsilon_n) - \frac12 k(2\varepsilon_n)} = q^{\frac12 k_1 - \frac12 k_2}, \\
t_0^{\frac12} &= \tau_0 =  q^{\frac12\kappa_0} 
= t_{-\varepsilon_1+\frac12 \delta}^{\frac12}
t_{-2\varepsilon_1+\delta}^{\frac12}
= q^{\frac12 k(-\varepsilon_1+\frac12\delta) + \frac12 k(-2\varepsilon_1+\delta)} 
= q^{\frac12 k_3 + \frac12 k_4}, \\
u_0^{\frac12} &= \tau'_0 =  q^{\frac12\kappa_0'} 
= t_{-\varepsilon_1+\frac12\delta}^{\frac12}
= q^{\frac12 k(-\varepsilon_1+\frac12\delta) - \frac12 k(-2\varepsilon_1+\delta)} = q^{\frac12 k_3 - \frac12 k_4},
\qquad\hbox{and} \\
t^{\frac12} &= \tau_i  =  q^{\frac12\kappa}
= t_{\varepsilon_i-\varepsilon_{i+1}}^{\frac12} = q^{\frac12k_5},
\quad\hbox{for $i\in \{1,\ldots, n-1\}$,}
\end{align*}
and the formulas in~\cite[(1.5.1)]{Mac03} correspond to interchanging $\kappa_0$ and
$\kappa_n'$.
\qed
\end{remark}

\begin{remark} \textbf{Askey-Wilson parameters.}
In type $(C_1^\vee, C_1)$, the bosonic Macdonald polynomials $P_\lambda(q, t_1,u_1,t_0,u_0)$
are also known as the Askey-Wilson polynomials. Following~\cite[\S3]{Nou95},
the correspondence to the original Askey-Wilson parameters is given by
\begin{equation}
q=q, \quad
a = q^{\frac12}t_0^{\frac12}u_0^{\frac12}, \quad
b = - q^{\frac12}t_0^{\frac12}u_0^{-\frac12}, \quad
c = t_n^{\frac12}u_n^{\frac12}, \quad
d = -t_n^{\frac12}u_n^{-\frac12}.
\label{AWparams}
%\tag{AWparams}
\end{equation} 
These conversions are equivalent to
$$t_0 = -q^{-1}ab, \qquad t_n = -cd, \qquad u_0 = -ab^{-1}, \qquad u_n = -cd^{-1},$$
and it is useful to note that
$$a+b = q^{\frac12}t_0^{\frac12}(u_0^{\frac12}-u_0^{-\frac12})
\qquad\hbox{and}\qquad
c+d = t_n^{\frac12}(u_n^{\frac12}-u_n^{-\frac12}).
$$
Up to permutations of $a,b,c,d$, these parameters are used 
in~\cite[(1)]{Sah00}, 
in~\cite[(17)]{CGdGW16}, 
in~\cite[Def. 2.2]{CMW23} 
and, with different notation, in ~\cite[(5.1.14)]{Mac03}.  
\qed
\end{remark}

\newpage

\section{The affine Weyl group and root system for type $CC_n$}\label{sec: affineWG}

The affine root system of type $CC_n$ (in Macdonald's notation $(C^\vee_n, C_n)$) is the structure that
holds the combinatorics of Koornwinder polynomials in place as they are the Macdonald polynomials for this affine root system.  
The affine Weyl group $W$ 
plays the role of 
the group of symmetries of the affine root system.  In this section we introduce the definitions and notations
for working with the affine Weyl group $W$ and the affine root system of type $CC_n$.

The coroots $S^\vee$ and the roots $S$ for the affine root system of type $CC_n$ play just slightly different roles in the
theory, especially in the computations involving Koornwinder polynomials. One of the challenges in this work is to keep these
two mirror worlds in proper focus.  For this purpose, in Section~\ref{subsection: W_X and W_Y}, we carefully lay out two versions of the affine Weyl group $W$,
one denoted $W_X$ which acts on the coroots (with 5 orbits), and one denoted $W_Y$ which acts on the roots (with 5 orbits).
While the groups $W_X$, $W_Y$ and $W$ are all isomorphic, being pedantic about the notation at this
early stage prevents future headaches.  

To conclude this section we present a brief explanation of the reasoning for how the Macdonald polynomials
of other classical types (such as $B_n$, $C_n$, $BC_n$, $A_{2n-1}^{(2)}$, $D_{n+1}^{(2)}$ etc.) 
are derived from the Koornwinder polynomials by specialization.
This specialization process is summarized in the diagram included in Section~\ref{affrootsystemsfig}.

\subsection{The affine Weyl group $W$ and the finite Weyl group $W_{\mathrm{fin}}$}

Use a graphical notation for relations so that
$$
\begin{array}{cl}
\begin{matrix}
\begin{tikzpicture}[every node/.style={inner sep=1}, scale=1]
	\node[wV, label=above:{$g_i$} ] (0) at (0,0) {};
	\node[wV, label=above:{$g_j$} ] (1) at (1,0) {};
	\end{tikzpicture}\end{matrix}
    &\hbox{means $g_ig_j =g_jg_i$,} 
\\[0.4cm]
\begin{matrix} \begin{tikzpicture}[every node/.style={inner sep=1}, scale=1]
	\foreach \x in {0}{\foreach \y in {0}{\draw (\x,\y*.05) to (\x+1,\y*.05);}}
	\node[wV, label=above:{$g_i$} ] (0) at (0,0) {};
	\node[wV, label=above:{$g_j$} ] (1) at (1,0) {};
	\end{tikzpicture}\end{matrix}
&\hbox{means $g_ig_jg_i = g_jg_ig_j$, and} 
\\[0.4cm]
\begin{matrix} \begin{tikzpicture}[every node/.style={inner sep=1}, scale=1]
	\node at (.5,0) {\tikz \draw[black, double distance=2] (0,0) -- +(-.127,0);};
	\foreach \x in {0}{\foreach \y in {-1,1}{\draw (\x,\y*.05) to (\x+1,\y*.05);}}
	\node[wV, label=above:{$g_i$} ] (0) at (0,0) {};
	\node[wV, label=above:{$g_j$} ] (1) at (1,0) {};
	\end{tikzpicture}\end{matrix}
    &\hbox{means $g_ig_jg_ig_j = g_jg_ig_jg_i$.} 
\end{array}
$$

The \emph{affine Weyl group} is the group  $W$ presented by generators $s_0,s_1, \ldots, s_n$ and relations
$$s_i^2 = 1
\qquad\hbox{and}\qquad
	\begin{matrix}\begin{tikzpicture}[every node/.style={inner sep=1}, scale=1]
	\node at (.5,0) {\tikz \draw[black, double distance=2] (0,0) -- +(.127,0);};
	\node at (5.5,0) {\tikz \draw[black, double distance=2] (0,0) -- +(-.127,0);};
	\foreach \x in {0,5}{\foreach \y in {-1,1}{\draw (\x,\y*.05) to (\x+1,\y*.05);}}
	\node[wV, label=above:{$s_0$}] (0) at (0,0) {};
	\node[wV, label=above:{$s_n$}] (6) at (6,0) {};
	\foreach \x/\y in {1/1,2/2,4/n\!-\!2,5/n\!-\!1}{
		\node[wV] (\x) at (\x,0) {};}
	\draw (1) to (2) (4) to (5);
	\draw[dashed] (2) to (4);
	\node[wV, label=above:{$s_1$}] (1) at (1,0) {};
	\node[wV, label=above:{$s_2$}] (2) at (2,0) {};
	\node[wV, label=above:{$s_{n-2}$}] (4) at (4,0) {};
	\node[wV, label=above:{$s_{n-1}$}] (5) at (5,0) {};
	\end{tikzpicture}\end{matrix}
$$

The \emph{finite Weyl group} is the subgroup $W_{\mathrm{fin}}$ generated by $s_1, \ldots, s_n$.

Let $w\in W$.  The \emph{length of $w$}, $\ell(w)$, is the minimal $\ell\in \ZZ_{\ge 0}$ such that
$$w= s_{i_1}\ldots s_{i_\ell}
\quad\hbox{with $i_1, \ldots, i_\ell\in \{0,1, \ldots, n-1, n\}$}.$$
The expression $w= s_{i_1}\ldots s_{i_\ell}$ is a \emph{reduced word for $w$} and
any other expression of the form $w= s_{j_1}\ldots s_{j_k}$, with $j_1, \ldots, j_k\in \{0,1, \ldots, n-1, n\}$, has $k\geq \ell(w)$.

\subsubsection{Translation presentation of $W$}

Define $h_{\varepsilon_1}, \ldots, h_{\varepsilon_n}\in W$ by
$$h_{\varepsilon_1} = s_0s_1s_2\cdots s_{n-1}s_n s_{n-1}\cdots s_2s_1,
\qquad\hbox{and}\qquad
h_{\varepsilon_j} = s_j h_{\varepsilon_{j-1}}s_j,\ \hbox{for $j\in \{2, \ldots, n\}$.}
$$
For $\mu = (\mu_1, \ldots, \mu_n)\in \ZZ^n$ define the \emph{translation} $h_\mu$ by
\begin{equation}
h_\mu = (h_{\varepsilon_1})^{\mu_1}\cdots (h_{\varepsilon_n})^{\mu_n}
\label{tmudefn}
%\tag{tmudefn}
\end{equation}
and define $u_\mu\in W$ and $v_\mu\in W_{\mathrm{fin}}$ by the equation
\begin{equation}
h_\mu = u_\mu v_\mu,
\qquad\hbox{where $v_\mu\in W_{\mathrm{fin}}$ and 
$u_\mu$ is minimal length in the coset $h_\mu W_{\mathrm{fin}}$.}
\label{umuvmudefn}
%\tag{umuvmudefn}
\end{equation}

Define an action of $W_{\mathrm{fin}}$ on $\ZZ^n$ by 
\begin{align}
s_i (\mu_1, \ldots, \mu_n) &= (\mu_1, \ldots, \mu_{i-1}, \mu_{i+1}, \mu_i, \mu_{i+2}, \ldots, \mu_n),
\qquad\hbox{for $i\in \{1, \ldots, n-1\}$,}
\nonumber
\\
s_n (\mu_1, \ldots, \mu_n) &= (\mu_1, \ldots, \mu_{n-1}, -\mu_n).
\label{Wfinaction}
%\tag{Wfinaction}
\end{align}
Then 
$$W = \ZZ^n \ltimes W_{\mathrm{fin}}.$$
In other words,  if $\mu, \nu\in \ZZ^n$ and $v\in W_{\mathrm{fin}}$ then
$$vh_\mu = h_{v\mu} v, \qquad h_\mu h_\nu = h_{\mu+\nu}
\qquad\hbox{and}\qquad W = \{ h_\mu v\ |\ \mu\in \ZZ^n, v\in W_{\mathrm{fin}}\}.
$$

\begin{remark}
   In~\cite{CR22}, we use $t_\mu$ to denote the translations in the type $A$ case. We use the notation $h_\mu$ 
   for type $CC_n$ to avoid conflict with the set of parameters for Koornwinder polynomials (specified in Section~\ref{ssec:params}).
\end{remark}

\subsection{Affine coroots, affine roots and the groups $W_X$ and $W_Y$}\label{subsection: W_X and W_Y}

In this subsection we set up the notation for the affine Weyl groups $W_X$ and $W_Y$.
Both groups $W_X$ and $W_Y$ are isomorphic to $W$, 
but they serve slightly different roles and it is necessary to set up the notation to distinguish them.

\subsubsection{The dual lattices $\fa_\ZZ$ and $\fa_\ZZ^*$}\label{aZaZ*}

Let $\varepsilon_1, \ldots, \varepsilon_n$ and $\varepsilon_1^\vee, \ldots, \varepsilon_n^\vee$ be symbols and define
dual lattices (i.e. dual free $\ZZ$-modules)
$$\fa_\ZZ^* = \{ \gamma_1\varepsilon_1+ \cdots + \gamma_n \varepsilon_n\ |\ \gamma_1,\ldots, \gamma_n\in \ZZ\}
\qquad\hbox{and}\qquad
\fa_\ZZ = \{ \mu_1\varepsilon^\vee_1+ \cdots + \mu_n \varepsilon^\vee_n\ |\ \mu_1,\ldots, \mu_n\in \ZZ\}$$
with $\ZZ$-bilinear pairing
$$\langle \ , \ \rangle\colon \fa^*_\ZZ\times \fa_\ZZ \to \ZZ
\qquad\hbox{given by}\qquad
\langle \varepsilon_i, \varepsilon_j^\vee\rangle = \delta_{ij}.
$$
Both $\fa_\ZZ$ and $\fa_\ZZ^*$ are isomorphic to $\ZZ^n$.

\subsubsection{The affine coroots for type $CC_n$}

Let
$$Q^\vee = \hbox{$\ZZ$-span}\{ \varepsilon_1^\vee, \ldots, \varepsilon_n^\vee, \hbox{$\frac12$} K\}
$$
be the $\ZZ$-vector space spanned by symbols $\varepsilon_1^\vee, \ldots, \varepsilon_n^\vee$ and $\frac12 K$.
The affine Weyl group $W_X$ is the group of $\ZZ$-linear transformations of $Q^\vee$
generated by the transformations $s_0^\vee, s_1^\vee, \ldots, s_n^\vee$ given as follows:
If $\lambda^\vee = \lambda_1 \varepsilon^\vee_1 + \cdots + \lambda_n \varepsilon^\vee_n + \frac{k}{2}K$ then
\begin{align}
s_0^\vee \lambda^\vee &= - \lambda_1\varepsilon^\vee_1 + \lambda_2\varepsilon^\vee_2
+\cdots+\lambda_n \varepsilon^\vee_n + \left(\hbox{$\frac{k}{2}$}+\lambda_1\right) K,
\nonumber  \\
s_n^\vee \lambda^\vee &= \lambda_1\varepsilon^\vee_1+\cdots +\lambda_{n-1}\varepsilon^\vee_{n-1}
-\lambda_n \varepsilon^\vee_n + \hbox{$\frac{k}{2}$}K, \qquad\hbox{and}
\label{lvl0onaZ}
%\tag{lvl0onaZ}
\\
s_i^\vee \lambda^\vee &= \lambda_1\varepsilon^\vee_1 + \cdots + \lambda_{i-1}\varepsilon^\vee_{i-1}
+\lambda_{i+1}\varepsilon^\vee_i+\lambda_i \varepsilon_{i+1}^\vee + \lambda_{i+2}\varepsilon^\vee_{i+2}+\cdots +
\lambda_n\varepsilon^\vee_n + \hbox{$\frac{k}{2}$} K,
\nonumber
\end{align}
for $i\in \{1, \ldots, n-1\}$.  
Let
$$s^\vee_{\varepsilon_1} = s^\vee_1\cdots s^\vee_n\cdots s^\vee_1, \qquad
h_{\varepsilon_1} = s_0^\vee  s^\vee_{\varepsilon_1} = s_0^\vee s^\vee_1\cdots s^\vee_n\cdots s^\vee_1
\qquad\hbox{and}\qquad
h_{\varepsilon_{i+1}} = s_i^\vee h_{\varepsilon_i}s_i^\vee,
$$
for $i\in \{1, \ldots, n-1\}$.  Then
$h_{\varepsilon_1}\lambda^\vee 
= s_0^\vee s_{\varepsilon_1}^\vee\lambda^\vee
= \lambda_1\varepsilon^\vee_1
+\cdots+\lambda_n \varepsilon^\vee_n + \left(\hbox{$\frac{k}{2}$}-\lambda_1\right) K$ and
$$
h_{\varepsilon_i}\lambda^\vee 
= \lambda_1\varepsilon^\vee_1
+\cdots+\lambda_n \varepsilon^\vee_n + \left(\hbox{$\frac{k}{2}$}-\lambda_i\right) K,
\qquad\hbox{for $i\in \{1, \ldots, n\}$.}
$$
If $\gamma = \gamma_1\varepsilon_1+\cdots \gamma_n\varepsilon_n$ and 
$\lambda^\vee = \lambda_1 \varepsilon^\vee_1 + \cdots + \lambda_n \varepsilon^\vee_n + \frac{k}{2}K$ then
\begin{align*}
h_\gamma \lambda^\vee &= h_{\varepsilon_1}^{\gamma_1}\cdots h_{\varepsilon_n}^{\gamma_n} \lambda^\vee
=\lambda^\vee + (\hbox{$\frac{k}{2}$}-(\lambda_1\gamma_1+\cdots+\lambda_n\gamma_n))K 
=\lambda^\vee + (\hbox{$\frac{k}{2}$}-\langle \gamma,\lambda^\vee\rangle )K,
\end{align*}
and special cases of this last formula are
\begin{equation}
h_\gamma \varepsilon_i^\vee 
=\varepsilon_i^\vee -\gamma_i K
\qquad\hbox{and}\qquad
h_\gamma K = K.
\label{lvl0trWY}
%\tag{lvl0trWY}
\end{equation}

The set of coroots $S^\vee$ for type $CC_n$ is the union of the five $W_X$-orbits given by
\begin{align*}
O^\vee_1 
&= W_X\cdot\alpha_n^\vee = W_X\cdot \varepsilon_n^\vee 
= \{ \pm \varepsilon_i^\vee+rK \ |\  i\in \{1, \ldots, n\}, r\in \ZZ\}, \\
2O^\vee_1 &= W_X\cdot 2\alpha^\vee_n = W_X\cdot 2\varepsilon_n^\vee
= \{ \pm 2\varepsilon_i^\vee+2rK \ |\  i\in \{1, \ldots, n\}, r\in \ZZ\} , \\
O^\vee_3 
&= W_X\cdot \alpha_0^\vee = W_X\cdot (-\varepsilon_1^\vee+\hbox{$\frac12$}K) 
= \{ \pm(\varepsilon_i^\vee+\hbox{$\frac12$}(2r+1)K \ |\  i\in \{1, \ldots, n\}, r\in \ZZ\}, \\
2O^\vee_3 
&= W_X\cdot 2\alpha_0^\vee = W_Y\cdot (-2\varepsilon_1^\vee+K)
= \{ \pm 2\varepsilon_i^\vee+(2r+1)K \ |\  i\in \{1, \ldots, n\}, r\in \ZZ\} , \\
O^\vee_5  
&= W_X\cdot \alpha^\vee_1 = W_X\cdot (\varepsilon_1^\vee-\varepsilon_2^\vee)
= 
\left\{ \begin{array}{l}
\pm(\varepsilon_i^\vee+\varepsilon_j^\vee)+rK,  \\
\pm(\varepsilon_i^\vee-\varepsilon_j^\vee)+rK  
\end{array}\ \Big\vert\ 
i,j\in \{1, \ldots, n\}, i<j, r\in \ZZ \right\},
\end{align*}
where
\begin{equation}
\begin{matrix}
\begin{tikzpicture}[every node/.style={inner sep=1}, scale=1]
	\node at (.5,0) {\tikz \draw[black, ->, >=implies, double distance=2] (0,0) -- +(-.127,0);};
	\node at (5.5,0) {\tikz \draw[black, ->, >=implies, double distance=2] (0,0) -- +(.127,0);};
	\foreach \x in {0,5}{\foreach \y in {-1,1}{\draw (\x,\y*.05) to (\x+1,\y*.05);}}
	\node[wV, label=below right:{$2\alpha_n^\vee = 2\varepsilon_n^\vee$}, label=above right:{$\alpha_n^\vee = \varepsilon_n^\vee$}] (6) at (6,0) {};
	\node[wV, label=below left:{$2\alpha_0^\vee = -2\varepsilon_1^\vee+K$}, label=above left:{$\alpha_0^\vee = -\varepsilon_1^\vee+\frac12 K$}] (0) at (0,0) {};
	\foreach \x/\y in {1/1,2/2,4/n\!-\!2,5/n\!-\!1}{
		\node[wV] (\x) at (\x,0) {};}
	\node[wV,label=above right:{$\alpha_i^\vee = \varepsilon_i^\vee-\varepsilon_{i+1}^\vee$}] (2) at (2,0) {};
	\draw (1) to (2) (4) to (5);
	\draw[dashed] (2) to (4);
	\end{tikzpicture}\end{matrix}
\label{scoroots}
%\tag{scoroots}
\end{equation}

\begin{remark}
    Throughout this paper we present several diagrams imitating the Dynkin diagram with the labeling related to the coroots (like the one above), roots, or parameters. These are merely intended for conceptual association rather than to specify relations between the objects.
\end{remark}

\subsubsection{The affine roots for type $CC_n$}\label{sec:affroots}

Let
$$Q= \hbox{$\ZZ$-span}\{ \varepsilon_1, \ldots, \varepsilon_n, \hbox{$\frac12$} \delta\}
$$
be the $\ZZ$-vector space spanned by symbols $\varepsilon_1, \ldots, \varepsilon_n$ and $\frac12 \delta$.
The affine Weyl group $W_Y$ is the group of $\ZZ$-linear transformations of $Q$
generated by the transformations $s_0, s_1, \ldots, s_n$ given as follows:
If $\mu = \mu_1 \varepsilon_1 + \cdots + \mu_n \varepsilon_n + \frac{k}{2} \delta$ then
\begin{align}
s_0 \mu &= - \mu_1\varepsilon_1 + \mu_2\varepsilon_2
+\cdots+\mu_n \varepsilon_n + \left(\hbox{$\frac{k}{2}$}+\lambda_1\right) \delta,
\nonumber  \\
s_n\mu &= \mu_1\varepsilon_1+\cdots +\mu_{n-1}\varepsilon_{n-1}
-\mu_n \varepsilon_n + \hbox{$\frac{k}{2}$} \delta, \qquad\hbox{and}
\label{lvl0onaZvee}
%\tag{lvl0onaZvee}
\\
s_i\mu &= \mu\varepsilon_1 + \cdots + \mu_{i-1}\varepsilon_{i-1}
+\mu_{i+1}\varepsilon_i+\mu_i \varepsilon_{i+1} + \mu_{i+2}\varepsilon_{i+2}+\cdots +
\mu_n\varepsilon_n + \hbox{$\frac{k}{2}$} \delta,
\nonumber
\end{align}
for $i\in \{1, \ldots, n-1\}$.  
Let
$$s_{\varepsilon_1} = s_1\cdots s_n\cdots s_1, \qquad
h_{\varepsilon_1^\vee} = s_0  s_{\varepsilon_1} = s_0 s_1\cdots s_n\cdots s_1
\qquad\hbox{and}\qquad
h_{\varepsilon_{i+1}^\vee} = s_i t_{\varepsilon_i^\vee}s_i,
$$
for $i\in \{1, \ldots, n-1\}$.  Then
$h_{\varepsilon_1^\vee}\mu
= s_0 s_{\varepsilon_1}\mu
= \mu_1\varepsilon_1
+\cdots+\mu_n \varepsilon_n + \left(\hbox{$\frac{k}{2}$}-\mu_1\right) \delta$ and
$$
h_{\varepsilon_i^\vee}\mu 
= \mu_1\varepsilon_1
+\cdots+\mu_n \varepsilon_n + \left(\hbox{$\frac{k}{2}$}-\mu_i\right) \delta,
\qquad\hbox{for $i\in \{1, \ldots, n\}$.}
$$
If $\nu^\vee = \nu_1\varepsilon_1^\vee+\cdots \nu_n\varepsilon_n^\vee$ and 
$\mu = \mu_1 \varepsilon_1 + \cdots + \mu_n \varepsilon_n + \frac{k}{2}\delta$ then
\begin{align*}
h_{\nu^\vee} \mu &= h_{\varepsilon_1^\vee}^{\nu_1}\cdots h_{\varepsilon_n^\vee}^{\nu_n} \mu
=\mu - (\hbox{$\frac{k}{2}$}+(\mu_1\nu_1+\cdots+\mu_n\nu_n))\delta 
=\mu + (\hbox{$\frac{k}{2}$}-\langle \mu,\nu^\vee\rangle )\delta,
\end{align*}
and special cases of this last formula are
$$h_{\nu^\vee} \varepsilon_i = \varepsilon_i - \nu_i \delta
\qquad\hbox{and}\qquad
h_{\nu^\vee} \delta = \delta.
$$

The set of roots $S$ for type $CC_n$ is the union of the five $W_Y$-orbits given by
\begin{align*}
O_1 
&= W_Y\cdot\alpha_n = W_Y\cdot \varepsilon_n   
= \{ \pm \varepsilon_i+r\delta \ |\  i\in \{1, \ldots, n\}, r\in \ZZ\}, \\
2O_1 &= O_2 = W_Y\cdot 2\alpha_n = W_Y\cdot 2\varepsilon_n
= \{ \pm 2\varepsilon_i + 2r\delta \ |\  i\in \{1, \ldots, n\}, r\in \ZZ\} , \\
O_3 
&= W_Y\cdot \alpha_n = W_Y\cdot (-\varepsilon_1+ \hbox{$\frac12$} \delta)
= \{ \pm(\varepsilon_i+\hbox{$\frac12$}(2r+1)\delta \ |\  i\in \{1, \ldots, n\}, r\in \ZZ\}, \\
2O_3 &= O_4 
= W_Y\cdot 2\alpha_n = W_Y\cdot (-2\varepsilon_1 + \delta)
= \{ \pm 2\varepsilon_i+(2r+1) \delta \ |\  i\in \{1, \ldots, n\}, r\in \ZZ\} , \\
O_5  
&= W_Y\cdot \alpha_1 = W_Y\cdot (\varepsilon_1-\varepsilon_2)
= 
\left\{ \begin{array}{l}
\pm(\varepsilon_i+\varepsilon_j)+r\delta  \\
\pm(\varepsilon_i-\varepsilon_j)+r\delta  
\end{array}\ \Big\vert\ 
i,j\in \{1, \ldots, n\}, i<j, r\in \ZZ \right\},
\end{align*}
where
\begin{equation}
\begin{matrix}
\begin{tikzpicture}[every node/.style={inner sep=1}, scale=1]
	\node at (.5,0) {\tikz \draw[black, ->, >=implies, double distance=2] (0,0) -- +(-.127,0);};
	\node at (5.5,0) {\tikz \draw[black, ->, >=implies, double distance=2] (0,0) -- +(.127,0);};
	\foreach \x in {0,5}{\foreach \y in {-1,1}{\draw (\x,\y*.05) to (\x+1,\y*.05);}}
	\node[wV, label=below right:{$2\alpha_n = 2\varepsilon_n$}, label=above right:{$\alpha_n = \varepsilon_n$}] (6) at (6,0) {};
	\node[wV, label=below left:{$2\alpha_0 = -2\varepsilon_1+K$}, label=above left:{$\alpha_0 = -\varepsilon_1+\frac12 \delta$}] (0) at (0,0) {};
	\foreach \x/\y in {1/1,2/2,4/n\!-\!2,5/n\!-\!1}{
		\node[wV] (\x) at (\x,0) {};}
	\node[wV,label=above right:{$\alpha_i = \varepsilon_i-\varepsilon_{i+1}$}] (2) at (2,0) {};
	\draw (1) to (2) (4) to (5);
	\draw[dashed] (2) to (4);
	\end{tikzpicture}\end{matrix}
\label{sroots}
%\tag{sroots}
\end{equation}

\subsection{Other classical types}\label{ssec:paramspec}

With the notations as in Section~\ref{sec:affroots}, each affine root system of classical type is a subset of $Q$.
The irreducible affine root systems of classical type (and the appropriate specializations for obtaining the 
Macdonald polynomials of each type from the Macdonald polynomials of type $(C^\vee_n,C_n)$) are given by the
diagram in Section~\ref{affrootsystemsfig}.
The middle notation for each root system is the notation in Macdonald~\cite[\S 1.3]{Mac03}, the right notation
is that of Bruhat-Tits~\cite{BT72} and the left notation is that of Kac~\cite[Ch.\ 6]{Kac}.

To determine the specializations, we look at the kernel for the inner product (see Section~\ref{ssec:kernel}).  Following~\cite[(5.1.3)]{Mac03} the orthogonal polynomials are determined by the inner product which, in turn, is
determined by factors of the form in~\eqref{bothorbits} corresponding to orbits of roots.

For instance, when both $O_1$ and $2O_1$
are present then the factor corresponding to the root $\varepsilon_n$ is
\begin{equation}
\frac{1}{\Delta_{\varepsilon_n} \Delta_{2\varepsilon_n} } 
= \frac{ (1-t_n^{\frac12}u_n^{\frac12}x_n)(1+t_n^{\frac12}u_n^{-\frac12}x_n)}{1-x_n^2}.
\label{bothorbits}
%\tag{bothorbits}
\end{equation}
(The notation $\frac{1}{\Delta_{\varepsilon_n}\Delta_{2\varepsilon_n}}$ for this factor is as in~\cite[(5.1.3)]{Mac03};
the notation that we use for this factor in Section~\ref{ssec:kernel} is $\kappa_{\varepsilon_n}^X$.)
If only the orbit $O_1$ is present then the factor is
$$\frac{1}{\Delta_{\varepsilon_n}} 
= \frac{1-t_nx_n}{1-x_n}
\qquad\hbox{which is obtained by specializing $t_n^{\frac12} = u_n^{\frac12}$ in~\eqref{bothorbits}.}
$$
If only the orbit $2O_1$ is present then the factor is
$$\frac{1}{\Delta_{2\varepsilon_n}} 
= \frac{1-t_nx_n^2}{1-x_n^2}
\qquad\hbox{which is obtained by specializing $u_n^{\frac12}=1$ in~\eqref{bothorbits}.}
$$
In this way, the parameter specializations in the diagram in Section~\ref{affrootsystemsfig} are determined by
which orbits of roots are present in the root system.

\section{$c$-functions and DAHA relations}\label{sec: c-funct and DAHA}

This section collects the tools for working with Koornwinder polynomials as polynomials in $x_1^{\pm1}, \ldots, x_n^{\pm1}$ that depend on 6 parameters.  The root system
of type $CC_n$ provides the structure for organizing the many symmetries between the variables and
the various parameters and this section specifies carefully the links to the root system.  The $c$-functions, introduced in Section~\ref{sec:cfnY}, are a core structure to providing explicit formulas for Koornwinder expansions, specializations and
norm formulas.  With the notation
for the $c$-functions in hand, Section~\ref{sec:DAHA}
describes briefly the relations of the double affine Hecke algebra.  These relations provide a convenient
summary of the calculus of the operators on polynomials which are used in the rest of the paper.

\subsection{Parameters}\label{ssec:params}

Let
$$
q, t^{\frac12}, t_0^{\frac12}, u_0^{\frac12}, t_n^{\frac12}, u_n^{\frac12} \quad\hbox{be independent parameters,}
$$
and let
$$\KK = \CC(q, t^{\frac12}, t_0^{\frac12}, u_0^{\frac12}, t_n^{\frac12}, u_n^{\frac12} )
\quad\hbox{be the field of fractions of $\CC[q, t^{\frac12}, t_0^{\frac12}, u_0^{\frac12}, t_n^{\frac12}, u_n^{\frac12} ]$.}  
$$
The field $\KK$ will be the base field for most algebras in this paper.

Recalling the simple coroots and simple roots from~\eqref{scoroots} and~\eqref{sroots}, set
\begin{align}
&\begin{matrix}
\begin{tikzpicture}[every node/.style={inner sep=1}, scale=1]
	\node at (.5,0) {\tikz \draw[black, ->, >=implies, double distance=2] (0,0) -- +(.127,0);};
	\node at (5.5,0) {\tikz \draw[black, ->, >=implies, double distance=2] (0,0) -- +(.127,0);};
	\foreach \x in {0,5}{\foreach \y in {-1,1}{\draw (\x,\y*.05) to (\x+1,\y*.05);}}
	\node[wV, label=below right:{$u_{\alpha_n^\vee} = t_0$}, label=above right:{$t_{\alpha_n^\vee} = t_n$}] (6) at (6,0) {};
	\node[wV, label=below left:{$u_{\alpha_0^\vee} = u_0$}, label=above left:{$t_{\alpha_0^\vee} = u_n$}] (0) at (0,0) {};
	\foreach \x/\y in {1/1,2/2,4/n\!-\!2,5/n\!-\!1}{
		\node[wV] (\x) at (\x,0) {};}
	\node[wV,label=above right:{$t_{\alpha_i^\vee} = u_{\alpha_i^\vee} = t$}] (2) at (2,0) {};
	\draw (1) to (2) (4) to (5);
	\draw[dashed] (2) to (4);
	\end{tikzpicture}
\end{matrix}
\label{cortparams}
%\tag{cortparams}
 \intertext{and}
&\begin{matrix}
\begin{tikzpicture}[every node/.style={inner sep=1}, scale=1]
	\node at (.5,0) {\tikz \draw[black, ->, >=implies, double distance=2] (0,0) -- +(.127,0);};
	\node at (5.5,0) {\tikz \draw[black, ->, >=implies, double distance=2] (0,0) -- +(.127,0);};
	\foreach \x in {0,5}{\foreach \y in {-1,1}{\draw (\x,\y*.05) to (\x+1,\y*.05);}}
	\node[wV, label=below right:{$u_{\alpha_n} = u_n$}, label=above right:{$t_{\alpha_n} = t_n$}] (6) at (6,0) {};
	\node[wV, label=below left:{$u_{\alpha_0} = u_0$}, label=above left:{$t_{\alpha_0} = t_0$}] (0) at (0,0) {};
	\foreach \x/\y in {1/1,2/2,4/n\!-\!2,5/n\!-\!1}{
		\node[wV] (\x) at (\x,0) {};}
	\node[wV,label=above right:{$t_{\alpha_i} = u_{\alpha_i} = t$}] (2) at (2,0) {};
	\draw (1) to (2) (4) to (5);
	\draw[dashed] (2) to (4);
	\end{tikzpicture}\end{matrix}
\label{rtparams}
%\tag{rtparams}
\end{align}

Let $S^\vee$ be the set of affine coroots and let $S$ be the set of affine roots.
Define $t_{\alpha^\vee}$, $u_{\alpha^\vee}$ and $t_\alpha$, $u_\alpha$ for arbitrary coroots $\alpha^\vee$ and 
arbitrary roots $\alpha$ by requiring
\begin{equation}
t_{w\alpha} = t_{\alpha} \qquad\hbox{and}\qquad
t_{w\alpha^\vee} = t_{\alpha^\vee},
\qquad
\hbox{for $w\in W$, $\alpha^\vee\in S^\vee$ and $\alpha\in S$.}
\label{paraminv}
%\tag{paraminv}
\end{equation}
The difference between the parameters for the coroots in~\eqref{cortparams} and the parameters for the roots in~\eqref{rtparams}
corresponds exactly to introduction of the ``dual labels'' in~\cite[(1.5.1)]{Mac03}.

\subsection{The polynomial rings $\KK[X]$ and $\KK[Y]$}

Let $\KK[X] = \KK[X_1^{\pm1}, \ldots, X_n^{\pm1}]$ be the Laurent polynomial
ring in the variables $X_1, \ldots, X_n$.  
Identify $\KK[X]$ with the group algebra of $Q= \hbox{$\ZZ$-span}\{ \varepsilon_1, \ldots, \varepsilon_n, \hbox{$\frac12$} \delta\}$
via the notations
$$q^{\frac12} = X^{\frac12\delta} \quad\hbox{and}\quad
X_i = X^{\varepsilon_i},
\quad\hbox{and}\quad
q^{\frac{k}{2}} X_1^{\mu_1}\cdots X_n^{\mu_n}
=X^{\frac{k}{2}\delta+\mu_1\varepsilon_1+\cdots+\mu_n\varepsilon_n}=X^\mu,
$$
for $i\in \{1, \ldots, n\}$ and $\mu = \mu_1\varepsilon_1+\cdots+\mu_n\varepsilon_n+\frac{k}{2}\delta
\in Q$. 
The image of the simple roots in $\KK[X]$ is given by
\begin{equation}
\begin{matrix}
\begin{tikzpicture}[every node/.style={inner sep=1}, scale=1]
	\node at (.5,0) {\tikz \draw[black, ->, >=implies, double distance=2] (0,0) -- +(.127,0);};
	\node at (5.5,0) {\tikz \draw[black, ->, >=implies, double distance=2] (0,0) -- +(.127,0);};
	\foreach \x in {0,5}{\foreach \y in {-1,1}{\draw (\x,\y*.05) to (\x+1,\y*.05);}}
	\node[wV, label=above right:{$X^{\alpha_n} = X_n$}] (6) at (6,0) {};
	\node[wV, label=above left:{$X^{\alpha_0} = q^{\frac12} X_1^{-1}$ }] (0) at (0,0) {};
	\foreach \x/\y in {1/1,2/2,4/n\!-\!2,5/n\!-\!1}{
		\node[wV] (\x) at (\x,0) {};}
	\node[wV,label=above right:{$X^{\alpha_i} = X_iX_{i+1}^{-1}$}] (2) at (2,0) {};
	\draw (1) to (2) (4) to (5);
	\draw[dashed] (2) to (4);
	\end{tikzpicture}\end{matrix}
\label{Xsroots}
%\tag{Xsroots}
\end{equation}

Let $\KK[Y] = \KK[Y_1^{\pm1}, \ldots, Y_n^{\pm1}]$ be the Laurent polynomial
ring in the variables $Y_1, \ldots, Y_n$.  
Identify $\KK[Y]$ with the group algebra of $Q^\vee 
= \hbox{$\ZZ$-span}\{ \varepsilon_1^\vee, \ldots, \varepsilon_n^\vee, \hbox{$\frac12$} K\}$
via the notations
$$q^{\frac12} = Y^{-\frac12 K} \quad\hbox{and}\quad
Y_i = Y^{\varepsilon_i^\vee},
\quad\hbox{and}\quad
q^{-\frac{k}{2}} Y_1^{\lambda_1}\cdots Y_n^{\lambda_n}
=Y^{\frac{k}{2}K+\lambda_1\varepsilon_1^\vee+\cdots+\lambda_n\varepsilon_n^\vee}=Y^{\lambda^\vee},
$$
for $i\in \{1, \ldots, n\}$ and $\lambda^\vee 
= \lambda_1\varepsilon_1^\vee+\cdots+\lambda_n\varepsilon_n^\vee+\frac{k}{2}K
\in  Q^\vee$. 
The image of the simple coroots in $\KK[Y]$ is given by
\begin{equation}
\begin{matrix}
\begin{tikzpicture}[every node/.style={inner sep=1}, scale=1]
	\node at (.5,0) {\tikz \draw[black, ->, >=implies, double distance=2] (0,0) -- +(.127,0);};
	\node at (5.5,0) {\tikz \draw[black, ->, >=implies, double distance=2] (0,0) -- +(.127,0);};
	\foreach \x in {0,5}{\foreach \y in {-1,1}{\draw (\x,\y*.05) to (\x+1,\y*.05);}}
	\node[wV, label=above right:{$Y^{\alpha_n^\vee} = Y_n$}] (6) at (6,0) {};
	\node[wV, label=above left:{$Y^{\alpha_0^\vee} = q^{-\frac12} Y_1^{-1}$ }] (0) at (0,0) {};
	\foreach \x/\y in {1/1,2/2,4/n\!-\!2,5/n\!-\!1}{
		\node[wV] (\x) at (\x,0) {};}
	\node[wV,label=above right:{$Y^{\alpha_i^\vee} = Y_iY_{i+1}^{-1}$}] (2) at (2,0) {};
	\draw (1) to (2) (4) to (5);
	\draw[dashed] (2) to (4);
	\end{tikzpicture}\end{matrix}
\label{Ysroots}
%\tag{Ysroots}
\end{equation}

\subsection{$c$-functions} \label{sec:cfnY}

Let $\KK(Y)$ be the field of fractions of the Laurent polynomial ring
$\KK[Y]$. % = \KK[Y_1^{\pm1}, \ldots, Y_n^{\pm1}]$.
For a coroot $\alpha^\vee$ let (see~\cite[(4.2.2) and (4.3.9)]{Mac03}),
\begin{align}
c_{\alpha^\vee}^Y = t_{\alpha^\vee}^{-\frac12}
\frac{(1-t_{\alpha^\vee}^{\frac12}u_{\alpha^\vee}^{\frac12}Y^{\alpha^\vee})
(1+t_{\alpha^\vee}^{\frac12}u_{\alpha^\vee}^{-\frac12}Y^{\alpha^\vee})}
{(1-Y^{2\alpha^\vee})}
\qquad &\hbox{and}\qquad
\kappa_{\alpha^\vee}^Y = t_{\alpha^\vee}^{\frac12} c_{\alpha^\vee}^Y.
\label{cfnYdefn}
%\tag{cfnYdefn} 
\\
\hbox{If $t_{\alpha^\vee}^{\frac12}=u_{\alpha^\vee}^{\frac12}$ \quad then}\qquad 
c_{\alpha^\vee}^Y = \frac{t_{\alpha^\vee}^{-\frac12}-t_{\alpha^\vee}^{\frac12}Y^{\alpha^\vee}}{1-Y^{\alpha^\vee}}
\qquad &\hbox{and}\qquad
\kappa_{\alpha^\vee}^Y = \frac{1 - t_{\alpha^\vee}Y^{\alpha^\vee}}{1-Y^{\alpha^\vee}}. \notag
\end{align}
(More accurately, the function $c_{\alpha^\vee}^Y$ should be considered as a local factor of a $c$-function, see
\cite{Sto11}.)
The expression $\kappa_{\alpha^\vee}^Y$ is a slightly renormalized version of the $c$-function $c_{\alpha^\vee}^Y$ which,
although not technically necessary, is immensely helpful for making the formulas more palatable.

In general, for arbitrary $t_{\alpha^\vee}^{\frac12}$ and $u_{\alpha^\vee}^{\frac12}$,
\begin{align}
c_{\alpha^\vee}^Y &+c_{-\alpha^\vee}^Y
= 
\frac{(1-t_{\alpha^\vee}^{\frac12}u_{\alpha^\vee}^{\frac12}Y^{\alpha^\vee})
(1+t_{\alpha^\vee}^{\frac12}u_{\alpha^\vee}^{-\frac12}Y^{\alpha^\vee})}
{t_{\alpha^\vee}^{\frac12}(1-Y^{2\alpha^\vee})}
+\frac{(1-t_{\alpha^\vee}^{\frac12}u_{\alpha^\vee}^{\frac12}Y^{-\alpha^\vee})
(1+t_{\alpha^\vee}^{\frac12}u_{\alpha^\vee}^{-\frac12}Y^{-\alpha^\vee})}
{t_{\alpha^\vee}^{\frac12}(1-Y^{-2\alpha^\vee})}  \nonumber  \\
&= 
\frac{(1-t_{\alpha^\vee}^{\frac12}u_{\alpha^\vee}^{\frac12}Y^{\alpha^\vee})(1+t_{\alpha^\vee}^{\frac12}u_{\alpha^\vee}^{-\frac12}Y^{\alpha^\vee})}
{t_{\alpha^\vee}^{\frac12}(1-Y^{2\alpha^\vee})}
+\frac{t_{\alpha^\vee}(1-t_{\alpha^\vee}^{-\frac12}u_{\alpha^\vee}^{-\frac12}Y^{\alpha^\vee})
(1+t_{\alpha^\vee}^{-\frac12}u_{\alpha^\vee}^{\frac12}Y^{\alpha^\vee})}
{t_{\alpha^\vee}^{\frac12}(1-Y^{2\alpha^\vee})} \nonumber  \\
&= 
\frac{1-t_{\alpha^\vee}^{\frac12}u_{\alpha^\vee}^{\frac12}Y^{\alpha^\vee}+t_{\alpha^\vee}^{\frac12}u_{\alpha^\vee}^{-\frac12}Y^{\alpha^\vee}
-t_{\alpha^\vee}Y^{2\alpha^\vee} +t_{\alpha^\vee} - t_{\alpha^\vee}^{\frac12}u_{\alpha^\vee}^{-\frac12}Y^{\alpha^\vee}
+t_{\alpha^\vee}^{\frac12}u_{\alpha^\vee}^{\frac12}Y^{\alpha^\vee}- Y^{2\alpha^\vee})}
{t_{\alpha^\vee}^{\frac12}(1-Y^{2\alpha^\vee})} \nonumber  \\
&= 
\frac{(1+t_{\alpha^\vee})(1-Y^{2\alpha^\vee})}
{t_{\alpha^\vee}^{\frac12}(1-Y^{2\alpha^\vee})}  = t_{\alpha^\vee}^{\frac12}+t_{\alpha^\vee}^{-\frac12}.
\label{cplusc}
%\tag{cplusc}
\end{align}

Let $w\in W$ and let $w=s_{i_1}\cdots s_{i_\ell}$ be a reduced word for $w$.
The \emph{coroot sequence} of the reduced word $w = s_{i_1}\cdots s_{i_\ell}$ is
\begin{equation}
\hbox{the sequence $(\beta^\vee_k\ |\ \hbox{$k\in \{1, \ldots, \ell\}$ and $i_k\ne \pi$} \})$ 
given by}\qquad
\beta^\vee_k = s_{i_\ell}^{-1}\cdots s^{-1}_{i_{k+1}}\alpha^\vee_{i_k}.
\end{equation}
Then define
\begin{equation}
\kappa_w^Y = \prod_{\beta^\vee\in \mathrm{Inv}(w)} \kappa^Y_{\beta^\vee},
\qquad\hbox{where}\qquad
\mathrm{Inv}(w) = \{\beta_1^\vee, \ldots, \beta_\ell^\vee\}
\label{kappawdef}
%\tag{kappawdef}
\end{equation}
is the set of elements in a coroot sequence for a reduced word for $w$.
If $w=uv$ with $\ell(u)+\ell(v) = \ell(w)$ then the 
coroot sequence of $w$ is 
$v^{-1}$ times the coroot sequence of $u$
followed by 
the coroot sequence of $v$ 
so that
$$\mathrm{Inv}(uv) = v^{-1}\mathrm{Inv}(u) \cup \mathrm{Inv}(v).$$

The inversions of elements of $W_{\mathrm{fin}}$ come in two types: `droite' and `standard'.  Using the indicators
$d$ for `droite' and $s$ for `standard', let
$$(S^\vee_{0,d})^+ = \{ \varepsilon_1^\vee, \ldots, \varepsilon_n^\vee\}
\qquad\hbox{and}\qquad
(S^\vee_{0,s})^+ = \left\{ \begin{array}{l} \varepsilon_i^\vee-\varepsilon_j^\vee, \\ \varepsilon_i^\vee+\varepsilon_j^\vee \end{array}\  \Big\vert\ 
\hbox{$i,j\in \{1, \ldots, n\}$ with $i<j$}\right\},
$$
and for $v\in W_{\mathrm{fin}}$, define
$$\mathrm{Inv}_d(v) = \mathrm{Inv}(v)\cap (S^\vee_{0,d})^+
\quad\hbox{and}\quad
\mathrm{Inv}_s(v) = \mathrm{Inv}(v)\cap (S^\vee_{0,s})^+.
$$

For $v\in W_{\mathrm{fin}}$ define 
\begin{align}
\kappa^+_v 
&= \left(\prod_{\beta^\vee\in \mathrm{Inv}_s(v)} \kappa_{\beta^\vee}^Y\right) 
\left(\prod_{\beta^\vee\in \mathrm{Inv}_d(v)} \kappa_{\beta^\vee}^Y\right),
&\kappa^{\pm}_v 
&= \left(\prod_{\beta^\vee\in \mathrm{Inv}_s(v)} \kappa_{\beta^\vee}^Y\right) 
\left(\prod_{\beta^\vee\in \mathrm{Inv}_d(v)} \kappa_{\beta^\vee}^{Y^{-1}}\right),
\nonumber \\
\kappa^-_v 
&= \left(\prod_{\beta^\vee\in \mathrm{Inv}_s(v)} \kappa_{\beta^\vee}^{Y^{-1}}\right) 
\left(\prod_{\beta^\vee\in \mathrm{Inv}_d(v)} \kappa_{\beta^\vee}^{Y^{-1}}\right),
&\kappa^{\mp}_v
&= \left(\prod_{\beta^\vee\in \mathrm{Inv}_s(w)} \kappa_{\beta^\vee}^{Y^{-1}}\right) 
\left(\prod_{\beta^\vee\in \mathrm{Inv}_d(w)} \kappa_{\beta^\vee}^{Y}\right).
\label{kXi}
%\tag{kXi}
\end{align}

Finally, define  also
$$
\kappa_\mathrm{dr}^{Y} = \prod_{i=1}^n \kappa_{\varepsilon_i^\vee}^{Y}
=  \prod_{i=1}^n \frac{(1-t_n^{\frac12}t_0^{\frac12}Y_i)(1+t_n^{\frac12}t_0^{-\frac12}Y_i) }{ (1-Y_i^{2})}
$$
and 
\begin{equation}
\kappa_\mathrm{st}^{Y} 
= \prod_{1\le i<j\le n}^n \kappa_{\varepsilon_i^\vee-\varepsilon_j^\vee}^{Y} \kappa_{\varepsilon_i^\vee+\varepsilon_j^\vee}^{Y} 
= \prod_{1\le i<j\le n}^n \frac{(1-t Y_iY_j^{-1})}{(1-Y_iY_j^{-1})}
\frac{(1-t Y_iY_j)}{(1-Y_iY_j)}
\label{cstcdr}
%\tag{cstcdr}
\end{equation}
so that
\begin{equation}
\kappa_{w_0}^Y = \kappa_{\mathrm{dr}}^Y \kappa_{\mathrm{st}}^Y.
\label{kw0Y}
%\tag{kw0Y}
\end{equation}

This subsection has presented the $c$-functions and related functions $\kappa_{\beta^\vee}^Y$ 
 in terms of the $\{Y_1,\cdots,Y_n\}$ variables. We will also consider these notions in other sets of variables, 
 like $\{Y_1^{-1},\cdots,Y_n^{-1}\}$, $\{X_1,\cdots,X_n\}$, and $\{X^{-1}_1,\cdots,X^{-1}_n\}$,
and use notations like $\kappa_{\beta^\vee}^{Y^{-1}}$, $\kappa_\beta^X$, $\kappa_\beta^{X^{-1}}$, respectively.
For example
\begin{equation}
\kappa^{X^{-1}}_{w_0}
= \left(\prod_{i=1}^n \frac{(1-t_n^{\frac12} u_n^{\frac12}X_i^{-1})(1+t_n^{\frac12}u_n^{-\frac12}X_i^{-1})}
{(1-X_i^{-2})}\right)
\left( \prod_{1\le i<j\le n} \frac{(1-tX_i^{-1}X_j)(1-tX_i^{-1}X_j^{-1})}{(1-X_i^{-1}X_j)(1-X_i^{-1}X_j^{-1})}\right).
\label{kw0X}
%\tag{kw0X}
\end{equation}

\subsection{The algebras $\widetilde H_{\mathrm{loc}}$ and $\widetilde H_{\mathrm{int}}$}\label{sec:DAHA}

The Koornwinder polynomials are elements of the polynomial ring $\KK[x_1^{\pm1}, \ldots, x_n^{\pm1}]$ which
are characterized, up to normalization, by the fact that they are eigenvectors for the Cherednik-Dunkl
operators $Y_1, \ldots, Y_n$ (see~\eqref{Eeig}). However the $c$-functions, which form the core calculus
for working with Koornwinder polynomials, are elements of $\KK(X)$, the field of fractions of the polynomial ring. Thus extending from $\KK[X]$ to $\KK(X)$ is necessary for handling the tools.

The Cherednik-Dunkl operators are elements of the double affine Hecke algebra $\widetilde{H}_{\mathrm{int}}$, 
which is formed by pasting the two polynomial rings $\KK[X]$ and $\KK[Y]$ together with a finite Hecke algebra $H_{\mathrm{fin}}$.
However, the right home for the $c$-functions and the operators for creating Koornwinder polynomials
is a larger algebra $\widetilde{H}_{\mathrm{loc}}$ which extends the algebra $\widetilde{H}_{\mathrm{int}}$
by extending $\KK[X]$ and $\KK[Y]$ to the fraction fields $\KK(X)$ and $\KK(Y)$.  In this subsection we introduce the
algebras $\widetilde{H}_{\mathrm{loc}}$ and $\widetilde{H}_{\mathrm{int}}$ by generators and relations.

Let $\KK(X)$ and $\KK(Y)$ be the fraction fields of $\KK[X]$ and $\KK[Y]$ respectively.
Recall that $\KK[X]$ is the group algebra of $Q$ and $\KK[Y]$ is the group algebra of $Q^\vee$, that
$$X^{\frac12\delta} = q\qquad\hbox{and}\qquad Y^{-\frac12 K} = q^{\frac12}$$
and that
$W_X$ acts on $Q^\vee$ and $W_Y$ acts on $Q$ by the formulas given in~\eqref{lvl0onaZvee} and~\eqref{lvl0onaZ}
so that
\begin{equation}
wX^\mu = X^{w\mu}
\qquad\hbox{and}\qquad zY^{\lambda^\vee} = Y^{z\lambda^\vee},
\label{Wactions}
%\tag{Wactions}
\end{equation}
for $w\in W_Y$, $z\in W_X$, $\lambda^\vee\in Q^\vee$ and $\mu\in Q$.

Let $\widetilde{H}_{\mathrm{loc}}$ be the $\KK$-algebra generated by $\eta_{s^\vee_0}, \ldots, \eta_{s^\vee_n}$, $\xi_{s_0}, \ldots, \xi_{s_n}$,
$T_1, \ldots, T_n$, 
$T_{\alpha_0},  T_{\alpha_1}, \ldots, T_{\alpha_n}$ and $T_{\alpha_0^\vee},  T_{\alpha_1^\vee}, \ldots, T_{\alpha_n^\vee}$ 
and $T_{s_{\varepsilon_1}}$ and  $\KK(X)$ and $\KK(Y)$ with relations 
\begin{equation}
\begin{matrix}
	\begin{matrix}\begin{tikzpicture}[every node/.style={inner sep=1}, scale=1]
	\node at (.5,0) {\tikz \draw[black, double distance=2] (0,0) -- +(.127,0);};
	\node at (5.5,0) {\tikz \draw[black, double distance=2] (0,0) -- +(-.127,0);};
	\foreach \x in {0,5}{\foreach \y in {-1,1}{\draw (\x,\y*.05) to (\x+1,\y*.05);}}
	\node[wV, label=above:{$\eta_{s_0^\vee}$}] (0) at (0,0) {};
	\node[wV, label=above:{$\eta_{s_n^\vee}$}] (6) at (6,0) {};
	\foreach \x/\y in {1/1,2/2,4/n\!-\!2,5/n\!-\!1}{
		\node[wV] (\x) at (\x,0) {};}
	\draw (1) to (2) (4) to (5);
	\draw[dashed] (2) to (4);
	\node[wV, label=above:{$\eta_{s_1^\vee}$}] (1) at (1,0) {};
	\node[wV, label=above:{$\eta_{s_2^\vee}$}] (2) at (2,0) {};
	\node[wV, label=above:{$\eta_{s_{n-2}^\vee}$}] (4) at (4,0) {};
	\node[wV, label=above:{$\eta_{s_{n-1}^\vee}$}] (5) at (5,0) {};
	\end{tikzpicture}\end{matrix}
\qquad \qquad 
&\eta_{s_i}^2=1,\qquad
&\eta_{s_i}Y^{\lambda^\vee} = Y^{s_i\lambda^\vee}\eta_{s_i}
\\
	\begin{matrix}\begin{tikzpicture}[every node/.style={inner sep=1}, scale=1]
	\node at (.5,0) {\tikz \draw[black, double distance=2] (0,0) -- +(.127,0);};
	\node at (5.5,0) {\tikz \draw[black, double distance=2] (0,0) -- +(-.127,0);};
	\foreach \x in {0,5}{\foreach \y in {-1,1}{\draw (\x,\y*.05) to (\x+1,\y*.05);}}
	\node[wV, label=above:{$\xi_{s_0}$}] (0) at (0,0) {};
	\node[wV, label=above:{$\xi_{s_n}$}] (6) at (6,0) {};
	\foreach \x/\y in {1/1,2/2,4/n\!-\!2,5/n\!-\!1}{
		\node[wV] (\x) at (\x,0) {};}
	\draw (1) to (2) (4) to (5);
	\draw[dashed] (2) to (4);
	\node[wV, label=above:{$\xi_{s_1}$}] (1) at (1,0) {};
	\node[wV, label=above:{$\xi_{s_2}$}] (2) at (2,0) {};
	\node[wV, label=above:{$\xi_{s_{n-2}}$}] (4) at (4,0) {};
	\node[wV, label=above:{$\xi_{s_{n-1}}$}] (5) at (5,0) {};
	\end{tikzpicture}\end{matrix}
\qquad \qquad 
&\xi_{s_i}^2=1,\qquad
&\xi_{s_i}X^\mu = X^{s_i\mu}\xi_{s_i},
\end{matrix}
\label{Drels}
%\tag{Drels}
\end{equation}
for $i\in\{0,1,\ldots, n\}$, $\lambda^\vee\in Q^\vee$ and $\mu\in Q$, and 
\begin{gather}
T_{\alpha^\vee_i}+t_{\alpha^\vee_i}^{-\frac12} = (1+\eta_{s_i^\vee})c_{\alpha_i^\vee}^Y
\quad\hbox{and}\quad
T_{\alpha_i}+t_{\alpha_i}^{-\frac12} = (1+\xi_{s_i})c_{\alpha_i}^{X^{-1}},
\qquad\hbox{for $i\in \{0,1,\ldots, n\}$,}
\label{TdefH} 
\\[0.25cm]
%\tag{TdefH}
T_j = T_{\alpha_j} = T_{\alpha^\vee_j} \quad\hbox{for $j\in \{1, \ldots, n\}$}
\qquad\hbox{and}\qquad
T_{s_{\varepsilon_1}} = T_1\cdots T_{n-1}T_nT_{n-1}\cdots T_1,
\label{glue} \\[0.25cm]
%\tag{glue}
Y^{\varepsilon_1^\vee} = T_{\alpha_0}T_{s_{\varepsilon_1}}
\qquad
X^{\varepsilon_1} = (T_{\alpha_0^\vee})^{-1}T_{s_{\varepsilon_1}}^{-1}
\\[0.25cm]
Y^{\varepsilon_{j+1}^\vee} = T_{\alpha_j}^{-1}Y^{\varepsilon_j^\vee}T_{\alpha_j}^{-1},
\qquad\hbox{and}\qquad
X^{\varepsilon_{j+1}}
= T_{\alpha_j} X^{\varepsilon_j}T_{\alpha_j},
\qquad\hbox{for $j\in \{1, \ldots, n-1\}$,}
\label{ud}
\\[0.25cm]
%\tag{ud}
Y^{-\varepsilon_1^\vee}X^{\varepsilon_1}
= q^{\frac12}(u_0^{\frac12}-u_0^{-\frac12}) T_{s_{\varepsilon_1}}^{-1}
+ qT_{s_{\varepsilon_1}}^{-1} X^{-\varepsilon_1}  Y^{\varepsilon_1}T_{s_{\varepsilon_1}}^{-1}.
\label{Hbrg}
%\tag{Hbrg}
\end{gather}

This presentation of $\widetilde H_{\mathrm{loc}}$ is not minimal as there are many redundant generators and
many redundant relations.  It is designed to specify notations and list the relations that we will need, and to motivate
the operators on polynomials which are the main tools for working with Macdonald polynomials in general.
To be precise, since
$$(1+\xi_{s_i})c_{\alpha_i}^{X^{-1}}
=T_{\alpha_i}+t_{\alpha_i}^{-\frac12}
=(T_{\alpha_i}-t_{\alpha_i}^{\frac12})+t_{\alpha_i}^{-\frac12}+t_{\alpha_i}^{\frac12}
=  (T_{\alpha_i}-t_{\alpha_i}^{\frac12})+c_{\alpha_i}^{X^{-1}}+c_{-\alpha_i}^{X^{-1}}
$$
then
\begin{equation}
T_{\alpha_i}-t^{\frac12} = -c_{-\alpha_i}^{X^{-1}} + \xi_{s_i}c_{\alpha_i}^{X^{-1}}
= - c_{-\alpha_i}^{X^{-1}}(1-\xi_{s_i})
= - c_{\alpha_i}^{X}(1-\xi_{s_i}),
\label{Talt}
%\tag{Talt}
\end{equation}
and this is the formula used to define the action of the double affine Hecke algebra on $\KK[X]$ in~\eqref{TiopsinKX}.

The \emph{double affine Hecke algebra (DAHA)} is the subalgebra
$\widetilde H_{\mathrm{int}}$ inside $\widetilde H_{\mathrm{loc}}$ 
$$\hbox{generated by 
$X^{\pm1}_1, \ldots, X^{\pm1}_n$, $Y_1^{\pm1}, \ldots, Y_n^{\pm1}$, and
$T_1, \ldots, T_n$.}$$
The algebra $\widetilde H_{\mathrm{int}}$ is an integral form of $\widetilde H_{\mathrm{loc}}$ (alternatively, the algebra 
$\widetilde H_{\mathrm{loc}}$ is a completion, or localized form, of $\widetilde H_{\mathrm{int}}$).
A common definition of the DAHA uses the relations listed in the following proposition, which follow without difficulty
from the defining relations of $\widetilde{H}_{\mathrm{loc}}$.

\begin{prop} Let $i\in \{0,1,\ldots, n\}$, $\lambda^\vee\in Q^\vee$ and $\mu\in Q$.  Let
$T_0^\# = q^{-\frac12} T_{\alpha_0^\vee}^{-1}T_{s_{\varepsilon_1}}^{-1} T_{\alpha_0}^{-1}$.
Then
\begin{gather}
  \begin{tikzpicture}[every node/.style={inner sep=1}, scale=1]
	\node at (.5,0) {\tikz \draw[black, double distance=2] (0,0) -- +(.127,0);};
	\node at (5.5,0) {\tikz \draw[black, double distance=2] (0,0) -- +(-.127,0);};
	\foreach \x in {0,5}{\foreach \y in {-1,1}{\draw (\x,\y*.05) to (\x+1,\y*.05);}}
	\node[wV, label=above:{$T_{\alpha_0^\vee}$}] (0) at (0,0) {};
	\node[wV, label=above:{$T_{\alpha_n^\vee}$}] (6) at (6,0) {};
	\foreach \x/\y in {1/1,2/2,4/n\!-\!2,5/n\!-\!1}{
		\node[wV] (\x) at (\x,0) {};}
	\draw (1) to (2) (4) to (5);
	\draw[dashed] (2) to (4);
	\node[wV, label=above:{$T_{\alpha_1}^\vee$}] (1) at (1,0) {};
	\node[wV, label=above:{$T_{\alpha_2^\vee}$}] (2) at (2,0) {};
	\node[wV, label=above:{$T_{\alpha_{n-2}^\vee}$}] (4) at (4,0) {};
	\node[wV, label=above:{$T_{\alpha_{n-1}^\vee}$}] (5) at (5,0) {};
	\end{tikzpicture}
 \qquad\qquad 
    \begin{tikzpicture}[every node/.style={inner sep=1}, scale=1]
	\node at (.5,0) {\tikz \draw[black, double distance=2] (0,0) -- +(.127,0);};
	\node at (5.5,0) {\tikz \draw[black, double distance=2] (0,0) -- +(-.127,0);};
	\foreach \x in {0,5}{\foreach \y in {-1,1}{\draw (\x,\y*.05) to (\x+1,\y*.05);}}
	\node[wV, label=above:{$T_{\alpha_0}$}] (0) at (0,0) {};
	\node[wV, label=above:{$T_{\alpha_n}$}] (6) at (6,0) {};
	\foreach \x/\y in {1/1,2/2,4/n\!-\!2,5/n\!-\!1}{
		\node[wV] (\x) at (\x,0) {};}
	\draw (1) to (2) (4) to (5);
	\draw[dashed] (2) to (4);
	\node[wV, label=above:{$T_{\alpha_1}$}] (1) at (1,0) {};
	\node[wV, label=above:{$T_{\alpha_2}$}] (2) at (2,0) {};
	\node[wV, label=above:{$T_{\alpha_{n-2}}$}] (4) at (4,0) {};
	\node[wV, label=above:{$T_{\alpha_{n-1}}$}] (5) at (5,0) {};
	\end{tikzpicture}
\label{B}
%\tag{B} 
\\[0.25cm]
T_{\alpha_0^\vee} T_{\alpha_1}^{-1} T_{\alpha_0}T_{\alpha_1} = T_{\alpha_1}^{-1}T_{\alpha_0}T_{\alpha_1}T_{\alpha_0^\vee},
\qquad\qquad
(T_0^\#-u_0^{\frac12})(T_0^\#+u_0^{\frac12}) = 0,
\label{HS}
%\tag{HS}
\\[0.25cm]
(T_{\alpha^\vee_i}-t_{\alpha^\vee_i}^{\frac12})(T_{\alpha^\vee_i}+t_{\alpha^\vee_i}^{-\frac12})=0
\quad\hbox{and}\quad
T_{\alpha_i^\vee} Y^{\lambda^\vee} 
= Y^{s_i^\vee\lambda^\vee} T_{\alpha_i^\vee} 
+ (c_{\alpha_i^\vee}^Y-t_{\alpha_i^\vee}^{-\frac12}) (Y^{\lambda^\vee} - Y^{s_i^\vee\lambda^\vee} ),
\label{HY}
%\tag{HY}
\\[0.25cm]
(T_{\alpha_i}-t_{\alpha_i}^{\frac12})(T_{\alpha_i}+t_{\alpha_i}^{-\frac12})=0
\quad\hbox{and}\quad
T_{\alpha_i} X^{-\mu}
= X^{-s_i\mu} T_{\alpha_i} +  (c_{\alpha_i}^{X^{-1}}-t_{\alpha_i}^{-\frac12}) (X^{-\mu} - X^{-s_i \mu} ).
\label{HX}
%\tag{HX}
\end{gather}
\end{prop}

\begin{remark} In personal communication,
J.\ Stokman insightfully points out that the algebra $\widetilde{H}_{\mathrm{loc}}$ defined
in~\eqref{Drels}--\eqref{Hbrg} is fishy. One can Ore-localize DAHA in either the $X$-elements
(so that the normalized intertwiners $\xi_{s_j}$ are in the localized algebra), or in the $Y$-elements
(so that the $\eta_{s_j}$ are in the ones in the localized algebra)  but not simultaneously.  
This difficulty is alluded to in a different form in~\cite[Remark 3.5]{CR22}.
However, trying to set up the accurate formal framework for handling an
$X$-$Y$-localized algebra would be distracting from the combinatorial perspective of this paper.
We feel that, even if the right localization formalism is not in place, the concept of an algebra $\widetilde{H}_{\mathrm{loc}}$ 
that contains all the useful relations for computations with these operators is a healthy point of view and so we have chosen to include it.
The algebra $\widetilde{H}_{\mathrm{loc}}$ is not absolutely necessary for the results in this paper as the operators on the polynomial representation that we use are well-defined and the proofs that we give are valid.
\end{remark}

\subsection{Symmetrizers}\label{ssec:symm}

There are four ways of symmetrizing/antisymmetrizing in the Koornwinder polynomial context, corresponding to the
four 1-dimensional representations of $W_{\mathrm{fin}}$.  These 1-dimensional representations 
are the analogues of the usual sign of a permutation.
The four symmetrizers, and useful formulas for them, are presented in Sections~\ref{ssec:symm},~\ref{ssec:cforsymm} 
and~\ref{ssec:symmwpar}.  The
symmetrizers will be used in Section~\ref{sec:BFM} to construct and manage the bosonic (symmetric), fermionic (antisymmetric)
and the two mesonic (half symmetric-half antisymmetric) versions of the Koornwinder polynomials.
These four symmetrized/antisymmetrized versions of Koornwinder polynomials, and the relations between them, turn out to be 
fundamental in the proof of the norm formulas and constant term formulas that are established in Section~\ref{sec:orth}.

The finite Hecke algebra $H_{\mathrm{fin}}$ is the $\KK$-subalgebra of $\widetilde H_{\mathrm{int}}$ generated by $T_1, \ldots, T_{n-1}$ and $T_n$.  
The finite Hecke algebra
\begin{equation}
H_{\mathrm{fin}} \quad\hbox{has $\KK$-basis}\quad \{T_v\ |\ v\in W_{\mathrm{fin}}\},
\label{Hfinbasis}
%\tag{Hfinbasis}
\end{equation}
where $T_v = T_{i_1}\ldots T_{i_k}$ if $v = s_{i_1}\cdots s_{i_k}$ is a reduced word for $v$ in $W_{\mathrm{fin}}$.
The four one dimensional representations of $H_{\mathrm{fin}}$ are 
$$\chi^+\colon H_{\mathrm{fin}}\to \KK,
\qquad
\chi^\pm\colon H_{\mathrm{fin}}\to \KK,
\qquad
\chi^\mp\colon H_{\mathrm{fin}}\to \KK,
\qquad
\chi^-\colon H_{\mathrm{fin}}\to \KK,
$$
given by 
\begin{align*}
\chi^+(T_i) &= \begin{cases}
t^{\frac12}, &\hbox{if $i\in \{1, \ldots, n-1\}$,} \\
t_n^{\frac12}, &\hbox{if $i=n$,}
\end{cases}
\qquad
&\chi^\pm(T_i) &= \begin{cases}
t^{\frac12}, &\hbox{if $i\in \{1, \ldots, n-1\}$,} \\
(-t_n)^{-\frac12}, &\hbox{if $i=n$,}
\end{cases}
\\
\chi^\mp(T_i) &= \begin{cases}
(-t)^{-\frac12}, &\hbox{if $i\in \{1, \ldots, n-1\}$,} \\
t_n^{\frac12}, &\hbox{if $i=n$,}
\end{cases}
\qquad
&\chi^-(T_i) &= \begin{cases}
(-t)^{-\frac12}, &\hbox{if $i\in \{1, \ldots, n-1\}$,} \\
(-t_n)^{-\frac12}, &\hbox{if $i=n$.}
\end{cases}
\end{align*}

For $v\in W_{\mathrm{fin}}$ define $\ell_s(v)$ and $\ell_d(v)$ by
\begin{align*}
\chi^+(T_v) &= (t^{\frac12})^{\ell_s(v)} (t_n^{\frac12})^{\ell_d(v)}.
\end{align*}

The \emph{Hecke symmetrizers} are 
$$\hbox{the elements\quad $\varepsilon_+$,  $\varepsilon_\pm$,  $\varepsilon_\mp$ and $\varepsilon_-$ \quad
of $H_{\mathrm{fin}}$}
$$
which are defined such that, in terms of the basis in~\eqref{Hfinbasis}, the coefficient of $T_{w_0}$ is $1$ and for $w\in W_{\mathrm{fin}}$,
$$T_w \varepsilon_+ = \chi^+(T_w)\varepsilon_+,
\qquad
T_w \varepsilon_{\pm}= \chi^\pm(T_w)\varepsilon_{\pm},
\qquad
T_w \varepsilon_{\mp}= \chi^\mp(T_w)\varepsilon_{\mp},
\qquad
T_w \varepsilon_- = \chi^-(T_w)\varepsilon_-.
$$
In other words, if $\Xi\in \{+,\pm, \mp, -\}$ then $T_w \varepsilon_\Xi = \chi^\Xi(T_w)\varepsilon_\Xi$.

A reduced word for the longest element of $W_{\mathrm{fin}}$ is
$$w_0 = (s_1\cdots s_n\cdots s_1)(s_2\cdots s_n\cdots s_2)\cdots (s_{n-1}s_ns_{n-1})s_n
\qquad\hbox{and}\qquad
t^{\frac12\ell_s(w_0)}t_n^{\frac12\ell_d(w_0)} = t^{\frac12n(n-1)} t_n^{\frac12 n}.
$$
In terms of the basis in~\eqref{Hfinbasis}
the  symmetrizers  are given explicitly by
\begin{equation}
    \varepsilon_\Xi = \frac{1}{\chi^\Xi(T_{w_0})} \sum_{v\in W_{\mathrm{fin}}} \chi^{\Xi}(T_v) T_v,
\qquad\hbox{for $\Xi\in \{+,\pm, \mp, -\}$.}
\label{symme-Hfin}
%\tag{symme-Hfin}
\end{equation}

The \emph{Poincar\'e polynomial for $W_{\mathrm{fin}}$} is
\begin{equation}
W_0(t,t_n) = \sum_{w\in W_{\mathrm{fin}}} t^{\ell_s(w)} t_n^{\ell_d(w)}
= \sum_{w\in W_{\mathrm{fin}}} \chi^+(T_w)^2.
\label{Poindef}
%\tag{Poindef}
\end{equation}
Three alternate formulas for $W_0(t,t_n)$ are given in Proposition~\ref{Poinprop}.
Then
\begin{align}
\varepsilon_+^2 &= \frac{1}{\chi^+(T_{w_0})} W_0(t,t_n)\varepsilon_+, \qquad
&\varepsilon_\pm^2 &= \frac{1}{\chi^\pm(T_{w_0})} W_0(t,t_n^{-1})\varepsilon_\pm, 
\nonumber \\
\varepsilon_-^2 &= \frac{1}{\chi^-(T_{w_0})} W_0(t^{-1},t_n^{-1})\varepsilon_-,
&\varepsilon_\mp^2 &= \frac{1}{\chi^\mp(T_{w_0})} W_0(t^{-1},t_n)\varepsilon_\mp. 
\label{symmsq}
%\tag{symmsq}
\end{align}

\begin{example}
For $\Xi=\pm$, 
\begin{align*}
\varepsilon_\pm &=\frac{1}{\chi^\pm(T_{w_0})}\sum_{w\in W_{\mathrm{fin}}} \chi^\pm(T_w)T_w 
=  t^{-\frac12\ell_s(w_0)} (-t_n^{\frac12})^{\ell_d(w_0)}
\sum_{w\in W_{\mathrm{fin}}} t^{\frac12 \ell_s(w)}(-t_n^{-\frac12})^{\ell_d(w)} T_w,
\end{align*}
and since $T_w \varepsilon_\pm = \chi^{\pm}(T_w)\varepsilon_\pm$ then
$$\varepsilon_\pm^2 = \frac{1}{\chi^\pm(T_{w_0})}  \sum_{w\in W_{\mathrm{fin}}} \chi^\pm(T_w)^2\varepsilon_\pm
= \frac{1}{\chi^\pm(T_{w_0})} W_0(t,t_n^{-1})\varepsilon_\pm.
$$
\end{example}

\subsection{$c$-function formulas for symmetrizers}\label{ssec:cforsymm}

The definition of the Hecke symmetrizers and the formulas for them given in Section~\ref{ssec:symm},
are purely in terms of the $T_w$ in the double affine Hecke algebra $\widetilde{H}_{\mathrm{int}}$.  However the
Koornwinder polynomials are more naturally constructed and managed with the $\eta_w$ and $\xi_w$ that are
in $\widetilde{H}_{\mathrm{loc}}$ and so it becomes desirable to have expressions for the symmetrizers
that are in terms of the $\eta_w$ and $\xi_w$ (and $c$-functions).  Perhaps surprisingly, these conversion formulas,
presented in Proposition~\ref{cfcnsymm},
are compact and elegant (and useful!).  We will use them in Section~\ref{sec:BFM}
to provide formulas for the symmetric (bosonic), antisymmetric (fermionic) and half symmetric-half antisymmetric (mesonic)
Koornwinder polynomials.

For $w\in W_{\mathrm{fin}}$ let
$$\xi_w = \xi_{s_{i_1}}\cdots \xi_{s_{i_\ell}}
\qquad\hbox{and}\qquad
\eta_w = \eta_{s_{i_1}^\vee}\cdots \eta_{s_{i_\ell}^\vee},
$$
if $w=s_{i_1}\cdots s_{i_\ell}$ is a reduced word for $w$.  There are four $X$-symmetrizers
\begin{align*}
e_+^X &= \sum_{w\in W_{\mathrm{fin}}} \xi_w, \qquad
&e_-^X &= \sum_{w\in W_{\mathrm{fin}}} (-1)^{\ell_s(w)+\ell_d(w)} \xi_w, \\
e_\pm^X &= \sum_{w\in W_{\mathrm{fin}}} (-1)^{\ell_d(w)} \xi_w, \qquad
&e_\mp^X &= \sum_{w\in W_{\mathrm{fin}}} (-1)^{\ell_s(w)} \xi_w,
\end{align*}
and four $Y$-symmetrizers
\begin{align*}
e_+^Y &= \sum_{w\in W_{\mathrm{fin}}} \eta_w, 
&e_-^Y &= \sum_{w\in W_{\mathrm{fin}}} (-1)^{\ell_s(w)+\ell_d(w)} \eta_w, \\
e_\pm^Y &= \sum_{w\in W_{\mathrm{fin}}} (-1)^{\ell_d(w)} \eta_w,
&e_\mp^Y &= \sum_{w\in W_{\mathrm{fin}}} (-1)^{\ell_s(w)} \eta_w.
\end{align*}

The following proposition writes the Hecke symmetrizers in terms of the $X$-symmetrizers and the $Y$-symmetrizers.

\begin{prop} \label{cfcnsymm} 
Let $\kappa_{\mathrm{st}}^X$, $\kappa_{\mathrm{dr}}^X$ be as defined in~\eqref{cstcdr}.
\begin{equation*}
\begin{array}{ccc}
\chi^+(T_{w_0})\varepsilon_+ 
= e_+^X \kappa_{\mathrm{st}}^{X^{-1}}\kappa_{\mathrm{dr}}^{X^{-1}} = e_+^Y \kappa_{\mathrm{st}}^Y \kappa_{\mathrm{dr}}^Y
&  \quad\hbox{and}\quad &
\chi^+(T_{w_0})\varepsilon_- 
= \kappa_{\mathrm{st}}^X \kappa_{\mathrm{dr}}^X e_-^X = \kappa_{\mathrm{st}}^{Y^{-1}}\kappa_{\mathrm{dr}}^{Y^{-1}} e_-^Y, \\[0.2cm]
\chi^+(T_{w_0})\varepsilon_\pm 
= \kappa_{\mathrm{dr}}^X e_{\pm}^X \kappa_{\mathrm{st}}^{X^{-1}} = \kappa_{\mathrm{dr}}^{Y^{-1}} e^Y_{\pm} \kappa_{\mathrm{st}}^Y
& \quad\hbox{and}\quad &
\chi^+(T_{w_0})\varepsilon_\mp 
= \kappa_{\mathrm{st}}^X e_{\mp}^X \kappa_{\mathrm{dr}}^{X^{-1}} = \kappa_{\mathrm{st}}^{Y^{-1}} e^Y_{\mp} \kappa_{\mathrm{dr}}^Y.
\end{array}
\end{equation*}
\end{prop}
\begin{proof}
Let us prove the formula $\chi^+(T_{w_0})\varepsilon_{\pm} = \kappa_{\mathrm{dr}}^X e_{\pm}^X \kappa_{\mathrm{st}}^{X^{-1}}.$
The proof for the other cases is similar.

Let 
$$R_{\pm} = \kappa_{\mathrm{dr}}^X e_{\pm}^X \kappa_{\mathrm{st}}^{X^{-1}}.$$
For $i\in \{1, \ldots, n-1\}$, and using~\eqref{Talt}, 
\begin{align*}
(T_i-t^{\frac12})R_\pm 
&= (T_i-t^{\frac12})\kappa_{\mathrm{dr}}^X e_{\pm}^X \kappa_{\mathrm{st}}^{X^{-1}} 
= - c_{\alpha_i}^{X}(1-\xi_{s_i})\kappa_{\mathrm{dr}}^X e_{\pm}^X \kappa_{\mathrm{st}}^{X^{-1}}
\\
&= - \kappa_{\mathrm{dr}}^X  c_{\alpha_i}^{X}(1-\xi_{s_i})e_{\pm}^X \kappa_{\mathrm{st}}^{X^{-1}}
= - \kappa_{\mathrm{dr}}^X  c_{\alpha_i}^{X}\cdot 0\cdot  \kappa_{\mathrm{st}}^{X^{-1}} =0,
\qquad\hbox{so that}\quad
T_i R_{\pm} = t^{\frac12}R_{\pm}.
\end{align*}
Using~\eqref{TdefH},
\begin{align*}
(T_n+t_n^{-\frac12})R_\pm
&= (T_n+t_n^{-\frac12})\kappa_{\mathrm{dr}}^X e_{\pm}^X \kappa_{\mathrm{st}}^{X^{-1}} 
= (1+\xi_{s_n}) c_{\alpha_n}^{X^{-1}} \kappa_{\mathrm{dr}}^X  e_{\pm}^X \kappa_{\mathrm{st}}^{X^{-1}} 
= (1+\xi_{s_n}) \kappa_{\mathrm{dr}}^X c_{\alpha_n}^{X^{-1}} e_{\pm}^X \kappa_{\mathrm{st}}^{X^{-1}} 
\\
&= \kappa_{\mathrm{dr}}^X (1+\frac{c_{\alpha_n}^{X^{-1}}}{c_{\alpha_n}^X} \xi_{s_n})c_{\alpha_n}^{X^{-1}}  e_{\pm}^X \kappa_{\mathrm{st}}^{X^{-1}} 
= \kappa_{\mathrm{dr}}^X \left(c_{\alpha_n}^{X^{-1}}+\frac{c_{\alpha_n}^{X^{-1}}}{c_{\alpha_n}^X} c_{\alpha_n}^{X}\xi_{s_n}\right)  e_{\pm}^X \kappa_{\mathrm{st}}^{X^{-1}} 
\\
&= \kappa_{\mathrm{dr}}^X \left(c_{\alpha_n}^{X^{-1}}+c_{\alpha_n}^{X^{-1}}\xi_{s_n}\right)  e_{\pm}^X \kappa_{\mathrm{st}}^{X^{-1}} 
= \kappa_{\mathrm{dr}}^X \left(c_{\alpha_n}^{X^{-1}}-c_{\alpha_n}^{X^{-1}}\right)  e_{\pm}^X \kappa_{\mathrm{st}}^{X^{-1}} 
=0,
\end{align*}
so that $T_n R_{\pm} = - t_n^{-\frac12}R_{\pm}.$

Since  $T_i = \xi_{s_i}c_{\alpha_i}^{X^{-1}} + (c_{\alpha_i}^{X^{-1}}-t^{-\frac12})$
and the coefficient of $T_{w_0}$ is $1$ then there are rational functions $a_w^{X^{-1}}$ such that
\begin{align*}
\chi^+(T_{w_0}) \varepsilon_{\pm} 
&= \chi^+(T_{w_0}) \sum_{w\in W_{\mathrm{fin}}} \xi_w a_w^{X^{-1}}
= \chi^+(T_{w_0})\xi_{w_0}c_{w_0}^{X^{-1}} + \chi^+(T_{w_0})\sum_{w<w_0} \xi_w a_w^{X^{-1}}
\\
&= \kappa_{\mathrm{dr}}^{X}\xi_{w_0} \kappa_{\mathrm{st}}^{X^{-1}}  + \chi^+(T_{w_0}) \sum_{w<w_0} \xi_w a_w^{X^{-1}}.
\end{align*}
The element $\varepsilon_{\pm}\in H_{\mathrm{fin}}$ is determined by the conditions that  the coefficient of $T_{w_0}$ is $1$, 
and $T_n\varepsilon_{\pm} = -t_n^{-\frac12}\varepsilon_{\pm}$ 
and $T_i \varepsilon_{\pm} = t^{\frac12}\varepsilon_{\pm}$ for $i\in \{1, \ldots, n-1\}$.
So $\chi^+(T_{w_0})\varepsilon_{\pm} = R_{\pm}$.
\end{proof}

\subsection{Symmetrizers and stabilizers}\label{ssec:symmwpar}

The finite Weyl group $W_{\mathrm{fin}}$ acts on $\ZZ^n$ by the formulas in~\eqref{Wfinaction}.
Since the action of $W_{\mathrm{fin}}$ on $\ZZ^n$ is not free, there are elements of $\ZZ^n$ that have nontrivial
stabilizer, and one is forced to confront these stabilizers.  This subsection computes formulas for the symmetrizers
which take into account, and allow us to manage, the cases when the stabilizer is nontrivial.

Let
$$(\ZZ_{\ge0}^n)^+ = \{
\lambda=(\lambda_1, \ldots, \lambda_n)\in \ZZ^n \ |\ \lambda_1\ge \cdots \ge \lambda_n\ge 0\}.
$$
The set $(\ZZ_{\ge 0}^n)^+$ is a set of representatives of the $W_{\mathrm{fin}}$-orbits on $\ZZ^n$.
For $\lambda\in (\ZZ^n_{\ge0})^+$, let
$$W_\lambda = \{ w\in W_{\mathrm{fin}}\ |\ w\lambda = \lambda\}
\quad\hbox{and}\quad
W^\lambda = \left\{ \begin{array}{c} \hbox{minimal length representatives} \\
\hbox{of cosets in $W_{\mathrm{fin}}/W_\lambda$} \end{array} \right\}.
$$
Let $w_\lambda$ be the longest element of $W_\lambda$ and let $v_\lambda$ be the maximal element of $W^\lambda$.
As in~\eqref{umuvmudefn},  then $v_\lambda$ is the minimal length element of $W_{\mathrm{fin}}$ such that $v_\lambda \lambda$ is increasing with all entries $\le 0$.
Let
\begin{equation}
\varepsilon^+_\lambda = \frac{1}{\chi^+(T_{w_\lambda})} \sum_{w\in W_\lambda} \chi^+(T_w)T_w
\label{parsymm}
%\tag{parsymm}
\end{equation}
so that
$T_w \varepsilon^+_\lambda = \chi^+(T_w) \varepsilon^+_\lambda$,
for $w\in W_\lambda$, and the coefficient of $T_{w_\lambda}$ in $\varepsilon^+_\lambda$ is $1$.

Define $\rho, \omega, \pi\in (\ZZ_{\ge 0}^n)^+$ by 
\begin{equation}
\rho = (n,n-1, \ldots, 2,1), \qquad
\omega = (1,1,\ldots, 1), \qquad
\pi = (n-1, \ldots, 2,1,0).
\label{rhodef}
%\tag{rhodef}
\end{equation}
The statement of the following proposition 
is designed to stress the analogies between the four symmetrizers $\varepsilon_+$, $\varepsilon_\pm$,
$\varepsilon_\mp$ and $\varepsilon_-$.
As in Remark~\eqref{smstabs}, in practice, there are simplifications since the stabilizers $W_{\lambda+\rho}$ and
$W_{\lambda+\pi}$ have order 1 or 2.  Specifically,
$$
W_{\lambda+\rho}=1, \qquad
W^{\lambda+\rho} = W_{\mathrm{fin}}, \quad \chi^+(T_{w_{\lambda+\rho}}) = 1 \quad\hbox{and}\quad
\varepsilon^+_{\lambda+\rho} = 1;
$$
$$\hbox{if $\lambda_n=0$ then}\qquad
W_{\lambda+\pi} = \{1, s_n\}, \quad \chi^+(T_{w_{\lambda+\pi}}) = t_n^{\frac12} \quad\hbox{and}\qquad
\varepsilon^+_{\lambda+\pi} = T_n+t_n^{-\frac12};$$
and 
$$\hbox{if $\lambda_n\ne0$ then}\qquad
W_{\lambda+\pi} = 1, \quad
W^{\lambda+\pi} = W_{\mathrm{fin}},\quad 
\chi^+(T_{w_{\lambda+\pi}}) = 1 \quad\hbox{and}\qquad
\varepsilon^+_{\lambda+\pi}=1.
$$

\begin{prop} \label{symmexp}
Let $\lambda\in (\ZZ_{\ge 0}^n)^+$.
Then 
\begin{align*}
\chi^+(T_{w_0})\varepsilon_+ &= \chi^+(w_\lambda)
\left(\sum_{z\in W^\lambda} \kappa^+_{v_\lambda z} \eta_z \kappa^Y_z\right)\varepsilon^+_\lambda, 
\\
\chi^+(T_{w_0})\varepsilon_\pm &= \chi^+(w_{\lambda+\omega})
\left(\sum_{z\in W^{\lambda+\omega}} 
(-1)^{\ell_d(z)}\kappa^\pm_{v_{\lambda+\omega} z} \eta_z \kappa^Y_z\right)\varepsilon^+_{\lambda+\omega}, 
\\
\chi^+(T_{w_0})\varepsilon_\mp &= \chi^+(w_{\lambda+\pi})
\left(\sum_{z\in W^{\lambda+\pi} }
(-1)^{\ell_s(z)}\kappa^\mp_{v_{\lambda+\pi} z} \eta_z \kappa^Y_z\right)\varepsilon^+_{\lambda+\pi}, 
\\
\chi^+(T_{w_0})\varepsilon_- &= \chi^+(w_{\lambda+\rho})
\left(\sum_{z\in W^{\lambda+\rho}} 
(-1)^{\ell(z)} \kappa^-_{v_{\lambda+\rho} z} \eta_z \kappa^Y_z\right)\varepsilon^+_{\lambda+\rho}.
\end{align*}
\end{prop}

\begin{proof} We will prove the $\pm$ case.  The proof for the other cases is similar.

Let 
$$J_{\lambda+\omega} = \{ \alpha\in (S^\vee)_{0,s}^+\ |\ \alpha^\vee\not\in \mathrm{Inv}(w_{\lambda+\omega})\}$$
so that $J_{\lambda+\omega} $ is the complement of $\mathrm{Inv}(w_{\lambda+\omega})$ in $S_{0,s}^+$.
If $v\in W_{\lambda+\omega}$ then
$vJ_{\lambda+\omega} = J_{\lambda+\omega}$,
since $v$ permutes the elements of $J_{\lambda+\omega}$.  
$$\hbox{If $z\in W^{\lambda+\omega}$\quad then\quad} 
\kappa_{\mathrm{dr}}^{Y^{-1}} \eta_z \kappa_{J_{\lambda+\omega}}^Y =  \kappa^\pm_{v_{\lambda+\omega}z}\eta_z \kappa_z^Y.
$$
Then
\begin{align*}
\chi^+(T_{w_0})\varepsilon_\pm 
= \kappa_{\mathrm{dr}}^{Y^{-1}} e_\pm^Y \kappa_{\mathrm{st}}^Y
&= \kappa_{\mathrm{dr}}^{Y^{-1}} \left( \sum_{w\in W_{\mathrm{fin}}} (-1)^{\ell_d(w)} \eta_w\right) \kappa_{\mathrm{st}}^Y
\\
&= \kappa_{\mathrm{dr}}^{Y^{-1}} \left( 
\sum_{z\in W^{\lambda+\omega}} \sum_{v\in W_{\lambda+\omega}}  (-1)^{\ell_d(z)} \eta_z \eta_v\right) \kappa_{\mathrm{st}}^Y
\\ 
&= \kappa_{\mathrm{dr}}^{Y^{-1}} \left( 
\sum_{z\in W^{\lambda+\omega}}   (-1)^{\ell_d(z)} \eta_z\right)\left(\sum_{v\in W_{\lambda+\omega}} \eta_v\right) 
\kappa_{J_{\lambda+\omega}}^Y \kappa_{w_{\lambda+\omega}}^Y
\\
&= \kappa_{\mathrm{dr}}^{Y^{-1}} \left( 
\sum_{z\in W^{\lambda+\omega}}   (-1)^{\ell_d(z)} \eta_z\right)
\kappa_{J_{\lambda+\omega}}^Y
\left(\sum_{v\in W_{\lambda+\omega}} \eta_v\right)  \kappa_{w_{\lambda+\omega}}^Y
\\
&=  \left( 
\sum_{z\in W^{\lambda+\omega}}   (-1)^{\ell_d(z)} \kappa^\pm_{v_{\lambda+\omega}z} \eta_z \kappa_z^Y\right)
\left(\sum_{v\in W_{\lambda+\omega}} \eta_v\right)  \kappa_{w_{\lambda+\omega}}^Y
\\
&=  \left( 
\sum_{z\in W^{\lambda+\omega}}   (-1)^{\ell_d(z)} \kappa^\pm_{v_{\lambda+\omega}z} \eta_z\kappa_z^Y\right)
\chi^+(T_{w_{\lambda+\omega}})\varepsilon^+_{\lambda+\omega}.
\end{align*}
\end{proof}

\begin{remark} \label{rem:bigvmu}
With $\omega = (1,1,\ldots, 1)$ and $\pi = (n-1, \ldots, 2,1,0)$ and $v_\omega$ and $v_\pi$ as defined
in~\eqref{umuvmudefn} then
$$v_\omega(i) = -(n-i+1)\quad\hbox{so that}\quad
\mathrm{Inv}_s(v_\omega) = \{ \varepsilon_i+\varepsilon_j\ |\ i<j\} 
\quad\hbox{and}\quad
\mathrm{Inv}_d(v_\omega) = \{\varepsilon_1, \ldots, \varepsilon_n\}
$$
since, for example, $v_\omega(\varepsilon_1-\varepsilon_n) = -\varepsilon_n+\varepsilon_1 = \varepsilon_1-\varepsilon_n$.
Then
$$v_\pi(i) = \begin{cases} -i, &\hbox{if $i\ne n$,} \\ n, &\hbox{if $i=n$,} \end{cases}
\quad\hbox{so that}\quad
\mathrm{Inv}_s(v_\pi) = \{ \varepsilon_i\pm \varepsilon_j\}
\quad\hbox{and}\quad
\mathrm{Inv}_d(v_\pi) = \{ \varepsilon_1, \ldots, \varepsilon_{n-1}\}.
$$
If $\lambda\in (\ZZ_{\ge 0}^n)^+$ then 
$W_{\lambda+\omega}\subseteq W_\omega$ and $W_{\lambda+\pi}\subseteq W_\pi$
giving 
$$\mathrm{Inv}_d(v_{\lambda+\omega})\supseteq \mathrm{Inv}_d(v_\omega)
\qquad\hbox{and}\qquad
\mathrm{Inv}_s(v_{\lambda+\pi})\supseteq \mathrm{Inv}_s(v_\pi).$$
Thus
$$\mathrm{Inv}_d(v_{\lambda+\omega}) = \{\varepsilon_1, \ldots, \varepsilon_n\}
\qquad\hbox{and}\qquad
\mathrm{Inv}_s(v_{\lambda+\pi}) = \{\varepsilon_i\pm \varepsilon_j\ |\ i<j\}.$$
\qed\end{remark}

\section{Electronic Macdonald-Koornwinder polynomials}\label{sec: electronic MK}

The electronic Koornwinder polynomials $E_\mu$ form a basis of the polynomial ring
$\KK[x_1^{\pm1},\ldots, x_n^{\pm1}]$.  They are simultaneous eigenvectors for 
the Cherednik-Dunkl operators $Y_1, \ldots, Y_n$.  This is in analogy to the way that, in quantum mechanics, 
Hermite polynomials are eigenfunctions of a Hamiltonian operator.  In this section we set up the operators
$Y_1, \ldots, Y_n$ on the polynomial ring, characterize the electronic Macdonald polynomials $E_\mu$ as eigenvectors, and provide recursive formulas for
computing them.

\subsection{Operators on polynomials}

Let $\KK[x] = \KK[x_1^{\pm1}, \ldots, x_n^{\pm1}]$.  
Let $\ZZ^n$ denote the set of length $n$ sequences $\mu = (\mu_1, \ldots, \mu_n)$ of integers.
The ring 
$$\hbox{$\KK[x]$ has basis}\quad
\{ x^\mu\ |\ \mu = (\mu_1 ,\ldots, \mu_n) \in \ZZ^n\},
\qquad\hbox{where
\quad
$x^\mu = x_1^{\mu_1}\cdots x_n^{\mu_n}$.}
$$
 Define operators $\xi_{s_0}, \xi_{s_1}, \ldots, \xi_{s_n}$ on $\KK[x^{\pm1}_1, \ldots, x^{\pm1}_n]$ by
\begin{align*}
(\xi_{s_0}f)(x_1, \ldots, x_n) &= f(qx^{-1}_1, x_2, \ldots, x_n), \\
(\xi_{s_i}f)(x_1, \ldots, x_n) &= f(x_1, \ldots, x_{i-1}, x_{i+1}, x_i, x_{i+2}, \ldots, x_n), \\
(\xi_{s_n}f)(x_1, \ldots, x_n) &= f(x_1, \ldots, x_{n-1}, x^{-1}_n), 
\end{align*}
for $f\in \KK[x^{\pm1}_1, \ldots, x^{\pm1}_n]$ and $i\in \{1, \ldots, n-1\}$.  
Define operators $X_1, \ldots, X_n$ on $\KK[x_1^{\pm1}, \ldots, x_n^{\pm1}]$ by
\begin{equation}
X_j f = x_j f, \qquad \hbox{for  $j\in \{1, \ldots, n\}$.} 
\label{polyaction}
%\tag{polyaction}
\end{equation}

Consider the induced representation 
$$\Ind_{H_Y}^{\widetilde H}(\mathbf{1}_Y) = \widetilde H_{\mathrm{int}}\mathbf{1}_Y 
=\hbox{$\KK$-span}\{ X_1^{\mu_1}\cdots X_n^{\mu_n}\mathbf{1}_Y\ |\ \mu = (\mu_1, \ldots, \mu_n)\in \ZZ^n\}
$$
determined by 
$$T_{\alpha_0}\mathbf{1}_Y = t_0^{\frac12} \mathbf{1}_Y, \qquad\hbox{and}\qquad
T_i\mathbf{1}_Y = t^{\frac12}\mathbf{1}_Y, \qquad\hbox{for $i\in \{1, \ldots, n\}$.}
$$
Then the map
\begin{equation}
\begin{matrix}
\KK[x_1^{\pm1}, \ldots, x_n^{\pm1}] &\longrightarrow &\widetilde H_{\mathrm{int}} \mathbf{1}_Y \\
x_1^{\mu_1}\cdots x_n^{\mu_n} &\longmapsto &X_1^{\mu_1}\cdots X_n^{\mu_n}\mathbf{1}_Y
\end{matrix}
\qquad\hbox{is an $\widetilde H_{\mathrm{int}}$-module isomorphism.}
\label{indrep}
%\tag{indrep}
\end{equation}
We shall often identify $\KK[x] = \KK[x_1^{\pm1}, \ldots, x_n^{\pm1}]$ and $\widetilde H_{\mathrm{int}} \mathbf{1}_Y$ and
$\KK[X] = \KK[X_1^{\pm1},\ldots, X_n^{\pm1}]$ via this isomorphism.

\subsection{The operators $T_{\alpha_0}, \ldots, T_{\alpha_n}$ and $Y_1, \ldots, Y_n$}

Define operators $T_{\alpha_0}, T_{\alpha_1}, \ldots, T_{\alpha_n}$ by
\begin{equation}
T_{\alpha_i} = t_{\alpha_i}^{\frac12} - c_{\alpha_i}^X(1-\xi_{s_i}),
\qquad\hbox{for $i\in \{0,1,\ldots, n\}$.}
\label{TiopsinKX}
%\tag{TiopsinKX}
\end{equation}
Define
$$Y_j = T_{\alpha_{j-1}}^{-1}\cdots T_{\alpha_1}^{-1}T_{\alpha_0}T_{\alpha_1}\cdots T_{\alpha_n}\cdots T_{\alpha_j},
\qquad\hbox{for $j\in \{1, \ldots, n\}$.}
$$

Using~\eqref{rtparams} and~\eqref{Xsroots},
\begin{align}
t_0^{\frac12}T_{\alpha_0} &= t_0  - \left( \frac{(x_1+q^{\frac12}t^{\frac12}_0 u^{-\frac12}_0)
(x_1-q^{\frac12}t^{\frac12}_0u^{\frac12}_0)}
{x_1^2-q}\right) (1-\xi_{s_0}),
\nonumber \\
t_n^{\frac12}T_{\alpha_n} &= t_n  - \left( \frac{(1+t^{\frac12}_n u^{-\frac12}_n x_n)
(1- t^{\frac12}_n u^{\frac12}_n x_n)}
{1-x_n^2}\right) (1-\xi_{s_n}),
\label{DGaction}
%\tag{DGaction}
\\
t^{\frac12}T_{\alpha_i} &= t - \frac{x_{i+1}-tx_i}{x_{i+1}-x_i} (1-\xi_{s_i}),
\quad\hbox{for $i\in \{1, \ldots, n-1\}$}
\nonumber
\end{align}
(see~\cite[\S3]{Nou95} and~\cite[(13)]{Sah99} and~\cite[(73)]{CGdGW16}).

\subsection{Electronic Macdonald polynomials $E_\mu$}

For $\mu\in \ZZ^n$ let $v_\mu\in W_{\mathrm{fin}}$ be the minimal length signed permutation such that $v_\mu \mu$ is weakly increasing
with all entries $\le 0$.
The \emph{electronic Macdonald polynomials}
$$E_\mu(x_1,\ldots, x_n;q,t, t_0^{\frac12}, u_0^{\frac12}, t_n^{\frac12}, u_n^{\frac12})
\in \KK[x_1^{\pm1},\ldots, x_n^{\pm1}]
\qquad\hbox{are indexed by $\mu=(\mu_1, \ldots, \mu_n)\in \ZZ^n$}
$$
and are defined by the eigenvalue conditions
\begin{equation}
Y_j E_\mu = q^{-\mu_j}t^{-v_\mu(j)} (t^{\frac12}_0t^{\frac12}_n t^n)^{\mathrm{sgn}(v_\mu(j))}
E_\mu,
\label{Eeig}
%\tag{Eeig}
\end{equation}
for $\mu=(\mu_1, \ldots, \mu_n)\in \ZZ^n$ and $j\in \{1, \ldots, n\}$.
The normalization of $E_\mu$ is such that the coefficient of $x^\mu$ in $E_\mu$ is 1.

Let $\KK(Y)$ be the field of fractions of $\KK(Y)$.  For $\mu=(\mu_1, \ldots, \mu_n)\in \ZZ^n$ define homomorphisms
$\ev_\mu^t\colon \KK[Y]\to \KK[Y]$ by
\begin{equation}
\ev_\mu^t(Y_i) = q^{-\mu_i}t^{-v_\mu(j)}(t_0^{\frac12}t_n^{\frac12}t^n)^{\mathrm{sgn}(v_\mu(j))},
\qquad\hbox{for $i\in \{1, \ldots, n\}$.}
\label{evhom}
%\tag{evhom}
\end{equation}
Extend $\ev^t_\mu$ to those elements of the field $\KK(Y)$ for which the evaluated denominator is nonzero.  By~\eqref{Eeig}
\begin{equation}
f E_\mu = \ev^t_\mu(f) E_\mu, \qquad\hbox{for $f\in \KK[Y]$ and $\mu\in \ZZ^n$.}
\label{eveg}
%\tag{eveg}
\end{equation}

\subsection{The recursion for the $E_\mu$}

Although the eigenvalue conditions together with the normalization  completely characterize the Koornwinder polynomials
$E_\mu$, computing them by solving directly for eigenvectors is not efficient.  Fortunately, the operators
$\tau^\vee_i = \eta_{s_i^\vee} c^Y_{\alpha^\vee_i}$ from the algebra $\widetilde{H}_{\mathrm{loc}}$ provide a very nice
recursive way of computing the $E_\mu$.  This is analogous to the way that, in Schubert calculus, the Schubert polynomials
are constructed recursively using divided-difference operators.

Define operators $T_{\alpha_0^\vee}, T_{\alpha_1^\vee}, \ldots, T_{\alpha_n^\vee}$ on $\KK[X]$ by
\begin{equation}
(T_{\alpha_0^\vee})^{-1} = X_1T_{\alpha_1}\cdots T_{\alpha_n}\cdots T_{\alpha_1}
\qquad\hbox{and}\qquad
T_{\alpha^\vee_i} = T_{\alpha_i}
\qquad\hbox{for $i\in \{1, \ldots, n\}$.}
\label{TiveeonKX}
%\tag{TiveeonKX}
\end{equation}
and define
\begin{equation}
\tau^\vee_i = T_{\alpha_i^\vee}+(t_{\alpha_i^\vee}^{-\frac12} - c_{\alpha_i^\vee}^Y )
= (T_{\alpha_i^\vee})^{-1} +  (t_{\alpha_i^\vee}^{\frac12} - c_{\alpha_i^\vee}^Y),
\qquad\hbox{for $i\in \{0, 1, \ldots, n\}$.}
\label{tauvee}
%\tag{tauvee}
\end{equation}
By~\eqref{TdefH}, $\tau_i^\vee = \eta^Y_{s_i}c_{\alpha_i^\vee}^Y$ so that
\begin{equation}
\tau^\vee_i Y^{\lambda^\vee} = Y^{s_i^\vee\lambda^\vee} \tau^\vee_i,
\qquad\hbox{for $i\in \{0,1,\ldots, n\}$.}
\label{taupastY}
%\tag{taupastY}
\end{equation}
The group $W_X$ (generated by $s^\vee_0, \ldots, s^\vee_n$)
acts on $\ZZ^n$ by 
\begin{align*}
s_0^\vee(\mu_1, \ldots, \mu_n) &= (-\mu_1+1, \mu_2, \ldots, \mu_n), \\
s_i^\vee(\mu_1, \ldots, \mu_n) &= (\mu_1, \ldots, \mu_{i-1}, \mu_{i+1}, \mu_i, \mu_{i+2}, \ldots, \mu_n), 
\qquad\hbox{for $i\in \{1, \ldots, n-1\}$, and } \\
s_n^\vee(\mu_1, \ldots, \mu_n) &= (\mu_1, \ldots, \mu_{n-1}, -\mu_n).
\end{align*}
The relation~\eqref{taupastY} is the reason that 
the electronic Macdonald polynomials $E_\mu$
are equivalently defined by the following recursive relations:
\begin{enumerate}
\item[(E0)] $E_{(0,\ldots, 0)} = 1$;
\item[(E1)] if $\mu_1 \le 0$ then $E_{s_0^\vee\mu} = t^{n-1}t_n^{\frac12} \tau^\vee_0 E_\mu$;
\item[(E2)] if $i\in \{1, \ldots, n-1\}$ and $\mu_i>\mu_{i+1}$ then
$E_{s^\vee_i\mu} = t^{\frac12} \tau^\vee_i E_\mu$; and
\item[(E3)] 
if $\mu_n>0$ then $E_{s^\vee_n\mu} = t^{\frac12}_n \tau^\vee_n E_\mu$.
\end{enumerate}

\subsection{The creation formula for $E_\mu$}

The recursion of the previous subsection can be packaged nicely as a single formula for creating the Koornwinder polynomial
$E_\mu$.  This is the creation formula in~\eqref{Ecr}.

Let $\mu\in \ZZ^n$ and let $h_\mu\in W_X$ denote the corresponding translation.
Let $u_\mu\in W_X$ and $v_\mu\in W_{\mathrm{fin}}$ be as defined in~\eqref{umuvmudefn}, so that
$$u_\mu\in W_X\quad\hbox{and}\quad v_\mu\in W_{\mathrm{fin}}\quad\hbox{and}\quad
h_\mu = u_\mu v_\mu,\quad\hbox{with}\quad \ell(h_\mu) = \ell(u_\mu)+\ell(v_\mu).$$
Using the identification of $\KK[x_1^{\pm1}, \ldots, x_n^{\pm1}]$ with $\widetilde{H}_{\mathrm{int}}\mathbf{1}_Y$ as in~\eqref{indrep},
the creation formula for $E_\mu$ is
\begin{equation}
E_\mu 
= \frac{1}{\chi^+(T_{v_\mu^{-1}})} \tau_{u_\mu}^\vee \mathbf{1}_Y,
\qquad\hbox{where}\quad
\tau^\vee_{u_\mu} = \tau^\vee_{i_1}\cdots \tau^\vee_{i_\ell}
\label{Ecr}
%\tag{Ecr}
\end{equation}
if $u_\mu = s_{i_1}\cdots s_{i_\ell}$ is a reduced word for $u_\mu$.
\begin{proof}
Using
\begin{align*}
Y_i  \mathbf{1}_Y 
&= T_{\alpha_{i-1}}^{-1}\cdots T_{\alpha_1}^{-1}T_{\alpha_0}T_{\alpha_1}\cdots T_{\alpha_n}\cdots T_{\alpha_i} \mathbf{1}_Y
= t^{-\frac12(i-1)}t_0^{\frac12}t^{\frac12(n-1)}t_n^{\frac12}t^{\frac12(n-i)} \mathbf{1}_Y
=t^{-i} (t_0^{\frac12} t_n^{\frac12}t^n)  \mathbf{1}_Y,
\end{align*}
and $h_\mu^{-1}\varepsilon_j^\vee = h_{-\mu}\varepsilon_j^\vee = \varepsilon_j^\vee - \langle\mu, -\varepsilon_j^\vee\rangle K
=\varepsilon_j^\vee+\mu_j K$ 
gives
\begin{align*}
Y_j \tau^\vee_{u_\mu} \mathbf{1}_Y
&= Y^{\varepsilon_j^\vee} \tau^\vee_{u_\mu} \mathbf{1}_Y
= Y^{\varepsilon_j^\vee}  \tau^\vee_{u_\mu}\mathbf{1}_Y
= \tau^\vee_{u_\mu}Y^{u_\mu^{-1}\varepsilon_j^\vee}  \mathbf{1}_Y
= \tau^\vee_{u_\mu}Y^{v_\mu h_\mu^{-1}\varepsilon_j^\vee} \mathbf{1}_Y
\\
&= \tau^\vee_{u_\mu}Y^{v_\mu(\varepsilon_j^\vee + \mu_j  K)} \mathbf{1}_Y
= q^{-\mu_j} \tau^\vee_{u_\mu}Y^{\varepsilon_{v_\mu(j)}^\vee } \mathbf{1}_Y
\\
&= q^{-\mu_j} \tau^\vee_{u_\mu}Y_{v_\mu(j)} \mathbf{1}_Y
= q^{-\mu_j}t^{-v_\mu(j)} (t^{\frac12}_0t^{\frac12}_n t^n)^{\mathrm{sgn}(v_\mu(j))}
\tau^\vee_{u_\mu}\mathbf{1}_Y.
\end{align*}
Thus $\tau^\vee_{u_\mu}\mathbf{1}_Y$ is an eigenvector of $Y_j$
with eigenvalue $q^{-\mu_j}t^{-v_\mu(j)} (t^{\frac12}_0t^{\frac12}_n t^n)^{\mathrm{sgn}(v_\mu(j))}$.

Using the formulas~\eqref{tauvee} the product $\tau^\vee_{u_\mu}$ can be expanded in terms of the elements
$$\{ X^\gamma T_v\ |\ \gamma\in \ZZ^n, v\in W_{\mathrm{fin}}\}.$$
Since $u_\mu = t_\mu v_\mu^{-1}$ then the top term in this expansion is $X^\mu T_{v_\mu^{-1}}$ and
$$X^\mu T_{v_\mu^{-1}} \mathbf{1}_Y = X^\mu \chi^+(T_{v_\mu^{-1}}) \mathbf{1}_Y,
\qquad\hbox{where}\quad
\chi^+(T_{v_\mu^{-1}}) = t^{\frac12\ell_s(v_\mu^{-1})}t_n^{\frac12\ell_d(v_\mu^{-1})}.$$
Thus multiplying $\tau^\vee_{u_\mu}\mathbf{1}_Y$ by $\chi^+(T_{v_\mu^{-1}})^{-1}$ makes the
coefficient of $X^\mu$ equal to 1.
\end{proof}

\section{Bosonic, Fermionic and Mesonic}\label{sec:BFM}

The Weyl character formula is the formula that expresses the Schur function (a symmetric polynomial) as
a quotient of two determinants (antisymmetric polynomials).  There are Weyl character formulas in the
Koornwinder context as well (see Section~\ref{ssec:WCF}).  However, in the Koornwinder context, one finds that there are
\emph{four} Weyl character formulas, corresponding to the four symmetrizers $\varepsilon_+, \varepsilon_\pm, \varepsilon_\mp,
\varepsilon_-$ introduced in Section~\ref{ssec:symm}.

This section sets up the components for Weyl character formulas in the Koornwinder context.
There are four types of symmetrized Koornwinder polynomials:
the bosonic (symmetric) Koornwinder polynomials, the fermionic (antisymmetric) Koornwinder polynomials
and two types of mesonic (half symmetric-half antisymmetric) Koornwinder polynomials.  

The denominator in the classical Weyl character formula is the Vandermonde determinant, an antisymmetric polynomial
with a magical factorization.  In the Koornwinder case the Weyl denominators also have magical factorizations.
These Weyl denominators are presented in Section~\ref{ssec:Wdenoms}.

Every antisymmetric function can be obtained by multiplying a symmetric function by the Weyl denominators.
In~\cite[\S 4.3]{CR22} we viewed this correspondence between symmetric functions and antisymmetric functions as an
analog of the Boson-Fermion correspondence relating the symmetric algebra realization of Fock space
and exterior algebra realization of Fock space (a representation of a Heisenberg algebra, see~\cite[\S14.10]{Kac}).
In the Koornwinder context there are four spaces: the bosonic space (symmetric functions), fermionic space (antisymmetrized functions)
and two mesonic spaces (half symmetric-half antisymmetric functions).  These four spaces are all isomorphic as vector spaces, the isomorphisms
being given by multiplying by the different Weyl denominators.  This structure is explained Section~\ref{ssec:BFMspaces}.

In Sections~\ref{ssec:Ppoly},~\ref{ssec:Eexp} and~\ref{ssec:princspec},
we use the symmetrizers to give formulas for the Poincar\'e polynomial of $W_{\mathrm{fin}}$, for the $E$-expansions
of bosonic, fermionic and mesonic Koornwinder polynomials and formulas for the principal specializations.  These results
are Koornwinder analogues of the formulas in~\cite[Propositions 4,6 and 4.7 and Theorem 5.1]{CR22}.  
All of these formulas are given, in an even more general setting, in~\cite[(5,5,16), (5,7,8), (5.2.14),(5.3.9)]{Mac03}.

\subsection{Bosonic, Fermionic and Mesonic Macdonald-Koornwinder polynomials}

Let
$$(\ZZ_{\ge0}^n)^+ = \{
\lambda=(\lambda_1, \ldots, \lambda_n)\in \ZZ^n \ |\ \lambda_1\ge \cdots \ge \lambda_n\ge 0\}.
$$

For $\lambda\in (\ZZ^n_{\ge0})^+$, let
$$W_\lambda = \{ v\in W_{\mathrm{fin}}\ |\ v\lambda = \lambda\}
\qquad\hbox{and}\qquad
W_\lambda(t,t_n) = \sum_{v\in W_\lambda} t^{\ell_s(v)}t_n^{\ell_d(v)} = \sum_{v\in W_\lambda} \chi^+(T_v)^2.
$$
Let $\rho, \omega, \pi\in (\ZZ_{\ge 0}^n)^+$ be as defined in~\eqref{rhodef}.
Then, for $\lambda\in (\ZZ^n_{\ge0})^+$, define
the \emph{bosonic and fermionic Macdonald-Koornwinder polynomials} are
\begin{equation}
P_\lambda = \frac{\chi^+(T_{w_0}) }{ W_\lambda(t,t_n)} \varepsilon_+ E_\lambda
\qquad\hbox{and}\qquad
A_{\lambda+\rho} = \frac{\chi^+(T_{w_0}) }{ W_{\lambda+\rho}(t,t_n)}\varepsilon_- E_{\lambda+\rho},
\label{PAdef}
%\tag{PAdef}
\end{equation}
and the  \emph{mesonic Koornwinder polynomials} are
\begin{equation}
A^{\pm}_{\lambda+\omega} = \frac{\chi^+(T_{w_0}) }{ W_{\lambda+\omega}(t, t_n)}\varepsilon_{\pm} E_{\lambda+\omega}
\qquad\hbox{and}\qquad
A^{\mp}_{\lambda+\pi} = \frac{\chi^+(T_{w_0}) }{ W_{\lambda+\pi}(t,t_n)}\varepsilon_{\mp} E_{\lambda+\pi}.
\label{Mdef}
%\tag{Mdef}
\end{equation}

\begin{remark} \label{smstabs}
The notation in~\eqref{PAdef} and~\eqref{Mdef} displays the parallelism among the  expressions.  For computation it is useful to
note that the denominators can be given very explicitly:
\begin{align*}
W_{\lambda+\rho}(t,t_n) = 1, \qquad 
W_{\lambda+\omega}(t,t_n) = \prod_{i} \frac{(1-t^{m_i})}{1-t},
\qquad
W_{\lambda+\pi}(t,t_n) = \begin{cases}
t_n+1, &\hbox{if $\lambda_n=0$, } \\
1, &\hbox{if $\lambda_n\ne0$, }
\end{cases}
\end{align*}
where $m_i$ is the number of parts of size $i$ in $\lambda = (\lambda_1, \ldots, \lambda_n)$.
In particular, $W_{\lambda+\pi}(t,t_n)$ depends only on $t_n$ and $W_{\lambda+\omega}(t,t_n)$ depends only on $t$.
The factor $\chi^+(T_{w_0})$ guarantees that the coefficient of $X^{w_0\lambda} = X^{-\lambda}$ is equal to $1$ in $P_\lambda$.
Similarly, the coefficient of $X^{-(\lambda+\rho)}$ is equal to 1 in $A_{\lambda+\rho}$, 
the coefficient of $X^{-(\lambda+\omega)}$ is equal to 1 in $A^\pm_{\lambda+\omega}$ 
and the coefficient of $X^{-(\lambda+\pi)}$ is equal to 1 in $A^\mp_{\lambda+\pi}$.
\qed
\end{remark}

\subsection{Weyl denominators}\label{ssec:Wdenoms}

Define $a^\pm_\omega, a^\mp_\pi, a_\rho, A^\pm_\omega, A^\mp_\pi, A_\rho \in \KK[X]$ by
\begin{align}
a^{\pm}_\omega &=  x^{-\omega} \prod_{i=1}^n  (1-x_i^2), 
&A^\pm_\omega &= x^{-\omega} \left(\prod_{i=1}^n (1-t_n^{\frac12}u_n^{\frac12}x_i)(1+t_n^{\frac12}u_n^{-\frac12}x_i)\right),
\nonumber \\
a^{\mp}_\pi &= x^{-\pi} \prod_{1\le i<j\le n} (1-x_ix_j^{-1})(1-x_ix_j), 
&A^\mp_\pi &= x^{-\pi} \left(\prod_{1\le i<j\le n} (1-tx_ix_j)(1-tx_ix_j^{-1})\right),
\nonumber \\
a_\rho &=  a^{\mp}_\pi a^\pm_\omega,
&A_\rho &= A^\mp_\pi A^\pm_\omega.
\label{Wdenoms}
%\tag{Wdenoms}
\end{align}
Then
\begin{align*}
a^\pm_\omega &= \hbox{$\frac{1}{n!}$} e_\pm^X x^\omega, 
&a^\mp_\pi &= \hbox{$\frac{1}{2}$} e_\mp^X x^\pi,
&a_\rho &= e_-  x^\rho,
\\
A^\pm_\omega &= \frac{t^{\frac12 n(n-1)}t_n^{\frac12 n} }{ [n]! }\, \varepsilon_\pm^X x^\omega, 
&A^\mp_\pi &= \frac{ t^{\frac12 n(n-1)}t_n^{\frac12 n} }{(1+t_n)} \varepsilon_\mp^X x^\pi,
&A_\rho &= t^{\frac12 n(n-1)}t_n^{\frac12 n}\varepsilon_- x^\rho,
\end{align*}
where
$$t^{\frac12 n(n-1)}t_n^{\frac12 n} = \chi^+(T_{w_0}), \qquad
[n]! = \prod_{i=1}^n \frac{(1-t^i)}{1-t} = W_\omega(t,t_n), \qquad
1 + t_n = W_\pi(t,t_n).
$$
By~\eqref{cstcdr},
$$
\frac{A^{\pm}_\omega}{a^\pm_\omega} = \kappa_{\mathrm{dr}}^X,
\qquad
\frac{A^{\mp}_\pi}{a^\mp_\pi} = \kappa_{\mathrm{st}}^X, 
\qquad
\frac{A_\rho}{a_\rho} = \kappa_{w_0}^X.
$$
Since
$$A_\rho = \frac{\chi^+(T_{w_0})}{W_{\lambda+\rho}(t,t_n)} \varepsilon_- x^\rho 
= \frac{\chi^+(T_{w_0})}{W_{\lambda+\rho}(t,t_n)} \varepsilon_- E_\rho,$$
$$A^\pm_\omega = \frac{\chi^+(T_{w_0})}{W_\omega(t,t_n)} \varepsilon_\pm x^\omega 
= \frac{\chi^+(T_{w_0})}{W_\omega(t,t_n)} \varepsilon_\pm E_\omega, 
\qquad\hbox{and}\qquad
A^\mp_\pi = \frac{\chi^+(T_{w_0})}{W_\pi(t,t_n)} \varepsilon_\mp x^\pi
= \frac{\chi^+(T_{w_0})}{W_\pi(t,t_n)} \varepsilon_\mp E_\pi,
$$
there is no conflict of notation with the mesonic Macdonald polynomials introduced in~\eqref{Mdef}.

\subsection{Bosonic, fermionic and mesonic spaces}\label{ssec:BFMspaces}

The polynomial ring $\KK[X]$ is a module for the action of $\KK[X]^{W_{\mathrm{fin}}}$ and the structure of
$\KK[X]$ as a $\KK[X]^{W_{\mathrm{fin}}}$-module is of classical importance in the theory of reflection groups.
In fact, there are two commuting actions on $\KK[X]$, the action of $W_{\mathrm{fin}}$ and the action of
$\KK[X]^{W_{\mathrm{fin}}}$.  The part of this picture that is captured by the $X$-symmetrizers can be stated as follows.

Define
\begin{align*}
\KK[X]^{W_{\mathrm{fin}}} &= \{ f\in \KK[X]\ |\ \hbox{if $w\in W_{\mathrm{fin}}$ then $wf=f$}\}, \\
\KK[X]^{\pm} &= \{ f\in \KK[X]\ |\ \hbox{if $w\in W_{\mathrm{fin}}$ then $wf=(-1)^{\ell_d(w)} f$}\}, \\
\KK[X]^{\mp} &= \{ f\in \KK[X]\ |\ \hbox{if $w\in W_{\mathrm{fin}}$ then $wf=(-1)^{\ell_s(w)} f$}\}, \\
\KK[X]^{\det} &= \{ f\in \KK[X]\ |\ \hbox{if $w\in W_{\mathrm{fin}}$ then $wf=(-1)^{\ell_s(w)+\ell_d(w)} f$}\}.
\end{align*}
Then
\begin{align}
e_+^X\KK[X] &= \KK[X]^{W_{\mathrm{fin}}},
&e^X_\pm\KK[X] &= \KK[X]^{\pm} = a^{\pm}_{\omega} \KK[X]^{W_{\mathrm{fin}}},
\nonumber \\
e^X_-\KK[X] &= \KK[X]^{\mathrm{det}} = a_\rho \KK[X]^{W_{\mathrm{fin}}},
&e^X_\mp\KK[X] &= \KK[X]^{\mp} = a^{\mp}_\pi \KK[X]^{W_{\mathrm{fin}}}.
\label{Wharm}
%\tag{Wharm}
\end{align}

Now we proceed to a $t$-analogue of the equalities in~\eqref{Wharm}.
In this case Hecke algebra $H_{\mathrm{fin}}$ replaces the finite Weyl group, and the
actions of $H_{\mathrm{fin}}$ and $\KK[X]^{W_{\mathrm{fin}}}$ are commuting actions on 
$\KK[X]$.  The part of this picture captured by the Hecke symmetrizers is the following.

The bosonic, fermionic, and mesonic spaces are
\begin{align*}
\KK[X]^{\mathrm{Bos}} &= \{ f\in \KK[X]\ |\ \hbox{$T_nf = t_n^{\frac12}f$ and $T_if = t^{\frac12} f$ for $i\in \{1,\ldots, n\}$}\}, \\
\KK[X]^{\mathrm{Fer}} &= \{ f\in \KK[X]\ |\ \hbox{$T_nf = -t_n^{-\frac12}f$ and $T_if = -t^{-\frac12} f$ for $i\in \{1,\ldots, n\}$}\}, \\
\KK[X]^{\mathrm{Mes}\pm} &= \{ f\in \KK[X]\ |\ \hbox{$T_nf = -t_n^{-\frac12}f$ and $T_if = t^{\frac12} f$ for $i\in \{1,\ldots, n\}$}\}, \\
\KK[X]^{\mathrm{Mes}\mp} &= \{ f\in \KK[X]\ |\ \hbox{$T_nf = t_n^{\frac12}f$ and $T_if = -t^{-\frac12} f$ for $i\in \{1,\ldots, n\}$}\},
\end{align*}
With these definitions, the following proposition establishes $t$-analogues of the equalities in~\eqref{Wharm}.  The Weyl denominators
of Section~\ref{ssec:Wdenoms} are a key part of the structure.

\begin{prop} \label{BosFer}
Let $\varepsilon_+,  \varepsilon_\pm, \varepsilon_\mp, \varepsilon_-$ be the symmetrizers defined in~\eqref{symme-Hfin}.
\begin{align*}
\KK[X]^{\mathrm{Bos}} &= \varepsilon_+ \KK[X] =  \KK[X]^{W_{\mathrm{fin}}}, \qquad
&\KK[X]^{\mathrm{Fer}} &= \varepsilon_- \KK[X] = A_\rho \KK[X]^{W_{\mathrm{fin}}}, 
\\
\KK[X]^{\mathrm{Mes}\pm} &= \varepsilon_\pm \KK[X] = A^\pm_\omega\KK[X]^{W_{\mathrm{fin}}}, 
&\KK[X]^{\mathrm{Mes}\mp} &= \varepsilon_\mp \KK[X] = A^\mp_\pi \KK[X]^{W_{\mathrm{fin}}}, 
\end{align*}

Moreover, with $P_\lambda$, $A^\pm_{\lambda+\omega}$, $A^\mp_{\lambda+\pi}$ and $A_{\lambda+\rho}$ as in~\eqref{PAdef} and~\eqref{Mdef},
\begin{align*}
\{ P_\lambda\ |\ \lambda\in (\ZZ^n_{\ge0})^+\}
\quad\hbox{is a basis of $\varepsilon_+\KK[X]$,}
&&\{ A_{\lambda+\rho} \ |\ \lambda\in (\ZZ^n_{\ge0})^+\}
\quad\hbox{is a basis of $\varepsilon_-\KK[X]$,} \\
\{ A^\pm_{\lambda+\omega} \ |\ \lambda\in (\ZZ^n_{\ge0})^+\}
\quad\hbox{is a basis of $\varepsilon_{\pm}\KK[X]$,}
&&\{ A^\mp_{\lambda+\pi} \ |\ \lambda\in (\ZZ^n_{\ge0})^+\}
\quad\hbox{is a basis of $\varepsilon_{\mp}\KK[X]$,}
\end{align*}
\end{prop}
\begin{proof} 
We will give the proof for the $\pm$ case. The proofs for the other cases are similar.

Assume $f\in \varepsilon_\pm \KK[X]$.  Then there exists $g\in \KK[X]$ such that $f = \varepsilon_\pm g$ and
$$
T_n f = T_n \varepsilon_\pm g = -t_n^{-\frac12} \varepsilon_\pm g = -t_n^{-\frac12} f
\qquad\hbox{and}\qquad
T_i f = T_i\varepsilon_\pm g = t^{\frac12}\varepsilon_\pm g = t^{\frac12} f,
$$
for $i\in \{1, \ldots, n-1\}$.
So $f\in \KK[X]^{\mathrm{Mes}\pm}$ and $\varepsilon_\pm \KK[X] \subseteq \KK[X]^{\mathrm{Mes}\pm}$.

If $f\in \KK[X]^{\mathrm{Mes}\pm}$ then 
$$f = \frac{\chi^\pm(T_{w_0})}{W_0(t,t_n^{-1})}\varepsilon_\pm f 
= \frac{\chi^\pm(T_{w_0})}{W_0(t,t_n^{-1})}c_{\mathrm{dr}}^X e^X_\pm c_{\mathrm{st}}^{X^{-1}} f 
= \frac{\chi^\pm(T_{w_0})}{W_0(t,t_n^{-1})}\frac{A^\pm_\omega}{a^\pm_\omega}  e^X_\pm c_{\mathrm{st}}^{X^{-1}} f 
\in  
\frac{A^\pm_\omega}{a^\pm_\omega} \KK[X]^\pm.
$$
Since $\dfrac{A^\pm_\omega}{a^\pm_\omega} \KK[X]^\pm 
= \dfrac{A^\pm_\omega}{a^\pm_\omega} a^\pm_\omega \KK[X]^{W_{\mathrm{fin}}}
= A^\pm_\omega \KK[X]^{W_{\mathrm{fin}}}$ then
$\KK[X]^{\mathrm{Mes}\pm}\subseteq A^\pm_\omega \KK[X]^{W_{\mathrm{fin}}}$.

Assume $f\in A^\pm_\omega \KK[X]^{W_{\mathrm{fin}}}$.  Then there exists $g\in \KK[X]^{\mathrm{W_{\mathrm{fin}}}}$ such that
$f = A^\pm_\omega g$.  Then
$$f = A^\pm_\omega g = \left(\frac{\chi^+(T_{w_0})}{W_\omega(t,t_n)}\varepsilon_\pm x^\omega\right) g = \varepsilon_\pm \left(\frac{\chi^+(T_{w_0})}{W_\omega(t,t_n)}x^\omega g\right)\in \varepsilon_\pm \KK[X].$$
So $A^\pm_\omega \KK[X]^{W_{\mathrm{fin}}}\subseteq \varepsilon_\pm \KK[X]$.
This completes the proof that $\KK[X]^{\mathrm{Mes}\pm} = \varepsilon_\pm \KK[X] = A^\pm_\omega\KK[X]^{W_{\mathrm{fin}}}$.

Define
$$M^\pm_\mu = \varepsilon_\pm E_\mu\qquad\hbox{for $\mu\in \ZZ^n$.}$$
Let $i\in \{1, \ldots, n\}$.  
If $\mu_n=0$ then $s_n\mu = \mu$ and $T_nE_\mu = t_n^{\frac12}E_\mu$ so that
$$
M^\pm_\mu = \varepsilon_\pm T_{\alpha_n^\vee} E_\mu 
= (-t_n)^{\frac12} \varepsilon_\pm T_{\alpha_n^\vee} E_\mu = -M_\mu,
\quad\hbox{which forces\quad $M^\pm_\mu = 0$\quad when $\mu_n=0$.}
$$
If $i\in \{1, \ldots, n\}$ and $s_i\mu>\mu$ then
\begin{align*}
M^\pm_{s_i\mu} 
&= \varepsilon_\pm E_{s_i\mu} 
= \varepsilon_\pm t^{\frac12}\tau^\vee_i E_{\mu} 
= \varepsilon_\pm t^{\frac12}(T_{\alpha_i^\vee}+(t^{-\frac12}-c_{\alpha_i^\vee}^Y) E_{\mu} 
\\
&= \varepsilon_\pm t^{\frac12}(t^{\frac12}+t^{-\frac12}-c_{\alpha_i^\vee}^Y) E_{\mu} 
= \varepsilon_\pm t^{\frac12}c_{-\alpha_i^\vee}^Y E_{\mu} 
= \ev^t_\mu(\kappa_{-\alpha_i^\vee}^Y) M_{\mu},
\end{align*}
so that $M^\pm_\mu$ and $M^\pm_{s_i\mu}$ are linearly dependent.
It follows that
$$\varepsilon_\pm \KK[X] = \hbox{$\KK$-span}\{ M^\pm_\lambda\ |\ \lambda\in \ZZ^n\}
=\hbox{$\KK$-span}\{ M^\pm_{\lambda+\omega}\ |\ \lambda\in (\ZZ^n_{\ge 0})^+\}
=\hbox{$\KK$-span}\{ A^\pm_{\lambda+\omega}\ |\ \lambda\in (\ZZ^n_{\ge 0})^+\}
$$
Since $A^\pm_{\lambda+\omega}$ has top coefficient $x^{-(\lambda+\omega)}$ (in the DBlex order, see Section~\ref{sec:EPchar})
and the monomials are linear independent in $\KK[X]$ then the set
$\{ A^\pm_{\lambda+\omega}\ |\ \lambda\in (\ZZ^n_{\ge 0})^+\}$ is linearly independent.
\end{proof}

\subsection{Formulas for the Poincar\'e polynomial}\label{ssec:Ppoly}

Recall from~\eqref{Poindef} and~\eqref{evhom} that the Poincar\'e polynomial for $W_{\mathrm{fin}}$ is 
$$W_0(t,t_n) = \sum_{w\in W_{\mathrm{fin}}} \chi^+(T_w)^2
\qquad\hbox{and}\qquad
\ev^t_0(Y_i) = t^{n-i}t_0^{\frac12}t_n^{\frac12}$$
defines the evaluation homomorphism $\ev^t_0\colon \KK[Y]\to \KK$.
The following Proposition gives four ways of looking at the Poincar\'e polynomial: as a sum,
as a product, as an evaluation of $\kappa_{w_0}^Y$, and as a symmetrization of $\kappa_{w_0}^{X^{-1}}$.

\begin{prop} \label{Poinprop}
The group $W_{\mathrm{fin}}$ acts on 
$\KK[X]$ as in~\eqref{Wactions}.  Let $w_0$ denote the longest element of $W_{\mathrm{fin}}$
and let $\kappa^Y_{w_0}$ and $\kappa^{X^{-1}}_{w_0}$ be the noramlized $c$-functions given in 
\eqref{kw0Y} and~\eqref{kw0X}.
Then 
$$W_0(t,t_n)  = \prod_{i=1}^n \frac{(1-t^i)(1+t_nt^{i-1})}{(1-t)}
= \ev^t_0(\kappa_{w_0}^Y) = \sum_{w\in W_{\mathrm{fin}}} w(\kappa_{w_0}^{X^{-1}}).
$$
\end{prop}
\begin{proof}
Since $\chi^+(T_w)^2 = t^{\ell_s(w)}t_n^{\ell_d(w)}$ then
$$ 
\chi^+(T_{w_0})\varepsilon_+ \mathbf{1}_Y
= \sum_{w\in W_{\mathrm{fin}}} \chi^+(T_w) T_w \mathbf{1}_Y
= \sum_{w\in W_{\mathrm{fin}}} (\chi^+(T_w))^2 \mathbf{1}_Y
= W_0(t,t_n) \mathbf{1}_Y.
$$
Using the first formula in Proposition~\ref{cfcnsymm},
\begin{align*}
\chi^+(T_{w_0}) \varepsilon_+ \mathbf{1}_Y=  e_+^X \kappa_{w_0}^{X^{-1}} \mathbf{1}_Y
= \sum_{w\in W_{\mathrm{fin}}} \xi_w \kappa_{w_0}^{X^{-1}}\mathbf{1}_Y
= \sum_{w\in W_{\mathrm{fin}}} w(\kappa_{w_0}^{X^{-1}}).
\end{align*}
Using the second formula in Proposition~\ref{cfcnsymm},
$\varepsilon_+ \mathbf{1}_Y = e_+^Y \kappa^Y_{w_0}\mathbf{1}_Y$ and
\begin{align*}
\chi^+(T_{w_0}) \varepsilon_+ \mathbf{1}_Y
&= \left( \sum_{w\in W_{\mathrm{fin}}} \eta_w\right) \kappa_{w_0}^Y \mathbf{1}_Y
= \ev_0^t \left(\kappa_{w_0}^Y\right)  \left( 1+ \sum_{w\in W_{\mathrm{fin}}, w\ne 1} \eta_w \right) \mathbf{1}_Y 
\\
&= \ev_0^t\left(\kappa_{w_0}^Y\right) \left(1+0\right) \mathbf{1}_Y 
= \ev_0^t \left(\kappa_{w_0}^Y\right) \mathbf{1}_Y .
\end{align*}
Finally,
\begin{align*}
\ev_0^t \left(\kappa_{w_0}^Y\right)
&= \ev_0^t\left(\prod_{i<j} \frac{(1-t Y_iY_j^{-1})(1-t Y_iY_j) }{ (1-Y_iY_j^{-1})(1-Y_iY_j)}\right)
\ev_0^t\left( \prod_{i=1}^n \frac{(1-t_n^{\frac12}t_0^{\frac12}Y_i)(1+t_n^{\frac12}t_0^{-\frac12}Y_i)}{(1-Y_i^2)}
\right)
\\
&= 
\left( \prod_{i<j} \frac{(1-t t^{n-i}t_0^{\frac12}t_n^{\frac12} t^{-(n-j)} t_n^{-\frac12}t_0^{-\frac12})
(1-t t^{n-i}t_0^{\frac12}t_n^{\frac12}t^{n-j}t_0^{\frac12}t_n^{\frac12})}
{(1-t^{n-i}t_n^{\frac12}t_0^{\frac12}t^{-(n-j)}t_n^{-\frac12}t_0^{-\frac12})
(1-t^{n-i}t_n^{\frac12}t_0^{\frac12}t^{n-j}t_0^{\frac12}t_n^{\frac12})}\right)
\\
&\qquad
\left( \prod_{i=1}^n \frac{(1-t_n^{\frac12}t_0^{\frac12}t_0^{\frac12}t_n^{\frac12}t^{n-i})
(1+t_n^{\frac12}t_0^{-\frac12}t_0^{\frac12}t_n^{\frac12}t^{n-i})}{(1-t_0t_nt^{2(n-i)})}\right)
\\
&= 
\left( \prod_{i<j} \frac{(1-t^{j-i+1})(1-t^{2n-j-i+1}t_n t_0)}{(1-t^{j-i})(1-t^{2n-j-i}t_n t_0)}\right)
\left( \prod_{i=1}^n \frac{(1-t_n t_0 t^{n-i})
(1+t_n t^{n-i})}{(1-t_0t_nt^{2(n-i)})}\right)
\\
&=
\left( \prod_{j=2}^n \frac{(1-t^j)(1-t^{2n-j}t_nt_0)}{(1-t)(1-t^{2(n-j)+1}t_nt_0)}\right)
\left( \prod_{i=1}^n \frac{(1-t_nt_0 t^{n-i})(1+t_n t^{n-i}) }{ (1-t_0t_nt^{2(n-i)})}\right)
\\
&=
\left( \prod_{i=1}^n \frac{(1-t^i)(1+t_nt^i)}{(1-t)}\right)
\left( \prod_{j=2}^n \frac{(1-t^{2n-j}t_nt_0)}{(1-t^{2(n-j)+1}t_nt_0)}\right)
\left( \prod_{i=1}^n \frac{(1-t_nt_0 t^{n-i})}{(1-t_0t_nt^{2(n-i)})}\right)
\\
&=
\left( \prod_{i=1}^n \frac{(1-t^i)(1+t_nt^{i-1})}{(1-t)}\right)
\left( \prod_{i=1}^{2n} \frac{(1-t^{i-1} t_nt_0)}{(1-t^{i-1} t_nt_0)}\right)
=
 \prod_{i=1}^n \frac{(1-t^i)(1+t_nt^{i-1})}{(1-t)}.
\end{align*}
\end{proof}

\subsection{E-expansions}\label{ssec:Eexp}

The following Proposition uses the formulas for symmetrizers in terms of $c$-functions from Proposition~\ref{symmexp}
to give explicit expansion of the bosonic, fermionic and mesonic Koornwinder polynomials in terms of the $E_\mu$.
The coefficients in these expansions are evaluations of $c$-functions.  This is an example of how the $c$-functions
(which live in the field of fractions) appear in the structure even when doing expansions of polynomials that
live in $\KK[x_1^{\pm1}, \ldots, x_n^{\pm1}]$.

\begin{prop} \label{Eexps} (E-expansion formulas)
Let $\lambda\in (\ZZ_{\ge0}^n)^+$ and let $P_\lambda$, $A^\pm_{\lambda+\omega}$, $A^\mp_{\lambda+\pi}$ and
$A_{\lambda+\rho}$ be as defined in~\eqref{PAdef} and~\eqref{Mdef}.  Let
$\ev^t_\mu$ be the evaluation homomorphisms defined in~\eqref{evhom} and let 
$\kappa_{v_\mu}^\Xi$ be the normalized $c$-functions defined in~\eqref{kXi} and~\eqref{umuvmudefn}.
Then
\begin{align*}
P_\lambda &= \sum_{\mu\in W\lambda} \ev^t_\mu(\kappa^+_{v_\mu})E_\mu,
&A^\pm_{\lambda+\omega} &= \sum_{\mu\in W(\lambda+\omega)} (-1)^{\ell_d(v_\mu)}\ev^t_\mu(\kappa^\pm_{v_\mu})E_\mu,
\\
A_{\lambda+\rho} &= \sum_{\mu\in W(\lambda+\rho)} (-1)^{\ell(v_\mu)}\ev^t_\mu(\kappa^-_{v_\mu})E_\mu,
&A^\mp_{\lambda+\pi} &= \sum_{\mu\in W(\lambda+\pi)} (-1)^{\ell_s(v_\mu)} \ev^t_\mu(\kappa^\mp_{v_\mu})E_\mu.
\end{align*}
\end{prop}
\begin{proof}
Let us do the case $A^{\pm}_{\lambda+\omega}$. The other cases are similar. 

Since $\tau_i^\vee = \eta_{s_i}^Y c_{\alpha^\vee_i}^Y$ and $t^{\frac12}\tau^\vee_i E_\mu = E_{s_i^\vee\mu}$ then
$$E_{s_i^\vee\mu} = t^{\frac12}\tau_i^\vee E_\mu 
= t^{\frac12} \eta_{s_i}^Y c_{\alpha^\vee_i}^Y E_\mu = \eta_{s_i}^Y \kappa_{\alpha_i}^Y E_\mu 
= \kappa_{\alpha_i}^{Y^{-1}} \eta_{s_i}^Y  E_\mu.
$$
If $z\in W^{\lambda+\omega}$ then
$$E_{z(\lambda+\omega)} = \eta_{s_{i_1}}\kappa_{\alpha_{i_1}}^Y \cdots \eta_{s_{i_\ell}}\kappa_{\alpha_{i_\ell}}^Y E_{\lambda+\omega}
= \eta_z \kappa^Y_z E_{\lambda+\omega} = \kappa_{z^{-1}}^{Y^{-1}}\eta_z E_{\lambda+\omega}.
$$
If $w\in W_{\lambda+\omega}$ then
$$T_w E_{\lambda+\omega} = \chi^+(T_w) E_{\lambda+\omega}
\qquad\hbox{(in the same way that $T_w\mathbf{1}_Y = \chi^+(T_w)\mathbf{1}_Y$ for $w\in W_{\mathrm{fin}}$).}
$$
This gives that
$$\varepsilon^+_{\lambda+\omega}E_{\lambda+\omega} 
= \frac{1}{\chi^+(T_{w_{\lambda+\omega}})} \sum_{w\in W_{\lambda+\omega}} \chi^+(T_w)^2 E_{\lambda+\omega}
= \frac{1}{\chi^+(T_{w_{\lambda+\omega}})} W_{\lambda+\omega}(t,t_n) E_{\lambda+\omega}.
$$
Using Proposition~\eqref{symmexp} gives 
\begin{align*}
A^\pm_{\lambda+\omega} 
&= \frac{\chi^+(T_{w_0}) }{W_{\lambda+\omega}(t,t_n)}
\varepsilon_\pm E_{\lambda+\omega}
= \frac{1}{W_{\lambda+\omega}(t,t_n)}
 \left(\sum_{v\in W^{\lambda+\omega}} (-1)^{\ell_d(z)} \kappa^\pm_{v_{\lambda+\omega}z} \eta^Y_z \kappa^Y_z \right)
\chi^+(T_{w_{\lambda+\omega}}) \varepsilon^+_{\lambda+\omega} E_{\lambda+\omega}
\\
&= 
 \left(\sum_{v\in W^{\lambda+\omega}} (-1)^{\ell_d(z)} \kappa^\pm_{v_{\lambda+\omega}z} \eta^Y_z \kappa^Y_z\right)
E_{\lambda+\omega}
= 
\sum_{v\in W^{\lambda+\omega}} (-1)^{\ell_d(z)} \kappa^\pm_{v_{\lambda+\omega}z}
E_{z(\lambda+\omega)}
\\
&= 
\sum_{z\in W_{\mathrm{fin}}} (-1)^{\ell_d(z)} \ev^t_{z(\lambda+\omega)}(\kappa_{v_{\lambda+\omega}z}^\pm) E_{z(\lambda+\omega)}.
\end{align*}
\end{proof}

\subsection{Principal specializations}\label{ssec:princspec}

One of the most pleasing combinatorial miracles in Lie theory is that principal specializations of Schur functions and
Weyl characters factor as products (see~\cite[\S10.9]{Kac} and~\cite[Ch.\ I \S3 Ex.\ 1]{Mac}).
This feature extends to Macdonald-Koornwinder  polynomials, and the result in this subsection shows that the
principal specializations of Macdonald-Koornwinder polynomials
are evaluations of $c$-functions which come naturally out of the recursive construction of the 
electronic Macdonald polynomial $E_\mu$.

Define ring homomorphisms $\ev^t_0\colon \KK[Y]\to \KK$ and $\ev^{t^{-1}}_0\colon \KK[Y]\to \KK$ by
$$\ev^t_0(Y_i) = t^{n-i}t_0^{\frac12}t_n^{\frac12}
\qquad\hbox{and}\qquad
\ev^t_0(Y^{-1}_i) = t^{-(n-i)}t_0^{-\frac12}t_n^{-\frac12},
\qquad\hbox{for $i\in \{1, \ldots, n\}$.}
$$

\begin{thm} Let $\mu, \lambda\in \ZZ^n$ with $\lambda_1\ge \cdots \ge \lambda_n \ge 0$.  Let $u_\mu$ and $h_\lambda$
be as defined in~\eqref{umuvmudefn} and the normalized $c$-functions $\kappa_{u_\mu}^{Y^{-1}}$ and $\kappa_{h_\lambda}^{Y^{-1}}$
as in~\eqref{kappawdef}.
Let
$$
a_1 = 
t^{-(n-1)}(t_0t_n)^{-\frac12},\quad
a_2 = 
t^{-(n-2)}(t_0t_n)^{-\frac12},\quad
\ldots, \quad
a_n = 
(t_0t_n)^{-\frac12}.
$$
Then 
\begin{align*}
E_\mu( a_1, \ldots, a_n;q,t,t_0^{\frac12},u_0^{\frac12}, t_n^{\frac12}, u_n^{\frac12})
&= \frac{1}{\chi^+(T_{v_\mu^{-1}})} \ev^t_0(c_{u_\mu}^{Y^{-1}})
\qquad\hbox{and}\qquad
\\
P_\lambda(a_1, \ldots, a_n;q,t,t_0^{\frac12},u_0^{\frac12}, t_n^{\frac12}, u_n^{\frac12})
&= \ev_0^{t^{-1}}(c_{h_\lambda}^{Y^{-1}}).
\end{align*}
\end{thm}
\begin{proof}
For this proof use the realization of the polynomial representation $\KK[X]$ as an induced module
$\widetilde H\mathbf{1}_Y$ via the $\widetilde H$-module isomorphism of~\eqref{indrep}.  
Let $\mathbf{1}_X$ be a formal symbol which satisfies $\mathbf{1}_XT_j = t^{\frac12}\mathbf{1}_X$ and
$\mathbf{1}_X T_0^\vee = t_0^{\frac12}\mathbf{1}_X$.
Using
$X_1 = (T_0^\vee)^{-1}T_1^{-1}\cdots T_n^{-1}\cdots T_1^{-1}$ and
$X_{i+1} = T_iX_iT_i$
gives
$$\mathbf{1}_X X_i = t^{-(n-1)}t_n^{-\frac12}t_0^{-\frac12} t^{i-1},
\qquad\hbox{for $i\in \{1, \ldots, n\}$.}
$$
Thus, if $\mu\in \ZZ^n$ then
$$\mathbf{1}_X E_\mu(x_1, \ldots, x_n; q,t, t_0^{\frac12}, u_0^{\frac12}, t_n^{\frac12}, u_n^{\frac12}) 
= \mathbf{1}_X  E_\mu(a_1, \ldots, a_n;q,t,t_0^{\frac12}, u_0^{\frac12}, t_n^{\frac12}, u_n^{\frac12}).
$$
For $i\in \{0,1,\ldots, n\}$
$$\mathbf{1}_X\tau^\vee_i
= \mathbf{1}_X\left( T_{\alpha^\vee_i} + (c_{\alpha^\vee_i}^{Y^{-1}}-t^{\frac12})\right)
= \mathbf{1}_X\left(t^{\frac12}+ (c_{\alpha^\vee_i}^{Y^{-1}}-t^{\frac12})\right)
=\mathbf{1}_X c_{\alpha_i^\vee}^{Y^{-1}}.
$$
By~\eqref{eveg},
$$c_{\alpha^\vee_i}^{Y^{-1}}\mathbf{1}_Y = \ev_0^t(c_{\alpha^\vee_i}^{Y^{-1}})\mathbf{1}_Y.$$
If $w\in W$ and $\ell(s_iw)> \ell(w)$ then
$$
\mathbf{1}_X \tau^\vee_i \tau^\vee_w \mathbf{1}_Y
=\mathbf{1}_X c_{\alpha^\vee_i}^{Y^{-1}} \tau^\vee_w \mathbf{1}_Y
=\mathbf{1}_X  \tau^\vee_w c_{w^{-1}\alpha^\vee_i}^{Y^{-1}}\mathbf{1}_Y
=\ev_0^t(c_{w^{-1}\alpha^\vee_i}^{Y^{-1}}) \mathbf{1}_X  \tau^\vee_w \mathbf{1}_Y.
$$
This is the induction step giving that if $w\in W$ and $w=s_{i_1}\cdots s_{i_\ell}$ is a reduced word for $w$ then
$$\mathbf{1}_X \tau^\vee_w \mathbf{1}_Y = \mathbf{1}_X \tau^\vee_{i_1}\cdots \tau^\vee_{i_\ell} \mathbf{1}_Y
= \mathbf{1}_X \ev_0^t(c_w^{Y^{-1}})\mathbf{1}_Y = \ev_0^t(c_w^{Y^{-1}}) \mathbf{1}_X\mathbf{1}_Y.
$$
Thus
$$\mathbf{1}_X E_\mu \mathbf{1}_Y 
= \frac{1}{\chi^+(T_{v_\mu^{-1}})} \mathbf{1}_X \tau_{u_\mu}^\vee \mathbf{1}_Y
= \frac{1}{\chi^+(T_{v_\mu^{-1}})} \ev^t_0(c_{u_\mu}^{Y^{-1}}) \mathbf{1}_X \mathbf{1}_Y.
$$

Using $\mathbf{1}_X \varepsilon_+ = \frac{1}{\chi^+(T_{w_0})}W_0(t,t_n) \mathbf{1}_X$ from~\eqref{symmsq} gives
\begin{align*}
P_\lambda&(a_1,\ldots,a_n;q,t,t_0^{\frac12}, u_0^{\frac12}, t_n^{\frac12}, u_n^{\frac12})\mathbf{1}_X\mathbf{1}_Y
= \mathbf{1}_X P_\lambda \mathbf{1}_Y 
= \mathbf{1}_X \frac{\chi^+(T_{w_0})}{W_\lambda(t,t_n)} \varepsilon_+ E_\lambda \mathbf{1}_Y
\\
&=  \frac{ W_0(t,t_n)}{W_\lambda(t,t_n)} \mathbf{1}_X E_\lambda \mathbf{1}_Y
= \frac{1}{\chi^+(T_{v_\lambda^{-1}})} \frac{W_0(t,t_n)}{W_\lambda(t,t_n)} \ev^t_0(c_{u_\lambda}^{Y^{-1}})
\mathbf{1}_X \mathbf{1}_Y.
\end{align*}
Let
$$w_\lambda\in W_{\mathrm{fin}}\qquad\hbox{be the longest element of}\quad
W_\lambda = \{ v\in W_{\mathrm{fin}}\ |\ v\lambda = \lambda\}.$$
Then by Proposition~\ref{Poinprop}
$$v_\lambda^{-1} = (w_0w_\lambda)^{-1} = w_\lambda w_0
\qquad\hbox{and}\qquad
W_\lambda(t,t_n) = \ev^t_0(\kappa_{w_\lambda}^Y) = \ev^{t^{-1}}_0(\kappa_{w_\lambda}^{Y^{-1}}).$$

Since $h_\lambda = u_\lambda v_\lambda$ then $u_\lambda = h_\lambda v_\lambda^{-1} = v_\lambda^{-1} h_{v_\lambda \lambda}
=v_\lambda^{-1} h_{v_\lambda \lambda} = v_\lambda^{-1} h_{w_0\lambda}$. 
Using this and $\ev^t_0(Y_i)=\ev^{t^{-1}}_0(Y_i^{-1}) = \ev^{t^{-1}}_0(Y_{-i}) = \ev^{t^{-1}}_0(Y_{w_0(i)})$ gives
$$\ev^t_0(c_{u_\lambda}^{Y^{-1}}) = \ev^{t^{-1}}_0(v_\lambda^{-1}  c_{u_\lambda}^{Y^{-1}}).$$

Therefore
\begin{align*}
\frac{1}{\chi^+(T_{v_\lambda^{-1}})} \frac{W_0(t,t_n)}{W_\lambda(t,t_n)} \ev^t_0(c_{u_\lambda}^{Y^{-1}})
&= \frac{\chi^+(T_{w_\lambda}) }{\chi^+(T_{w_0} )}  \ev^{t^{-1}}_0\left( \frac{\kappa_{w_0}^{Y^{-1}}}{\kappa_{w_\lambda}^{Y^{-1}}}\right) 
\ev^t_0(c_{u_\lambda}^{Y^{-1}})
= \ev^{t^{-1}}_0\left( \frac{ c_{w_0}^{Y^{-1}} }{ c_{w_\lambda}^{Y^{-1}}}\right) \ev^{t^{-1}}_0(v_\lambda^{-1}  c_{u_\lambda}^{Y^{-1}})
\\
&=\ev^{t^{-1}}_0(c_{v_\lambda}^{Y^{-1}}) \ev^{t^{-1}}_0(v_\lambda^{-1}  c_{u_\lambda}^{Y^{-1}})
=\ev^{t^{-1}}_0( c_{u_\lambda v_\lambda}^{Y^{-1}})
=\ev^{t^{-1}}_0( c_{h_\lambda}^{Y^{-1}}),
\end{align*}
which completes the proof of the second statement.
\end{proof}

\begin{remark}
The principal specializations of the fermionic and mesonic Macdonald polynomials are
$$
A_{\lambda+\rho}(a_1, \ldots, a_n;q,t,t_0^{\frac12},u_0^{\frac12}, t_n^{\frac12}, u_n^{\frac12}) =0, 
$$
$$A^{\pm}_{\lambda+\omega}(a_1, \ldots, a_n; q, t, t_0^{\frac12},u_0^{\frac12}, t_n^{\frac12}, u_n^{\frac12}) =0,
\qquad\hbox{and}\qquad
A^{\mp}_{\lambda+\pi}(a_1, \ldots, a_n; q, t, t_0^{\frac12},u_0^{\frac12}, t_n^{\frac12}, u_n^{\frac12}) =0.
$$
To establish this for $A^\pm_{\lambda+\omega}$,  
use $\mathbf{1}_X \varepsilon_\pm = 0$ to get
\begin{align*}
A^\pm_{\lambda+\omega}&(a_1, \ldots, a_n; q, t, t_0^{\frac12}, u_0^{\frac12}, t_n^{\frac12}, u_n^{\frac12})\mathbf{1}_X\mathbf{1}_Y
=\mathbf{1}_X A^\pm_{\lambda+\omega} \mathbf{1}_Y
=\frac{\chi^+(T_{w_0})}{W_{\lambda+\omega}(t,t_n)} \mathbf{1}_X \varepsilon_\pm E_{\lambda+\omega}  \mathbf{1}_Y=0.
\end{align*}
The proof for the other cases is similar.
\qed
\end{remark}

\section{Orthogonality}\label{sec:orth}

In this section we study the Koornwinder polynomials as a family of orthogonal polynomials
for a specific inner product.  The inner product $(\ ,\ )_+$ is defined via multiplication by a kernel and
taking the constant term (for those with an analytic bent, taking the constant term is
an integral and the kernel is what defines the measure for the integral).  
The Macdonald-Koornwinder inner product is defined in Sections~\ref{ssec:kernel} and~\ref{ssec:innprod}.
The kernel is a huge product of $c$-functions, one for each positive root in the affine root system of type $CC_n$.

The Koornwinder polynomials are characterized by orthogonality with respect to this inner product
and a triangular expansion in terms of monomials.
In order to use the Hecke algebra as a tool in the inner product setting it is crucial to establish
that the adjoints of operators that come from the Hecke algebra are tractable.  This is done in Section~\ref{ssec:adjoints}.
In particular, we find that the symmetrizers are self adjoint operators.

The proof of the Weyl character formulas and the norm formulas for Koornwinder polynomials
rely on a shift of parameters coming from multiplying by the Weyl denominators.
These going up a level formulas, derived in Section~\ref{ssec:levelup}, 
are the key to establishing recursive relations for computing norms.
The recursive relations are derived in Section~\ref{ssec:normred} and the norm formula for $(P_\lambda, P_\lambda)_+$
is established in Section~\ref{ssec:normformula}.

In the same way that there are four symmetrizers, there are four going up a level formulas, four Weyl character formulas,
four types of recursion relations.  In each case, one of the four formulas is usually a triviality, but we have included these
trivial formulas in our exposition each time in order to highlight the underlying symmetry of the structures.  In the end,
the various formulas combine and complement each other to provide the inductive structure for computing norms in terms of 
$c$-functions.

\subsection{The kernel $\Delta_{CC}^+$}\label{ssec:kernel}

For an affine root $\beta$, define
\begin{equation}
\kappa_\beta^X = t_\beta^\frac12 c_\beta^X  =
\frac{(1-t_{\beta}^{\frac12}u_{\beta}^{\frac12}X^{\beta})
(1+t_{\beta}^{\frac12}u_{\beta}^{-\frac12}X^{\beta})}
{(1-X^{2\beta})}.
\label{kpdefn}
%\tag{kpdefn}
\end{equation}
More specifically, if $i,j\in \{1, \ldots, n\}$ with $i<j$ and $r\in \ZZ_{\ge 0}$ then 
\begin{align*}
\kappa_{\varepsilon_i-\varepsilon_j+(r+1)\delta}^X &= \frac{1-tq^{r+1}X_iX_j^{-1}}{1-q^{r+1}X_iX_j^{-1} },
\quad
&\kappa_{\varepsilon_i+(r+1)\delta}^X 
&= \frac{(1-q^{r+1}t_n^{\frac12}u_n^{\frac12}X_i)(1+q^{r+1}t_n^{\frac12}u_n^{-\frac12}X_i) }{ 1-q^{2r+2}X_i^2},
\\
\kappa_{\varepsilon_i+\varepsilon_j+(r+1)\delta}^X &= \frac{1-tq^{r+1}X_iX_j }{1-q^{r+1}X_iX_j },
&\kappa_{\varepsilon_i+(r+\frac12)\delta}^X 
&= \frac{(1-q^{r+\frac12}t_0^{\frac12}u_0^{\frac12}X_i)(1+q^{r+\frac12}t_0^{\frac12}u_0^{-\frac12}X_i)}
{1-q^{2r+1}X_i^2}.
\end{align*}
Let $S^+$ be the set of positive roots for the affine root system so that
\begin{align*}
S^+ &= S_{g,+}^+\cup S_{g,-}^+ \cup S_{s,+}^+\cup S_{s,-}^+ \cup S_{d,+}^+ \cup S_{d,-}^+ \cup S_{0,s}^+ \cup S_{0,d}^+,
\end{align*}
where
\begin{align*}
S^+_{s,+} &= 
\left\{ 
\begin{array}{l}
(\varepsilon_i-\varepsilon_j) + (r+1)\delta, \\
(\varepsilon_i+\varepsilon_j) + (r+1)\delta
\end{array}
\ \Big\vert \ \begin{array}{l} \hbox{$i,j\in \{1, \ldots, n\}$,} \\ \hbox{$i<j$ and $r\in \ZZ_{\ge 0}$}\end{array} \right\},
\\
S^+_{s,-} &= 
\left\{ 
\begin{array}{l}
-(\varepsilon_i-\varepsilon_j) + (r+1)\delta, \\
-(\varepsilon_i+\varepsilon_j) + (r+1)\delta
\end{array}
\ \Big\vert \ \begin{array}{l} \hbox{$i,j\in \{1, \ldots, n\}$,} \\ \hbox{$i<j$ and $r\in \ZZ_{\ge 0}$}\end{array} \right\},
\end{align*}
\begin{align*}
S^+_{g,+} &= \left\{ \varepsilon_i + (r+\hbox{$\frac12$})\delta\ \Big\vert\  \begin{array}{c} i\in \{1, \ldots, n\}\\ r\in \ZZ_{\ge 0}\end{array} \right\},
&S^+_{d,+} &= \left\{ \varepsilon_i + (r+1)\delta\ \Big\vert\  \begin{array}{c} i\in \{1, \ldots, n\}\\ r\in \ZZ_{\ge 0}\end{array} \right\},
\\
S^+_{g,-} &= \left\{ -\varepsilon_i + (r+\hbox{$\frac12$})\delta\ \Big\vert\  \begin{array}{c} i\in \{1, \ldots, n\}\\ r\in \ZZ_{\ge 0}\end{array} \right\}, 
&S^+_{d,-} &= \left\{ - \varepsilon_i + (r+1)\delta\ \Big\vert\  \begin{array}{c} i\in \{1, \ldots, n\}\\ r\in \ZZ_{\ge 0}\end{array} \right\}, \\
S^+_{0,s} &= \{ \varepsilon_i \pm \varepsilon_j \ |\ \hbox{$i,j\in \{1, \ldots, n\}$ with $i<j$} \}, 
&S^+_{0,d} &= \{ \varepsilon_i \ |\ i\in \{1, \ldots, n\} \}.
\end{align*}
Then define
\begin{equation}
\Delta^+_{CC} = \prod_{\beta\in S^+} \kappa_\beta^X.
\label{rootformker}
%\tag{rootformker}
\end{equation}
and
\begin{align*}
\Delta^X_g &= \prod_{\beta\in S_{g,+}^+} \kappa_\beta^X, 
&\Delta^X_s &= \prod_{\beta\in S_{s,+}^+} \kappa_\beta^X, 
&\Delta^X_d &= \prod_{\beta\in S_{d,+}^+} \kappa_\beta^X, \qquad\quad
&\Delta^X_{0,s} &= \prod_{\beta\in S_{0,s}^+} \kappa_\beta^X,
\\
\Delta^{X^{-1}}_g &= \prod_{\beta\in S_{g,-}^+} \kappa_\beta^X, 
&\Delta^{X^{-1}}_s &= \prod_{\beta\in S_{s,-}^+} \kappa_\beta^X, 
&\Delta^{X^{-1}}_d &= \prod_{\beta\in S_{d,-}^+} \kappa_\beta^X, 
&\Delta^X_{0,d} &= \prod_{\beta\in S_{0,d}^+} \kappa_\beta^X,
\end{align*}
so that
$$\Delta^+_{CC} = \Delta_g^X \Delta_g^{X^{-1}} \Delta_s^X\Delta_s^{X^{-1}} 
\Delta_d^X\Delta_d^{X^{-1}} \Delta_{0,s}^X \Delta_{0,d}^X.$$

\begin{remark}\label{Deltaexp}
In terms of the Askey-Wilson parameters $a,b,c,d$ (see~\eqref{AWparams}), the expression $\Delta^+_{CC}$ 
used to define the inner product is
$$\Delta^+_{CC} = \Delta^{(1)}\Delta^{(2)},$$
where
$$\Delta^{(2)} = \prod_{1\le i< j\le n} 
\frac{(tx_ix_j^{-1};q)_\infty (tx_ix_j;q)_\infty (tqx_i^{-1}x_j;q)_\infty (tqx_i^{-1}x_j^{-1};q)_\infty}
{ (x_ix_j^{-1};q)_\infty (x_ix_j;q)_\infty (qx_i^{-1}x_j;q)_\infty (qx_i^{-1}x_j^{-1};q)_\infty}
$$
and
\begin{align*}
\Delta^{(1)}
&= \prod_{i=1}^n 
\frac {
(cx_i;q)_\infty (qc x_i^{-1};q)_\infty
(dx_i;q)_\infty (qd x_i^{-1};q)_\infty
(ax_i;q)_\infty (a x_i^{-1};q)_\infty
(bx_i;q)_\infty (b x_i^{-1};q)_\infty
}{ (x_i^2;q)_\infty (qx_i^{-2};q)_\infty}.
\end{align*}
This is verified by noting that $\Delta^{(2)} = \Delta_{s}^X \Delta_{0,s}^X\Delta_s^{X^{-1}}$, and that $\Delta^{(1)} = \Delta_{g}^X \Delta_{g}^{X^{-1}}\Delta_{d}^X \Delta_{d}^{X^{-1}}\Delta_{0,d}^X$ since
\begin{align*}
\prod_{i=1}^n \prod_{r\in \ZZ_{\ge 0}} \frac{(cx_i;q)_\infty (dx_i; q)_\infty }{(x_i^2;q^2)_\infty}
&= \prod_{i=1}^n \prod_{r\in \ZZ_{\ge 0}} \frac{(1-t_n^{\frac12}u_n^{\frac12} q^r x_i)(1+t_n^{\frac12}u_n^{-\frac12} q^r x_i)}{(1-q^{2r}x_i^2)}
=\Delta_{d}^X\Delta_{0,d}^X,
\\
\prod_{i=1}^n \prod_{r\in \ZZ_{\ge 0}} \frac{(qc x_i^{-1};q)_\infty (qd x_i^{-1}; q)_\infty }{(q^2x_i^{-2};q^2)_\infty}
&= \prod_{i=1}^n \prod_{r\in \ZZ_{\ge 0}} 
\frac{(1-t_n^{\frac12}u_n^{\frac12} q^{r+1} x_i^{-1})(1+t_n^{\frac12}u_n^{-\frac12} q^{r+1} x_i^{-1})}
{(1-q^{2r+2} x_i^{-2})}
=\Delta_d^{X^{-1}},
\\
\prod_{i=1}^n \prod_{r\in \ZZ_{\ge 0}} \frac{(ax_i;q)_\infty (bx_i; q)_\infty }{(qx_i;q^2)_\infty}
&= \prod_{i=1}^n \prod_{r\in \ZZ_{\ge 0}} 
\frac{(1-t_0^{\frac12}u_0^{\frac12} q^{r+\frac12} x_i)(1+t_0^{\frac12}u_0^{-\frac12} q^{r+\frac12} x_i)}
{(1-q^{2r+1}x_i^2)}
=\Delta_g^{X},
\\
\prod_{i=1}^n \prod_{r\in \ZZ_{\ge 0}} \frac{(a x_i^{-1};q)_\infty (b x_i^{-1}; q)_\infty }{(qx_i^{-2};q^2)_\infty}
&= \prod_{i=1}^n \prod_{r\in \ZZ_{\ge 0}} 
\frac{(1-t_0^{\frac12}u_0^{\frac12} q^{r+\frac12} x_i^{-1})(1+t_0^{\frac12}u_0^{-\frac12} q^{r+\frac12} x_i^{-1})}
{(1-q^{2r+1} x_i^{-2})}
=\Delta_g^{X^{-1}}.
\end{align*}
\qed\end{remark}

\subsection{Definition of the inner product}\label{ssec:innprod}

Let $\KK[x] = \KK[x_1^{\pm 1}, \ldots, x_n^{\pm1}]$.  Define an involution $\overline{\phantom{T}}\colon \KK[x] \to \KK[x]$ by
\begin{equation}
\overline{f}(x_1, \ldots, x_n;q,t,t_0^{\frac12},u_0^{\frac12}, t_n^{\frac12}, u_n^{\frac12}) 
= f(x_1^{-1}, \ldots, x_n^{-1}; q^{-1}, t^{-1}, t_0^{-\frac12}, u_0^{-\frac12}, t_n^{-\frac12}, u_n^{-\frac12}).
\label{invdefn}
%\tag{invdefn}
\end{equation}

Let $\Delta^+_{CC}$ be as defined in~\eqref{rootformker}.
Define a scalar product $(\ ,\ )_+\colon \KK[x]\times \KK[x] \to \KK$ by
\begin{equation}
(f_1, f_2)_+
= \mathrm{ct}\left( \frac{f_1 \overline{f_2}}
{\Delta^+_{CC}}\right),
\qquad\hbox{where}\quad
\mathrm{ct}(f) = (\hbox{constant term in $f$}), \quad \hbox{for $f\in \KK[x]$.}
\label{ipdef}
%\tag{ipdef}
\end{equation}

\subsection{The inner product characterization of $E_\mu$ and $P_\lambda$} \label{sec:EPchar}

Define
$$(\ZZ^n)^+ = \{ (\gamma_1, \ldots, \gamma_n)\in \ZZ^n
\ |\ \gamma_1\ge \cdots \ge \gamma_n\ge 0\}.$$
The elements of $(\ZZ^n)^+$ are partially ordered by the dominance order: For $\lambda, \mu\in (\ZZ_{\ge 0}^n)^+$,
$$\lambda < \mu \qquad\hbox{if}\qquad
\lambda_1+\cdots+\lambda_i \le \mu_1+\cdots +\mu_i,\quad\hbox{for $i\in \{1, \ldots, n\}$.}
$$
The elements of $\ZZ^n$ are partially ordered by the DBlex order: For $\lambda, \mu\in \ZZ^n$,
$$\lambda\le \mu 
\quad\hbox{if }\qquad
\begin{array}{c}
\hbox{$\lambda^+< \mu^+$ in dominance order} \\
\hbox{or} \\
\hbox{$\lambda^+=\mu^+$ and $z_\lambda<z_\mu$ in Bruhat order on $W_{\mathrm{fin}}$,}
\end{array}
$$
where 
$\lambda^+ \in W_{\mathrm{fin}}\lambda\cap  (\ZZ^n_{\ge 0})^+$ and
$z_\lambda$ is the minimal length element of $W_{\mathrm{fin}}$ so that $\lambda =z_\lambda \lambda^+$.

For $\mu = (\mu_1, \ldots, \mu_n) \in \ZZ^n$ write
$x^\mu = x_1^{\mu_1}\cdots x_n^{\mu_n}$
and for $\gamma\in (\ZZ^n)^+$, define the \emph{monomial symmetric function} $m_\gamma$ by
$$m_\gamma = \sum_{\mu\in W_{\mathrm{fin}} \gamma} x^\mu,
\qquad
\hbox{where the sum is over all elements of the orbit $W_{\mathrm{fin}}$-orbit of $\gamma$.}
$$
With these definitions we have the following characterizations of the $E_\mu$ and the $P_\lambda$.  The proofs of
Propositions~\ref{Ebyinnerprod} and~\ref{Pbyinnerprod} are exactly as in~\cite[(5.2.1) and (5.3.1)]{Mac03} and
\cite[Prop.\ 6.2 and 6.3]{CR22}.

\begin{prop}  \label{Ebyinnerprod} Let $\mu \in \ZZ^n$.
The electronic Macdonald polynomial $E_\mu$ is the unique element of
$\KK[x_1^{\pm1}, \ldots, x_n^{\pm1}]$ such that
\begin{enumerate}
\item[(a)] $E_\mu = x^\mu + (\hbox{lower terms})$;
\item[(b)] If $\nu\in \ZZ^n$ and  $\nu<\mu$ then $(E_\mu, x^\nu)_+ = 0$.
\end{enumerate}
\end{prop}

\begin{prop}  \label{Pbyinnerprod} Let $\lambda\in (\ZZ^n)^+$.
The bosonic Macdonald polynomial $P_\lambda$ is the unique element of
$\KK[x_1^{\pm1}, \ldots, x_n^{\pm1}]^{W_{\mathrm{fin}}}$ such that
\begin{enumerate}
\item[(a)] $P_\lambda = m_\lambda + (\hbox{lower terms})$;
\item[(b)] If $\gamma\in (\ZZ^n)^+$ and $\gamma < \lambda$ then $(P_\lambda, m_\gamma)_+ = 0$.
\end{enumerate}
\end{prop}

\subsection{Adjoints and orthogonality}\label{ssec:adjoints}

 For a linear operator $M\colon \KK[X]\to \KK[X]$, the \emph{adjoint of $M$}
is the linear operator $M^*\colon \KK[X] \to \KK[X]$ determined by 
$$\left( Mf_1, f_2\right)_+ = \left( f_1, M^*f_2\right)_+,
\qquad\hbox{for $f_1,f_2\in \KK[X]$,}
$$
where the inner product on $\KK[X]$ is as defined in~\eqref{ipdef}.

The following Proposition computes the adjoints of operators on $\KK[X]$ which come from $\widetilde{H}_{\mathrm{loc}}$.

\begin{prop} Let $i\in \{1, \ldots, n\}$ and $k\in \{0,1, \ldots, n\}$.  Then,
as operators on $\CC[x^{\pm1}_1, \ldots, x^{\pm1}_n]$,
$$
x_i^* = x^{-1}_i, \qquad 
T^*_{\alpha_k} = T_{\alpha_k}^{-1}, \qquad
Y_i^* = Y_i^{-1}, \qquad
\xi^*_{s_k} =\frac{\kappa^X_{\alpha_k}} {\kappa^X_{-\alpha_k}} \xi_{s_k}.
$$
\end{prop}
\begin{proof}
Let $J = \dfrac{1}{\Delta_{CC}^+}$.

\smallskip\noindent
\begin{itemize}
\item[$\circ$] Adjoint of multiplication by $x_i$: 

$$\left( x_if,g\right)_+ = \mathrm{ct}(x_if \cdot \overline{g}\cdot J)
= \mathrm{ct}(f\cdot \overline{x_i^{-1}g}\cdot J) = \left( f, x_i^{-1}g \right)_+.$$

\item[$\circ$] Adjoint of $\xi_{s_k}$: 
With $\kappa_\beta$ as in~\eqref{kpdefn},
\begin{align*}
\left(\kappa_\beta^X\right)^* = \overline{\kappa_\beta^X} 
&= \frac{(1-t_{\beta}^{-\frac12}u_{\beta}^{-\frac12}X^{-\beta})
(1+t_{\beta}^{-\frac12}u_{\beta}^{\frac12}X^{-\beta})}
{(1-X^{-2\beta})}
\\
&= \frac{t_\beta^{-\frac12}u_\beta^{-\frac12}X^{-\beta} t_\beta^{-\frac12}u_\beta^{\frac12}X^{-\beta}
(1-t_{\beta}^{\frac12}u_{\beta}^{\frac12}X^{\beta})
(1+t_{\beta}^{\frac12}u_{\beta}^{-\frac12}X^{\beta} ) }
{X^{-2\beta}(1-X^{2\beta})}
= t_\beta^{-1} \kappa_\beta.
\end{align*}
Since $s_k S^+ = (S^+-\{\alpha_k\}) \cup\{-\alpha_k\}$, we have
\begin{equation}
\xi_{s_k} \Delta^+_{CC} = \xi_{s_k}\left(\prod_{\beta\in S^+} \kappa_\beta\right) 
= \left( \prod_{\beta\in S^+-\{\alpha_k\}} \kappa_\beta\right)\kappa_{-\alpha_k}
=\Delta^+_{CC}\frac{\kappa_{-\alpha_k}}{\kappa_{\alpha_k}}.
\label{skonD}
%\tag{skonD}
\end{equation}
Using $\mathrm{ct}(\xi_{s_k}g) = \mathrm{ct}(g)$ and the formula from~\eqref{skonD},
\begin{align*}
\left( \xi_{s_k}f,g\right)_+
&= \mathrm{ct}\left( (\xi_{s_k}f)\overline{g} J  \right)
= \mathrm{ct}\left( \xi_{s_k}(f(\xi_{s_k}(\overline{g} J ))) \right)
= \mathrm{ct}\left( f(\xi_{s_k}(\overline{g} J )) \right)
\\
&= \mathrm{ct} \left(f(\xi_{s_k}\overline{g}) J \frac{\kappa_{\alpha_k} }{ \kappa_{-\alpha_k} } \right)
= \mathrm{ct} \left( f \overline{\frac{\kappa_{\alpha_k} }{ \kappa_{-\alpha_k} } (\xi_{s_k}g)} J \right)
=\left(f, \frac{\kappa_{\alpha_k} }{ \kappa_{-\alpha_k} } (\xi_{s_k} g)\,\right)_+.
\end{align*}
\item[$\circ$]  Adjoint of $T_{\alpha_k}$:  Using the formula for $T_{\alpha_k}$ in~\eqref{TdefH} and recalling that $(\ ,\ )_+$ is sesquilinear with respect to the involution $\overline{\phantom{T}}\colon \KK\to \KK$,
\begin{align*}
(T_{\alpha_k})^*
&= \left(-t_{\alpha_k}^{-\frac12}+(1+\xi_{s_k})c_{\alpha_k}^{X^{-1}} \right)^* 
= \left(-t_{\alpha_k}^{-\frac12}+(1+\xi_{s_k})t_{\alpha_k}^{-\frac12} \kappa_{-\alpha_k} \right)^* 
= -t_{\alpha_k}^{\frac12} + t_{\alpha_k}^{\frac12}\overline{\kappa_{-\alpha_k}}  (1+\xi_{s_k}^*)
\\
&= -t_{\alpha_k}^{\frac12} + t_{\alpha_k}^{\frac12}t_{\alpha_k}^{-1} \kappa_{-\alpha_k} 
\left(1+\frac{\kappa_{\alpha_k}}{\kappa_{-\alpha_k}} \xi_{s_k}\right) 
= -t_{\alpha_k}^{\frac12} + t_{\alpha_k}^{-\frac12} 
\left(\kappa_{-\alpha_k}+\kappa_{\alpha_k} \xi_{s_k}\right) 
\\
&=-t_{\alpha_k}^{\frac12} + (1+\xi_{s_k})c_{-\alpha_k}^X = T^{-1}_{\alpha_k}.
\end{align*}

\item[$\circ$]  Adjoint of $Y_j$: 
$$Y_1^* = (T_0T_1\cdots T_n \cdots T_1)^* = T_1^{-1}\cdots T^{-1}_n\cdots T_1^{-1} T_0^{-1} = 
(T_0 T_1\cdots T_n \cdots T_1)^{-1} = Y^{-1}_1,$$
and
if $j\in \{2, \ldots, n\}$ then
$$Y_j^* = (T^{-1}_{j-1} Y_{j-1} T^{-1}_{j-1})^* = T_{j-1}Y^{-1}_{j-1} T_{j-1} = 
(T^{-1}_{j-1} Y_{j-1} T^{-1}_{j-1})^{-1} = Y^{-1}_j.$$
\end{itemize}
\end{proof}

Since $T_j^{-1}\varepsilon_\pm^* = T_j^*\varepsilon_\pm^* = (\varepsilon_\pm T_j)^*$ for $j\in \{1, \ldots, n\}$
then $T_n^{-1}\varepsilon_\pm^* = (-t_n)^{\frac12}\varepsilon_\pm^*$
and $T_i^{-1}\varepsilon_\pm^* = t^{-\frac12}\varepsilon_\pm^*$ for $i\in \{1, \ldots, n-1\}$.  Since
$$\varepsilon_\pm^* = T_{w_0}^{-1} + (\hbox{lower terms}) = T_{w_0} + (\hbox{lower terms})
\qquad\hbox{then}\qquad
\varepsilon_\pm^* = \varepsilon_\pm.
$$
A similar argument applies to the other symmetrizers to show that 
\begin{equation}
\varepsilon_\Xi^* = \varepsilon_\Xi,
\qquad\hbox{for \quad $\Xi\in \{+,\pm,\mp,-\}$.}
\label{symmadj}
%\tag{symmadj}
\end{equation}

The relations $Y_i^*= Y^{-1}_i$ in combination with the knowledge of the eigenvalues for
the action of the $Y_i$ on the $E_\mu$ give the following orthogonality relations
for Macdonald polynomials.  The proof is exactly as in~\cite[(5.7.11)]{Mac03} and~\cite[Prop,\ 7.2]{CR22}.

\begin{prop} \label{orthg}
Let $(\ZZ^n_{\ge 0})^+ = \{ \lambda = (\lambda_1, \ldots, \lambda_n)\ |\ \lambda_1\ge \cdots \ge \lambda_n\ge 0\}$
and, for $\lambda\in (\ZZ^n_{\ge 0})^+$, 
let $P_\lambda$, $A^\pm_{\lambda+\omega}$, $A^\mp_{\lambda+\pi}$ and $A_{\lambda+\rho}$ be as defined in 
\eqref{PAdef} and~\eqref{Mdef}.
\item[(a)] Let $\lambda, \mu\in \ZZ^n$.  If $\mu \ne \lambda$ then 
$(E_\lambda, E_\mu)_+ = 0$.
\item[(b)] Let $\lambda, \mu\in (\ZZ^n_{\ge0})^+$.  If $\mu \ne \lambda$ then
\begin{align*}
(P_\lambda, P_\mu)_+ &= 0,
&(A^\pm_{\lambda+\omega}, A^\pm_{\mu+\omega})_+ &= 0, \\
(A_{\lambda+\rho}, A_{\mu+\rho})_+ &= 0, 
&(A^\mp_{\lambda+\pi}, A^\mp_{\mu+\pi})_+ &= 0.
\end{align*}
\end{prop}

\subsection{Going up a level}\label{ssec:levelup}

We describe four slightly different collections of 5 parameters by the brief notations
\begin{align*}
t^+ &= (t,t_0^{\frac12},u_0^{\frac12}, t_n^{\frac12}, u_n^{\frac12}),
&t^\pm &= (t,t_0^{\frac12},u_0^{\frac12}, q t_n^{\frac12}, u_n^{\frac12}),
\\
t^- &= (qt,t_0^{\frac12},u_0^{\frac12}, q t_n^{\frac12},  u_n^{\frac12}),
&t^\mp &= (qt,t_0^{\frac12},u_0^{\frac12}, t_n^{\frac12}, u_n^{\frac12}),
\end{align*}
and define
\begin{align*}
( f,g)_+ &= (f,g)_{q,t^+}  
&( f,g)_{\pm} &= (f,g)_{q,t^\pm} 
\\
( f,g)_{-} &= (f,g)_{q,t^-}
&(f,g)_{\mp} &= (f,g)_{q,t^\mp}. 
\end{align*}

The following Proposition shows that the norms of polynomials in the fermionic, bosonic and mesonic spaces
can be computed as norms of symmetric polynomials, but \emph{with shifted parameters}.  Alternatively, 
in the world of norms for symmetric polynomials,
the shifted parameters are a residue arising from the effect of 
multiplying by the Weyl denominators $A_\rho$, $A^\pm_\omega$, $A^\mp_\pi$.

\begin{prop} \label{levelup}
(Going up a level) Let $f,g\in \KK[X]^{W_{\mathrm{fin}}}$ so that $f$ and $g$ are symmetric functions.  
If $P_0=1$ and $A^\pm_\omega$, $A^\mp_\pi$ and $A_\rho$ are the Weyl denominators defined in~\eqref{Wdenoms} then
\begin{align*}
(f,g)_+
&= \frac{W_0(t,t_n)}{W_0(t,t_n)} ( P_0 f, P_0 g )_+,
&(f,g)_\pm
&= \frac{W_0(t,qt_n)}{W_0(t,t_n^{-1})} ( A^{\pm}_\omega f, A^{\pm}_\omega g )_+,
\\
(f,g)_-
&= \frac{W_0(qt,q t_n)}{W_0(t^{-1},t_n^{-1})} (A_\rho f, A_\rho g)_+,
&(f,g)_\mp
&= \frac{W_0(qt,t_n)}{W_0(t^{-1},t_n)} ( A^{\mp}_\pi f,  A^{\mp}_\pi g )_+.
\end{align*}
\end{prop}
\begin{proof} 
Let $t^\pm$ be the 5-tuple of parameters
$t^\pm = (t, t_0^{\frac12}, u_0^{\frac12}, q t_n^{\frac12}, u_n^{\frac12})$ and let
$$
\Delta^\pm_{CC}
= \Delta^X_g(t_0^{\frac12}, u_0^{\frac12})
\Delta^{X^{-1}}_g(t_0^{\frac12}, u_0^{\frac12})
\Delta^X_s(t)\Delta^{X^{-1}}_s(t)
\Delta^X_d(qt_n^{\frac12}, u_n^{\frac12})
\Delta^{X^{-1}}_d(qt_n^{\frac12}, u_n^{\frac12})
\Delta_{0,s}^X(t)
\Delta_{0,d}^X(qt_n^{\frac12}, u_n^{\frac12}).
$$
Since
\begin{align*}
\Delta_d^X(t_n^{\frac12}, u_n^{\frac12})
&= \Delta_d^X(qt_n^{\frac12}, u_n^{\frac12}) \cdot
\left(\prod_{1\le i\le n} (1-qt_n^{\frac12}u_n^{\frac12}x_i)(1+qt_n^{\frac12}u_n^{-\frac12}x_i)\right)
\\
&= \Delta_d^X(qt_n^{\frac12}, u_n^{\frac12}) \cdot x^\omega A^\pm_\omega(x,qt_n^{\frac12}, u_n^{\frac12})
\\
\Delta_{0,d}^X(t_n^{\frac12}, u_n^{\frac12})
&= \Delta_{0,d}^X(q t_n^{\frac12}, u_n^{\frac12}) 
\cdot \left(\prod_{1\le i\le n} \frac{(1-t_n^{\frac12}u_n^{\frac12}x_i)(1+t_n^{\frac12}u_n^{-\frac12}x_i)}
{(1-qt_n^{\frac12}u_n^{\frac12}x_i)(1+qt_n^{\frac12}u_n^{-\frac12}x_i)}\right)
\\
&= \Delta_{0,d}^X(q t_n^{\frac12}, u_n^{\frac12}) 
\frac{x^\omega A^\pm_\omega(x, t_n^{\frac12},u_n^{\frac12})}{ x^\omega A^\pm_\omega(x, qt_n^{\frac12},u_n^{\frac12})}
\\
\Delta_d^{X^{-1}}(t_n^{\frac12}, u_n^{\frac12})
&= \Delta_d^{X^{-1}}(qt_n^{\frac12}, u_n^{\frac12}) 
\cdot\left(\prod_{1\le i\le n} (1-qt_n^{\frac12}u_n^{\frac12}x^{-1}_i)(1+qt_n^{\frac12}u_n^{-\frac12}x^{-1}_i)\right)
\\
&= \Delta_d^{X^{-1}}(qt_n^{\frac12}, u_n^{\frac12}) x^{-\omega} A^\pm_\omega (x^{-1},qt_n^{\frac12}, u_n^{\frac12})
\end{align*}
then
\begin{align*}
\frac{\Delta^+_{CC}}{\Delta^\pm_{CC}}
&= A^\pm_\omega(x, t_n^{\frac12}, u_n^{\frac12})  A^\pm_\omega(x^{-1}, q t_n^{\frac12}, u_n^{\frac12}). 
\end{align*}

Using Proposition~\ref{Poinprop} and that $\text{Card}(W_{\mathrm{fin}}) = 2^nn!$ gives that for $h\in \KK[x_1^{\pm1}, \ldots, x_n^{\pm1}]^{W_{\mathrm{fin}}}$, 
\begin{align*}
\mathrm{ct}(&h\kappa_{w_0}^{X^{-1}})
= \frac{1}{2^nn!} \mathrm{ct}\left(\sum_{w\in W_{\mathrm{fin}}} 
w(h\kappa_{w_0}^{X^{-1}}) \right)
= \frac{1}{2^nn!} \mathrm{ct}\left(h \sum_{w\in W_{\mathrm{fin}}} w(\kappa_{w_0}^{X^{-1}}) \right)
= \frac{W_0(t,t_n)}{2^nn!} \mathrm{ct}(h).
\end{align*}
Let $\Delta_\infty^X = \Delta_g^X \Delta_s^X\Delta_d^X$ and $\Delta_0^X = \Delta_{0,s}(t)\Delta_{0,d}^X$ so that
$$\Delta^+_{CC} = \Delta_\infty^X \Delta_\infty^{X^{-1}}\Delta_0^X
\qquad\hbox{and}\qquad
\Delta_0^X = \frac{A^\pm_\omega(x,t_n^{\frac12}u_n^{\frac12}) }{a^\pm_\omega(x)}
\frac{A^\mp_\pi(x,t)}{a^\mp_\pi(x)}.
$$
Then
\begin{align*}
(A^\pm_\omega f, A^\pm_\omega g)_+ 
&= \mathrm{ct}\left(  \frac{f\overline{g}}{\Delta^+_{CC} }
A^\pm_\omega(x,t_n^{\frac12}, u_n^{\frac12})  \overline{A^\pm_\omega(x,t_n^{\frac12}, u_n^{\frac12})}  \right)
\\
&= \mathrm{ct}\left( \frac{f\overline{g}}{\Delta^X_\infty \Delta^{X^{-1}}_\infty}
\frac{ A^\pm_\omega(x,t_n^{\frac12},u_n^{\frac12})}{\Delta^X_0} A^\pm_\omega(x^{-1}, t_n^{-\frac12}, u_n^{-\frac12}) 
\right)
\\
&= \mathrm{ct}\left( \frac{f\overline{g}}{\Delta^X_\infty \Delta^{X^{-1}}_\infty }
\frac{ a^\pm_\omega(x)a^\mp_\pi(x) }{A^\mp_\pi(x,t)} A^\pm_\omega(x^{-1}, t_n^{-\frac12}, u_n^{-\frac12})  \right)
\\
&= \mathrm{ct}\left( \frac{f\overline{g}}{\Delta^X_\infty \Delta^{X^{-1}}_\infty }
\frac{ a^\pm_\omega(x)a^\mp_\pi(x) a^\pm_\omega(x^{-1})a^\mp_\pi(x^{-1}) }{ A^\mp_\pi(x,t) A^\mp_\pi(x^{-1},t) } 
\cdot
\frac{ A^\pm_\omega(x^{-1}, t_n^{-\frac12}, u_n^{-\frac12}) A^\mp_\pi(x^{-1},t) }{ a^\pm_\omega(x^{-1})a^\mp_\pi(x^{-1}) }  \right)
\\
&= \mathrm{ct}\left( \frac{f\overline{g}}{\Delta^X_\infty \Delta^{X^{-1}}_\infty }
\frac{ a^\pm_\omega(x)a^\mp_\pi(x) a^\pm_\omega(x^{-1})a^\mp_\pi(x^{-1}) }{ A^\mp_\pi(x,t) A^\mp_\pi(x^{-1},t) } 
\cdot \kappa^{X^{-1}}_{w_0}(t,t_n^{-\frac12}, u_n^{-\frac12})
\right)
\\
&= \frac{W_0(t,t_n^{-1})}{2^n n!}
\mathrm{ct}\left( \frac{f\overline{g}}{\Delta^X_\infty \Delta^{X^{-1}}_\infty}
\frac{ a^\pm_\omega(x)a^\mp_\pi(x) a^\pm_\omega(x^{-1})a^\mp_\pi(x^{-1}) }{ A^\mp_\pi(x,t) A^\mp_\pi(x^{-1},t) } 
\right)
\end{align*}
and
\begin{align*}
(f,g)_{\pm}
&= \mathrm{ct}\left( \frac{ f\overline{g} }{\Delta^\pm_{CC} }\right) 
= \mathrm{ct}\left( \frac{f\overline{g} }{ \Delta^+_{CC} } 
A^\pm_\omega(x,t_n^{\frac12},u_n^{\frac12}) A^\pm_\omega(x^{-1}, qt_n^{\frac12},u_n^{\frac12})  \right) 
\\
&= \mathrm{ct}\left( \frac{f\overline{g} }{ \Delta^X_\infty \Delta^{X^{-1}}_\infty } 
\frac{ A^\pm_\omega(x,t_n^{\frac12},u_n^{\frac12})}{ \Delta^X_0 } A^\pm_\omega(x^{-1}, qt_n^{\frac12},u_n^{\frac12}) \right) 
\\
&= \mathrm{ct}\left( \frac{f\overline{g} }{ \Delta^X_\infty \Delta^{X^{-1}}_\infty } 
\frac{ a^\pm_\omega(x)a^\mp_\pi(x) }{ A^\mp_\pi(x,t) } 
A^\pm_\omega(x^{-1}, qt_n^{\frac12},u_n^{\frac12}) \right) 
\\
&= \mathrm{ct}\left( \frac{f\overline{g} }{ \Delta^X_\infty \Delta^{X^{-1}}_\infty } 
\frac{ a^\pm_\omega(x)a^\mp_\pi(x) a^\pm_\omega(x^{-1})a^\mp_\pi(x^{-1}) }{ A^\mp_\pi(x,t) A^\mp_\pi(x^{-1},t) } 
\cdot
\frac{A^\pm_\omega(x^{-1}, qt_n^{\frac12},u_n^{\frac12})A^\mp_\pi(x^{-1},t) }{a^\pm_\omega(x^{-1})a^\mp_\pi(x^{-1})} \right) 
\\
&= \mathrm{ct}\left( \frac{f\overline{g} }{ \Delta^X_\infty \Delta^{X^{-1}}_\infty } 
\frac{ a^\pm_\omega(x)a^\mp_\pi(x) a^\pm_\omega(x^{-1})a^\mp_\pi(x^{-1}) }{ A^\mp_\pi(x,t) A^\mp_\pi(x^{-1},t) } 
\cdot \kappa^{X^{-1}}_{w_0}(t,qt_n^{\frac12},u_n^{\frac12}) \right)
\\
&= \frac{ W_0(t,qt_n) }{ 2^n n!}
\mathrm{ct}\left( \frac{f\overline{g} }{ \Delta^X_\infty \Delta^{X^{-1}}_\infty } 
\frac{ a^\pm_\omega(x)a^\mp_\pi(x) a^\pm_\omega(x^{-1})a^\mp_\pi(x^{-1}) }{ A^\mp_\pi(x,t) A^\mp_\pi(x^{-1},t) } 
\right).
\end{align*}
which proves the $\pm$ case.  The other cases are similar.
\end{proof}

\subsection{Weyl character formulas}\label{ssec:WCF}

As in~\eqref{rhodef}, let $0,\omega, \pi,\rho\in \ZZ^n$ by $0=(0,0,\ldots, 0)$,
$$
\omega = (1,1,\ldots, 1),\qquad
\pi = (n-1,n-2,\ldots, 2,1,0), \qquad
\rho = (n,n-1,\ldots, 3,2,1).
$$

\begin{thm} \label{WCFs} (Weyl character formulas)
Let $\lambda\in \ZZ^n$ with $\lambda_1\ge \lambda_2\ge \cdots \ge \lambda_n >0$.
Then
\begin{align*}
P_\lambda(q,t^+) &= \frac{P_{\lambda+0}(q,t^+)}{P_0(q,t^+)},
&P_\lambda(q,t^\pm) &= \frac{A^\pm_{\lambda+\omega}(q,t^+)}{A^{\omega}_\pi(q,t^+)},
\\
P_\lambda(q,t^-) &= \frac{A_{\lambda+\rho}(q,t^+)}{A_\rho(q,t^+)},
&P_\lambda(q,t^\mp) &= \frac{A^\mp_{\lambda+\pi}(q,t^+)}{A^{\mp}_\pi(q,t^+)},
\end{align*}
\end{thm}
\begin{proof}
The following is the proof for the case $P_\lambda(q,t^\pm)$.  The proof for the other cases is similar.
Since $A^\pm_{\lambda+\omega} = t^{\frac12\ell(w_0)} \varepsilon_\pm E_{\lambda+\omega}$ then 
$A^\pm_{\lambda+\omega}\in \CC[X]^{\mathrm{Mes}\pm}$.  
Thus, by Proposition~\ref{BosFer}, 
$$\hbox{there exists}\qquad f\in \KK[X]^{W_{\mathrm{fin}}}\qquad\hbox{such that}\qquad
A^\pm_{\lambda+\omega} = A^\pm_\omega f.
$$
Since $\lambda+\omega$ is dominant then $w_0(\lambda+\omega) = -(\lambda+\omega)$ is antidominant
(weakly increasing with all entries $\le 0$).
If $\mu\in \ZZ^n$ is such that the coefficient of $x^\mu$ in $A^\pm_{\lambda+\omega}$ is nonzero then 
$\mu\le -(\lambda+\omega)$ in the DBlex order.
Thus, using the notations of Section~\ref{sec:EPchar},
$$f = m_\lambda + (\hbox{lower terms}).$$
The $E$-expansion for $A^\pm_{\lambda+\omega}$ in Proposition~\ref{Eexps} gives that 
$$A^\pm_{\lambda+\omega} = \sum_{\mu\in W_{\mathrm{fin}}(\lambda+\pi)} d_{\lambda+\omega}^\mu E_\mu
= E_{w_0(\lambda+\omega)} + (\hbox{lower terms})
$$
and, from the definitions of $A^\pm_\omega$ and $m_\nu$, 
$$A^\pm_\omega m_\nu = x^{w_0(\nu+\omega)} + (\hbox{lower terms}).$$
Since $(E_{w_0(\lambda+\omega)}, x^\gamma)_+ = 0$ for $\gamma<w_0(\lambda+\omega)$ then
$$(A^\pm_\omega f, A^\pm_\omega m_\nu)_+ = (A^\pm_{\lambda+\omega}, A^\pm_\omega m_\nu)_+ = 0,
\qquad\hbox{for $\nu\in (\ZZ^n)^+$ with $\nu<\lambda$.}
$$
Using Proposition~\ref{levelup}, since $f\in \KK[X]^{W_{\mathrm{fn}}}$ and $m_\nu\in \KK[X]^{W_{\mathrm{fin}}}$ then 
$$(f, m_\nu)_\pm = \frac{ W_0(t, qt_n) }{ W_0(t, t_n^{-1})} (A^\pm_\omega f, A^\pm_\omega m_\nu)_+ = 0,
\qquad\hbox{for $\nu\in (\ZZ^n)^+$ with $\nu<\lambda$.}
$$
Thus, by Proposition~\ref{Pbyinnerprod}, $f = P_\lambda(q,t^\pm)$.
\end{proof}

\subsection{Reductions for norms}\label{ssec:normred}

The following proposition shows that the $c$-functions provide an effective framework for
describing the differences between the norms of the various flavors of Koornwinder polynomials.

\begin{prop} \label{normred}
\item[(a)] Let $\mu = (\mu_1, \ldots, \mu_n)\in \ZZ^n$.  Then
$$\frac{(E_\mu, E_\mu)_+}{(1,1)_+ } = \ev_0^t(c_{u_\mu}^Y c_{u_\mu}^{Y^{-1}}).$$
\item[(b)] Let $\lambda\in (\ZZ_{\ge 0}^n)^+$.
\begin{align*}
\frac{ (P_\lambda, P_\lambda)_+ }{ (E_\lambda, E_\lambda)_+  }
&= \frac{W_0(t,t_n)}{W_\lambda(t,t_n)} \ev^t_\lambda(\kappa_{v_\lambda}^+),
\\
\frac{(A^\pm_{\lambda+\omega},A^\pm_{\lambda+\omega})_+}{(E_{\lambda+\omega},E_{\lambda+\omega})_+ }
&=\frac{W_0(t,t_n^{-1})}{W_{\lambda+\omega}(t,t_n)} t_n^n \ev^t_{\lambda+\omega}(\kappa^\pm_{v_{\lambda+\omega}}),
\\
\frac{(A^\mp_{\lambda+\pi},A^\mp_{\lambda+\pi})_+}{(E_{\lambda+\pi},E_{\lambda+\pi})_+ }
&=\frac{W_0(t^{-1},t_n)}{W_{\lambda+\pi}(t,t_n)} t^{n(n-1)} \ev^t_{\lambda+\pi}(\kappa^\mp_{v_{\lambda+\pi}}),
\\
\frac{(A_{\lambda+\rho},A_{\lambda+\rho})_+}{(E_{\lambda+\rho},E_{\lambda+\rho})_+ }
&=\frac{W_0(t^{-1},t_n^{-1})}{W_{\lambda+\rho}(t,t_n)} t^{n(n-1)}t_n^n \ev^t_{\lambda+\rho}(\kappa^-_{v_{\lambda+\rho}}).
\end{align*}
\item[(c)] Let $\lambda\in (\ZZ_{\ge 0}^n)^+$.
\begin{align*}
\frac{(A^\pm_{\lambda+\omega},A^\pm_{\lambda+\omega})_+}{(P_{\lambda+\omega},P_{\lambda+\omega})_+ }
&=\frac{W_0(t,t_n^{-1}) }{ W_0(t,t_n)}
 t_n^n \ev^t_{\lambda+\omega}\left(\frac{ \kappa^\pm_{v_{\lambda+\omega}} }{\kappa^+_{v_{\lambda+\omega}}} \right),
\\
\frac{(A^\mp_{\lambda+\pi},A^\mp_{\lambda+\pi})_+}{(P_{\lambda+\pi},P_{\lambda+\pi})_+ }
&=\frac{W_0(t^{-1},t_n)}{W_0(t,t_n)} 
t^{n(n-1)} \ev^t_{\lambda+\pi}\left(\frac{ \kappa^\mp_{v_{\lambda+\pi}} }{\kappa^+_{v_{\lambda+\pi}}} \right),
\\
\frac{ (A_{\lambda+\rho}, A_{\lambda+\rho})_+ }{ (P_{\lambda+\rho}, P_{\lambda+\rho})_+  }
&=\frac{W_0(t^{-1},t_n^{-1})}{W_0(t,t_n)} 
t^{n(n-1)}t_n^n \ev_{\lambda+\rho}^t\left( \frac{ \kappa_{w_0}^- }{\kappa_{w_0}^+} \right).
\end{align*}
\item[(d)] Let $\lambda\in (\ZZ_{\ge 0}^n)^+$.  Then
\begin{align*}
\frac{ (P_\lambda(q,t^\pm), P_\lambda(q,t^\pm) )_\pm }{ (P_{\lambda+\omega}(q,t^+), P_{\lambda+\omega}(q,t^+))_+ }
&= \frac{W_0(t,qt_n)}{W_0(t,t_n)} 
t_n^{\ell_d(w_0)} \ev^t_{\lambda+\omega}\left(\frac{\kappa^\pm_{v_{\lambda+\omega}} }{ \kappa^+_{v_{\lambda+\omega}} }\right)
\\
\frac{ (P_\lambda(q,t^\mp), P_\lambda(q,t^\mp) )_\mp }{ (P_{\lambda+\pi}(q,t^+), P_{\lambda+\pi}(q,t^+))_+ }
&= \frac{W_0(qt,t_n)}{W_0(t,t_n)} 
t^{\ell_s(w_0)} \ev^t_{\lambda+\pi}\left(\frac{\kappa^\mp_{v_{\lambda+\pi}} }{ \kappa^+_{v_{\lambda+\pi}} }\right)
\\
\frac{ (P_\lambda(q,t^-), P_\lambda(q,t^-) )_- }{ ( P_{\lambda+\rho}(q,t^+), P_{\lambda+\rho}(q,t^+) )_+ }
&= \frac{W_0(qt,qt_n)}{W_0(t,t_n)} t^{\ell_s(w_0)}t_n^{\ell_d(w_0)} \ev_{\lambda+\rho}^t \left( \frac{ \kappa_{w_0}^- }{ \kappa_{w_0}^+ }\right).
\end{align*}
\end{prop}
\begin{proof}
(a) Using the creation formula for $E_\mu$,
\begin{align*}
(E_\mu, E_\mu)_+
&= (t^{-\frac12\ell(v_\mu^{-1})} \tau^\vee_{u_\mu} \mathbf{1}_Y, 
t^{-\frac12\ell(v_\mu^{-1})} \tau^\vee_{u_\mu} \mathbf{1}_Y)_+=
( \tau^\vee_{u_\mu^{-1}} \tau^\vee_{u_\mu} \mathbf{1}_Y, \mathbf{1}_Y)_+
\\
&= ( c_{u_\mu}^Y c_{u_\mu}^{Y^{-1}} \mathbf{1}_Y, \mathbf{1}_Y)_+
= \ev_0^t(c_{u_\mu}^Y c_{u_\mu}^{Y^{-1}})\cdot (1,1)_+.
\end{align*}
(b) Using Proposition~\ref{Eexps} ($E$-expansion formulas),
$$A^\pm_{\lambda+\omega} = \frac{\chi^\pm(T_{w_0}) }{W_{\lambda+\omega}(t,t_n)} \varepsilon_\pm  E_{\lambda+\omega} 
= \sum_{\mu\in W_{\mathrm{fin}} (\lambda+\omega)} d_{\lambda+\omega}^\mu E_\mu, 
\qquad \hbox{with}\qquad 
d_{\lambda+\omega}^{\lambda+\omega} = (-1)^{\ell_d(v_{\lambda+\omega})} \ev_{\lambda+\omega}^t(\kappa^\pm_{v_{\lambda+\omega}}),
$$
and~\eqref{symmsq}
$$\varepsilon_\pm^2 = \frac{1}{\chi^{\pm}(T_{w_0})} W_0(t,t_n^{-1})\varepsilon_\pm,$$
gives
\begin{align*}
(A^\pm_{\lambda+\omega}, A^\pm_{\lambda+\omega})_+
&= \left(\frac{ \chi^+(T_{w_0}) }{W_{\lambda+\omega}(t,t_n)} \varepsilon_\pm  E_{\lambda+\omega}, 
\frac{ \chi^+(T_{w_0}) }{W_{\lambda+\omega}(t,t_n)} \varepsilon_\pm  E_{\lambda+\omega} \right)_+
\qquad\hbox{(by~\eqref{Mdef})}
\\
&= \frac{1}{W_{\lambda+\omega}(t,t_n) W_{\lambda+\omega} (t^{-1}, t_n^{-1})} 
\left( \varepsilon_\pm^2 E_{\lambda+\omega},  E_{\lambda+\omega} \right)_+
\qquad\hbox{(by~\eqref{symmadj})}
\\
&= \frac{1}{\chi^\pm(T_{w_0})}  \frac{ W_0(t,t_n^{-1}) }
{W_{\lambda+\omega}(t,t_n)W_{\lambda+\omega}(t^{-1}, t_n^{-1}) } 
\left( \varepsilon_\pm E_{\lambda+\omega},  E_{\lambda+\omega} \right)_+
\qquad\hbox{(by~\eqref{symmsq})}
\\
&= \frac{\chi^+(T_{w_{\lambda+\omega}})^2}{\chi^\pm(T_{w_0}) \chi^+(T_{w_0}) }  \frac{ W_0(t,t_n^{-1}) }{W_{\lambda+\omega} (t, t_n)} 
\left( A^\pm_{\lambda+\omega},  E_{\lambda+\omega} \right)_+
\qquad\hbox{(by~\eqref{Mdef})}
\\
&= \frac{t^{\ell_d(w_0)} }{ (-1)^{\ell_d(w_0)} }
\frac{W_0(t,t_n^{-1})}{ W_{\lambda+\omega}(t^{-1},t_n^{-1})}
d_{\lambda+\omega}^{\lambda+\omega} (E_{\lambda+\omega}, E_{\lambda+\omega})_+
\\
&= \frac{t_n^n }{ (-1)^{\ell_d(w_0)}} 
\frac{W_0(t,t_n^{-1})}{W_{\lambda+\omega}(t,t_n)} 
(-1)^{\ell_d(v_{\lambda+\omega})}\ev_{\lambda+\omega}^t(\kappa^\pm_{v_{\lambda+\omega}})
(E_{\lambda+\omega}, E_{\lambda+\omega})_+
\\
&=
\frac{W_0(t,t_n^{-1})}{W_{\lambda+\omega}(t,t_n)} t_n^n
\ev_{\lambda+\omega}^t(\kappa^\pm_{v_{\lambda+\omega}})
(E_{\lambda+\omega}, E_{\lambda+\omega})_+,
\end{align*}
because $(-1)^{\ell_d(v_{\lambda+\omega})} = (-1)^{\ell_d(w_0)}$.
This proves the $\pm$ case.  The proof for the other cases is similar.

\medskip\noindent
(c) By part (b),
$$\frac{(A^\pm_{\lambda+\omega},A^\pm_{\lambda+\omega})_+}
{(E_{\lambda+\omega},E_{\lambda+\omega})_+}
=\frac{W_0(t,t_n^{-1})}{W_{\lambda+\omega}(t,t_n)} t_n^n
\ev_{\lambda+\omega}^t(\kappa^\pm_{v_{\lambda+\omega}})
\quad\hbox{and}\quad
\frac{(P_{\lambda+\omega}, P_{\lambda+\omega})_+}
{(E_{\lambda+\omega},E_{\lambda+\omega})_+}
=\frac{W_0(t,t_n)}{W_{\lambda+\omega}(t,t_n)} 
\ev_{\lambda+\omega}^t(\kappa^+_{v_{\lambda+\omega}})
$$
gives the relation between $(A^\pm_{\lambda+\omega},A^\pm_{\lambda+\omega})_+$ and
$(P_{\lambda+\omega}, P_{\lambda+\omega})_+$.
This proves the $\pm$ case.  The proof for the other cases is similar.

\smallskip\noindent
(d)  Using Proposition~\ref{levelup} (going up a level) and Proposition~\ref{WCFs} (Weyl character formula) gives
\begin{align*}
(P_\lambda(q,t^\pm), P_\lambda(q,t^\pm) )_\pm
&= \frac{W_0(t,qt_n)}{W_0(t,t_n^{-1})} (A^\pm_\omega P_\lambda(q,t^\pm), A^\pm_\omega P_\lambda(q,t^\pm))_+ 
\\
&= \frac{W_0(t,qt_n)}{W_0(t,t_n^{-1})} (A^\pm_{\lambda+\omega}(q,t^+), A^\pm_{\lambda+\omega}(q,t^+))_+ 
\\
&= \frac{W_0(t,qt_n)}{W_0(t,t_n)} 
t_n^n \ev^t_{\lambda+\omega}\left(\frac{\kappa^\pm_{v_{\lambda+\omega}} }{ \kappa^+_{v_{\lambda+\omega}} }\right)
(P_{\lambda+\omega}(q,t^+), P_{\lambda+\omega}(q,t^+))_+ 
\end{align*}
This proves the $\pm$ case.  The proof for the other cases is similar.
\end{proof}

\subsection{The symmetric inner product}

In this subsection we define, for symmetric functions,
 a slightly different inner product $\langle\ ,\ \rangle_+$ that has more symmetry
than the original inner product $(\ ,\ )_+$.  Fortunately, the difference between these inner products is
only a factor of $W_0(t,t_n)$, which makes the reduction relations of the previous subsection even simpler.
The symmetric inner product $\langle\ ,\ \rangle_+$
 has another very useful advantage: in terms of the Askey-Wilson parameters $a,b,c,d$, the inner product
$\langle\ ,\ \rangle_+$ is completely symmetric in $a,b,c$ and $d$.

\smallskip
Define involutions
$\overline{\phantom{T}}\colon \KK[X] \to \KK[X]$ and
${}^\sigma\colon \KK[X] \to \KK[X]$ and 
${}^t \colon \KK[X] \to \KK[X]$ by
\begin{align*}
\bar f(x_1, \ldots, x_n;q,t,t_0^{\frac12}, u_0^{\frac12}, t_n^{\frac12}, u_n^{\frac12}) 
&= f(x_1^{-1}, \ldots, x_n^{-1};q^{-1}, t^{-1}, t_0^{-\frac12}, u_0^{-\frac12}, t_n^{-\frac12}, u_n^{-\frac12}),
\\
f^\sigma(x_1, \ldots, x_n;q,t,t_0^{\frac12}, u_0^{\frac12}, t_n^{\frac12}, u_n^{\frac12}) 
&= f(x_1^{-1}, \ldots, x_n^{-1};q, t, t_0^{\frac12}, u_0^{\frac12}, t_n^{\frac12}, u_n^{\frac12}),
\\
f^t(x_1, \ldots, x_n;q,t,t_0^{\frac12}, u_0^{\frac12}, t_n^{\frac12}, u_n^{\frac12}) 
&= f(x_1, \ldots, x_n;q^{-1}, t^{-1}, t_0^{-\frac12}, u_0^{-\frac12}, t_n^{-\frac12}, u_n^{-\frac12}).
\end{align*}

Let $\kappa_{w_0}^{X^{-1}}$ be as given in~\eqref{kw0X}
and define
\begin{equation}
\nabla_{CC} = \Delta^+_{CC} \kappa_{w_0}^{X^{-1}}.
\label{nabladef}
%\tag{nabladef}
\end{equation}
Then define a new scalar product
$\langle\ ,\ \rangle_+\colon \KK[X]\times \KK[X]\to \KK$ by 
\begin{equation}
\langle f_1,f_2\rangle_+ = \frac{1}{2^n n!} \mathrm{ct}\left(f_1 \overline{f_2} \frac{1}{\nabla_{CC}}\right),
\qquad
\hbox{where\quad $2^n n! = \Card(W_{\mathrm{fin}})$.}
\label{syminp}
%\tag{syminp}
\end{equation}
The following result provides a comparison of $(\ ,\ )_+$ and $\langle\ ,\ \rangle_+$ as inner products on symmetric polynomials.

\begin{prop}
Let $f,g\in \KK[x_1^{\pm}, \ldots, x_n^{\pm1}]^{W_{\mathrm{fin}}}$.  Then
$$\langle f, g\rangle_+ = \frac{1}{W_0(t,t_n)} (f,g^t)_+.$$
\end{prop}
\begin{proof}
Let $f,g \in  \KK[x_1^{\pm}, \ldots, x_n^{\pm1}]^{W_{\mathrm{fin}}}$.  
Since $f$, $\bar g^t$ and $\nabla_{CC}$ are all invariant under the action of $W_{\mathrm{fin}}$ then
\begin{align*}
\langle f,g\rangle_+ 
&= \frac{1}{2^n n!} \mathrm{ct}\left( f g^\sigma \frac{1}{\nabla_{CC}}\right)
= \frac{1}{2^n n!\cdot W_0(t,t_n)} \mathrm{ct}\left( f \bar g^t \frac{1}{\nabla_{CC} } W_0(t,t_n) \right)
\\
&= \frac{1}{2^n n!\cdot W_0(t,t_n)} \mathrm{ct}\left( f \bar g^t \frac{1}{\nabla_{CC} }
\left(\sum_{w\in W_{\mathrm{fin}}} w(\kappa_{w_0}^{X^{-1}})\right) \right)
\\
&= \frac{1}{2^n n! \cdot W_0(t,t_n)} \mathrm{ct}\left( \sum_{w\in W_{\mathrm{fin}}} w\left(
f \bar g^t \frac{1}{\nabla_{CC} }
\kappa_{w_0}^{X^{-1}}\right) \right)
\\
&= \frac{1}{2^n n!\cdot W_0(t,t_n)} \mathrm{ct}\left( \sum_{w\in W_{\mathrm{fin}}} w\left(
f \bar g^t \frac{1}{\Delta^+_{CC} }
\right)\right)
\\
&= \frac{1}{W_0(t,t_n)} \mathrm{ct}\left( f \bar g^t \frac{1}{\Delta^+_{CC} } \right)
= \frac{1}{W_0(t,t_n)}(f,g^t)_+,
\end{align*}
where the third equality uses Proposition~\ref{Poinprop}.
\end{proof}

The following corollary of Proposition~\ref{normred} records the norm comparisons
in terms of the symmetric inner product $\langle\ ,\ \rangle_+$.

\begin{cor} \label{symmnormred}
Let $\lambda\in (\ZZ_{\ge 0}^n)^+$.  Then
\begin{align*}
\frac{ \langle P_\lambda(q,t^\pm), P_\lambda(q,t^\pm) \rangle_\pm }
{ \langle P_{\lambda+\omega}(q,t^+), P_{\lambda+\omega}(q,t^+)\rangle_+ }
&=  t_n^{\ell_d(w_0)} \ev^t_{\lambda+\omega}\left(\frac{\kappa^\pm_{v_{\lambda+\omega}} }{ \kappa^+_{v_{\lambda+\omega}} }\right) \\
\frac{ \langle P_\lambda(q,t^\mp), P_\lambda(q,t^\mp) \rangle_\mp }{ \langle P_{\lambda+\pi}(q,t^+), P_{\lambda+\pi}(q,t^+)\rangle_+ }
&= t^{\ell_s(w_0)} \ev^t_{\lambda+\pi}\left(\frac{\kappa^\mp_{v_{\lambda+\pi}} }{ \kappa^+_{v_{\lambda+\pi}} }\right) \\
\frac{ \langle P_\lambda(q,t^-), P_\lambda(q,t^-) \rangle_- }{ \langle P_{\lambda+\rho}(q,t^+), P_{\lambda+\rho}(q,t^+) )_+ }
&= t^{\ell_s(w_0)}t_n^{\ell_d(w_0)} \ev_{\lambda+\rho}^t \left( \frac{ \kappa_{w_0}^- }{ \kappa_{w_0}^+ }\right).
\end{align*}
\end{cor}

\begin{remark} \label{nablasymm}
Using Remark~\ref{Deltaexp} then
$$\nabla_{CC} = \nabla^{(1)}\nabla^{(2)},
\qquad\hbox{where}\quad
\nabla^{(2)}
= \prod_{1\le i< j\le n} 
\frac{(tx_ix_j^{-1};q)_\infty (tx_ix_j;q)_\infty (tx_i^{-1}x_j;q)_\infty (tx_i^{-1}x_j^{-1};q)_\infty}
{ (x_ix_j^{-1};q)_\infty (x_ix_j;q)_\infty (x_i^{-1}x_j;q)_\infty (x_i^{-1}x_j^{-1};q)_\infty}
$$
and
\begin{equation*}
\nabla^{(1)}
= \prod_{i=1}^n \frac{
(cx_i;q)_\infty (c x_i^{-1};q)_\infty
(dx_i;q)_\infty (d x_i^{-1};q)_\infty
(ax_i;q)_\infty (a x_i^{-1};q)_\infty
(bx_i;q)_\infty (b x_i^{-1};q)_\infty
}
{ (x_i^2;q)_\infty (x_i^{-2};q)_\infty} 
.
\end{equation*}
This formula shows that $\nabla_{CC}$, and thus $\langle\ ,\ \rangle_+$ and $P_\lambda(q,t^+)$, 
are completely symmetric in the Askey-Wilson parameters $a,b,c$ and $d$.
\qed\end{remark}

\begin{remark}Using Remark~\ref{rem:bigvmu},
$$
\frac{\kappa^\pm_{v_{\lambda+\omega}} }{\kappa^+_{v_{\lambda+\omega}} } 
=\frac{\kappa^\pm_{v_\omega}}{\kappa^+_{v_\omega}} 
= \prod_{i=1}^n \frac{\kappa^{Y^{-1}}_{\varepsilon_i} }{\kappa^Y_{\varepsilon_i}}
=\frac{\kappa^\pm_{w_0}}{\kappa^+_{w_0}} 
\qquad\hbox{and}\qquad
\frac{\kappa^\mp_{v_{\lambda+\pi}} }{ \kappa^+_{v_{\lambda+\pi} } }
=\frac{\kappa^\mp_{v_\pi}}{\kappa^+_{v_\pi}} 
= \prod_{1\le i<j\le n} 
\frac{\kappa^{Y^{-1}}_{\varepsilon_i-\varepsilon_j}\kappa^{Y^{-1}}_{\varepsilon_i+\varepsilon_j}  }
{\kappa^Y_{\varepsilon_i-\varepsilon_j} \kappa^Y_{\varepsilon_i+\varepsilon_j} } 
= \frac{\kappa^\mp_{w_0}}{\kappa^+_{w_0}} 
$$
and
$$\frac{\kappa^\pm_{v_\omega}}{\kappa^+_{v_\omega}} 
\cdot \frac{\kappa^\mp_{v_\omega}}{\kappa^+_{v_\omega}} 
=\frac{\kappa^-_{w_0}}{\kappa^+_{w_0}} 
= \left(\prod_{i=1}^n \frac{\kappa^{Y^{-1}}_{\varepsilon_i} }{\kappa^Y_{\varepsilon_i}}\right)
\left(\prod_{1\le i<j\le n} 
\frac{\kappa^{Y^{-1}}_{\varepsilon_i-\varepsilon_j}\kappa^{Y^{-1}}_{\varepsilon_i+\varepsilon_j}  }
{\kappa^Y_{\varepsilon_i-\varepsilon_j} \kappa^Y_{\varepsilon_i+\varepsilon_j} } 
\right)
$$
\qed\end{remark}

\subsection{The norm formula}\label{ssec:normformula}

In this section we use the recursive relations of Corollary~\ref{symmnormred} to derive a $c$-function formula
for $\langle P_\lambda(q,t^+), P_\lambda(q,t^+)\rangle_+$.

Introduce $Y$-versions of the $c$-function products from Remark~\ref{Deltaexp} as follows:
\begin{align*}
\Delta^Y_s(t^+)\Delta_{0,s}^Y(t^+)
&= \prod_{1\le i<j\le n} \frac{(tY_iY_j^{-1};q)_\infty (tY_iY_j;q)_\infty 
}{
(Y_iY_j^{-1};q)_\infty(Y_iY_j;q)_\infty} 
\\
\Delta_{d}^Y(t^+)\Delta_{0,d}^Y(t^+)
&=\prod_{i=1}^n \frac{(t_n^{\frac12}t_0^{\frac12}Y_i;q)_\infty (-t_n^{\frac12}t_0^{-\frac12}Y_i; q)_\infty }
{(Y_i^2;q^2)_\infty}
\\
\Delta_g^{Y}(t^+)
&=\prod_{i=1}^n  \frac{(u_n^{\frac12}u_0^{\frac12}q^{\frac12}Y_i;q)_\infty (-u_n^{\frac12}u_0^{-\frac12}q^{\frac12}Y_i; q)_\infty }
{(qY_i^2;q^2)_\infty}
\\
\Delta^{Y^{-1}}_s((t^+)^{-1})
&=\prod_{1\le i<j\le n} \frac{(t^{-1}qY_i^{-1}Y_j;q)_\infty (t^{-1}qY_i^{-1}Y_j^{-1};q)_\infty }{(qY_i^{-1}Y_j;q)_\infty(qY_i^{-1}Y_j^{-1};q)_\infty }
\\
\Delta_d^{Y^{-1}}((t^+)^{-1})
&=\prod_{i=1}^n \frac{(t_n^{-\frac12}t_0^{-\frac12} qY_i^{-1};q)_\infty (-t_n^{-\frac12}t_0^{\frac12} q Y_i^{-1}; q)_\infty }
{(q^2Y_i^{-2};q^2)_\infty}
\\
\Delta_g^{Y^{-1}}((t^+)^{-1})
&=\prod_{i=1}^n \frac{(u_n^{-\frac12}u_0^{-\frac12} q^{\frac12}Y_i^{-1};q)_\infty (-u_n^{-\frac12}u_0^{\frac12} q^{\frac12} Y_i^{-1}; q)_\infty }
{(qY_i^{-2};q^2)_\infty}.
\end{align*}

Define homomorphisms 
$\ev_{q^\lambda t^\rho}\colon \KK[Y]\to \KK$ and 
$\ev_{q^{-\lambda} t^{-\rho}}\colon \KK[Y]\to \KK$ by
$$\ev_{q^\lambda t^\rho}(Y_i) = q^{\lambda_i}t^{n-i}t_0^{\frac12}t_n^{\frac12}
\qquad\hbox{and}\qquad
\ev_{q^{-\lambda} t^{-\rho}}(Y_i) = q^{-\lambda_i}t^{-(n-i)}t_0^{-\frac12}t_n^{-\frac12}
=\ev_{q^{\lambda} t^{\rho}}(Y^{-1}_i).
$$
The correspondence between the parameters $q,t,t_0^{\frac12}, u_0^{\frac12}, t_n^{\frac12}, u_n^{\frac12}$ and
the Askey-Wilson parameters is given by
\begin{equation*}
q=q, \quad
a = q^{\frac12}t_0^{\frac12}u_0^{\frac12}, \quad
b = - q^{\frac12}t_0^{\frac12}u_0^{-\frac12}, \quad
c = t_n^{\frac12}u_n^{\frac12}, \quad
d = -t_n^{\frac12}u_n^{-\frac12}.
\end{equation*}
To control the spacing of the formulas introduce a notation
\begin{equation}
(z_1,z_2,\ldots, z_k;q)_\infty
=(z_1;q)_\infty(z_2;q)_\infty\cdots (z_k;q)_\infty
\label{Pochprod}
%\tag{Pochprod}
\end{equation}
Then
$\ev_{q^\lambda t^\rho}(\Delta_d^Y(t^+)\Delta_{0,d}^Y(t^+)\Delta_g^Y(t^+))$
is equal to
\begin{align}
&\prod_{i=1}^n \frac{
(t_n^{\frac12}t_0^{\frac12}q^{\lambda_i}t^{n-i}t_n^{\frac12}t_0^{\frac12},
-t_n^{\frac12}t_0^{-\frac12}q^{\lambda_i}t^{n-i}t_n^{\frac12}t_0^{\frac12},
u_n^{\frac12}u_0^{\frac12}q^{\frac12} q^{\lambda_i}t^{n-i}t_n^{\frac12}t_0^{\frac12},
-u_n^{\frac12}u_0^{-\frac12}q^{\frac12} q^{\lambda_i}t^{n-i}t_n^{\frac12}t_0^{\frac12};q)_\infty
}{
(q^{2\lambda_i}t^{2(n-i)}t_nt_0,
qq^{2\lambda_i}t^{2(n-i)}t_nt_0;q^2)_\infty
}
\nonumber \\
&=\prod_{i=1}^n \frac{
(t_nt_0q^{\lambda_i}t^{n-i},
-t_n q^{\lambda_i}t^{n-i},
t_n^{\frac12}t_0^{\frac12}u_n^{\frac12}u_0^{\frac12}q^{\frac12} q^{\lambda_i}t^{n-i},
-t_n^{\frac12}t_0^{\frac12}u_n^{\frac12}u_0^{-\frac12}q^{\frac12} q^{\lambda_i}t^{n-i};q)_\infty
}{
(q^{2\lambda_i}t^{2(n-i)}t_nt_0;q)_\infty
}
\nonumber \\
&=\prod_{i=1}^n \frac{
(q^{-1}abcd q^{\lambda_i}t^{n-i},
cd q^{\lambda_i}t^{n-i},
ac q^{\lambda_i}t^{n-i},
bc q^{\lambda_i}t^{n-i};q)_\infty
}{
(q^{-1}abcd q^{2\lambda_i}t^{2(n-i)};q)_\infty
}
\label{pospart}
%\tag{pospart}
\end{align}
and
$\ev_{q^{-\lambda} t^{-\rho}}(\Delta^{Y^{-1}}_d((t^+)^{-1})\Delta^{Y^{-1}}_g((t^+)^{-1})$ is equal to
\begin{align}
&\prod_{i=1}^n \frac{
(t_n^{-\frac12}t_0^{-\frac12}q q^{\lambda_i}t^{n-i}t_n^{\frac12}t_0^{\frac12},
-t_n^{-\frac12}t_0^{\frac12}q q^{\lambda_i}t^{n-i}t_n^{\frac12}t_0^{\frac12},
u_n^{-\frac12}u_0^{-\frac12}q^{\frac12} q^{\lambda_i}t^{n-i}t_n^{\frac12}t_0^{\frac12},
-u_n^{-\frac12}u_0^{\frac12}q^{\frac12} q^{\lambda_i}t^{n-i}t_n^{\frac12}t_0^{\frac12};q)_\infty
}{
(q^2 q^{2\lambda_i}t^{2(n-i)}t_nt_0,
qq^{2\lambda_i}t^{2(n-i)}t_nt_0;q^2)_\infty
}
\nonumber \\
&=\prod_{i=1}^n \frac{
(qq^{\lambda_i}t^{n-i},
-t_0 q q^{\lambda_i}t^{n-i},
t_n^{\frac12}t_0^{\frac12}u_n^{-\frac12}u_0^{-\frac12} q^{\frac12} q^{\lambda_i}t^{n-i},
-t_n^{\frac12}t_0^{\frac12}u_n^{-\frac12}u_0^{\frac12}q^{\frac12} q^{\lambda_i}t^{n-i};q)_\infty
}{
(q q^{2\lambda_i}t^{2(n-i)}t_nt_0;q)_\infty
}
\nonumber 
\\
&=\prod_{i=1}^n \frac{
(q q^{\lambda_i}t^{n-i},
ab q^{\lambda_i}t^{n-i},
bd q^{\lambda_i}t^{n-i},
ad q^{\lambda_i}t^{n-i};q)_\infty
}{
(abcd  q^{2\lambda_i}t^{2(n-i)};q)_\infty
}.\label{negpart}
%\tag{negpart}
\end{align}

\begin{thm}
Define $N_\lambda(q,t^+)$ to be the product
\begin{align*}
N_\lambda(q,t^+) 
&= 
\ev_{q^\lambda t^\rho}\left(\frac{1}{\Delta_g^Y(t^+)\Delta_s^Y(t^+)\Delta_d^Y(t^+)\Delta_{0,s}^Y(t^+)\Delta_{0,d}^Y(t^+)}\right) 
\\ 
&\qquad\cdot
\ev_{q^{-\lambda} t^{-\rho}}\left(\frac{1}{\Delta^{Y^{-1}}_g((t^+)^{-1})\Delta^{Y^{-1}}_s((t^+)^{-1})\Delta^{Y^{-1}}_d((t^+)^{-1})}\right).
\end{align*}
Then
$$\langle P_\lambda(q,t^+), P_\lambda(q,t^+)\rangle_+ = N_\lambda(q,t^+).$$
\end{thm}
\begin{proof}
The proof is accomplished by verifying the following properties
\begin{enumerate}
\item[(a)] $N_\lambda(q,1,1,1,1,1) = 1$,
\item[(b)] $N_\lambda(q,t^+)$ is symmetric in the parameters $a,b,c,d$,
\item[(c)] $N_\lambda(q,t^+)$ satisfies the recursions
$$\frac{N_{\lambda}(q,t^\pm)}{N_{\lambda+\omega}(q,t^+)} = t_n^n \ev^t_{\lambda+\omega}\left(
\frac{\kappa^\pm_{v_{\lambda+\omega}}}{\kappa^+_{v_{\lambda+\omega}}}\right)
\qquad\hbox{and}\qquad
\frac{N_{\lambda}(q,t^\mp)}{N_{\lambda+\pi}(q,t^+)} = t^{n(n-1)}\ev^t_{\lambda+\pi}\left(
\frac{\kappa^\mp_{v_{\lambda+\pi}}}{\kappa^+_{v_{\lambda+\pi}}}\right)
$$
\end{enumerate}
Property (a) follows from the fact that $\Delta_d^Y(1,1,1,1,1)$, $\Delta_{0,d}^Y(1,1,1,1,1) $, $\Delta_g^Y(1,1,1,1,1) $, 
$\Delta_d^{Y^{-1}}(1,1,1,1,1)$ and $\Delta_g^{Y^{-1}}(1,1,1,1,1)$ are all equal to 1.

\smallskip\noindent
Property (b) follows from the fact that the product of the expressions in~\eqref{pospart} and~\eqref{negpart} is symmetric in
the parameters $a,b,c,d$.

\smallskip\noindent
Property (c) 
Changing from the parameters
$t^+ = (t,t_0^{\frac12},u_0^{\frac12}, t_n^{\frac12}, u_n^{\frac12})$
to the parameters
$t^\pm = (t,t_0^{\frac12},u_0^{\frac12}, q t_n^{\frac12}, u_n^{\frac12})$
replaces $t_n^{\frac12}$ by $qt_n^{\frac12}$ so that $c$ gets replaced by $qc$ and $d$ gets replaced by $qd$ 
and $a$ and $b$ and $t$ stay fixed. 
Write
\begin{align*}
N^{gd}_\lambda(q,t^+) 
&= 
\ev_{q^\lambda t^\rho}\left(\frac{1}{\Delta_g^Y(t^+)\Delta_d^Y(t^+)\Delta_{0,d}^Y(t^+)}\right) 
\cdot
\ev_{q^{-\lambda} t^{-\rho}}\left(\frac{1}{\Delta^{Y^{-1}}_g((t^+)^{-1})\Delta^{Y^{-1}}_d((t^+)^{-1})}\right),
\\
N^{s}_\lambda(q,t^+) 
&= 
\ev_{q^\lambda t^\rho}\left(\frac{1}{\Delta_s^Y(t^+)\Delta_{0,s}^Y(t^+)}\right) 
\cdot
\ev_{q^{-\lambda} t^{-\rho}}\left(\frac{1}{\Delta^{Y^{-1}}_s((t^+)^{-1})}\right),
\end{align*}
so that $N_\lambda(q,t^+) = N^{gd}_\lambda(q,t^+)N^s_\lambda(q,t^+)$.
Using the notation of~\eqref{Pochprod},
$$N^{gd}_\lambda(q,t^\pm) = \prod_{i=1}^n 
\frac{
(
q^2 abcd q^{2\lambda_i}t^{2(n-i)},
q^2 q^{-1}abcd q^{2\lambda_i}t^{2(n-i)};
q)_\infty
}{
\left(\begin{array}{l}
qq^{\lambda_i}t^{n-i},
q^2 q^{-1}abcd q^{\lambda_i}t^{n-i},
ab q^{\lambda_i}t^{n-i},
q ac q^{\lambda_i}t^{n-i}, \\
q ad q^{\lambda_i}t^{n-i},
q bc q^{\lambda_i}t^{n-i},
q bd q^{\lambda_i}t^{n-i},
q^2 cd q^{\lambda_i}t^{n-i}
\end{array}
; q\right)_\infty
}
$$
and
$$N^{gd}_{\lambda+\omega}(q,t^+) = \prod_{i=1}^n 
\frac{
(
abcd q^{2(\lambda_i+1)}t^{2(n-i)},
q^{-1}abcd q^{2(\lambda_i+1)}t^{2(n-i)};
q)_\infty
}{
\left(
\begin{array}{l}
qq^{\lambda_i+1}t^{n-i},
q^{-1}abcd q^{\lambda_i+1}t^{n-i},
ab q^{\lambda_i+1}t^{n-i},
ac q^{\lambda_i+1}t^{n-i}, \\
ad q^{\lambda_i+1}t^{n-i},
bc q^{\lambda_i+1}t^{n-i},
bd q^{\lambda_i+1}t^{n-i},
cd q^{\lambda_i+1}t^{n-i}
\end{array};
q\right)_\infty
}
$$
and
$$\frac{N^{gd}_\lambda(q,t^\pm)}{N^{gd}_{\lambda+\omega}(q,t^+)}
=\prod_{i=1}^n
\frac
{(1-abcd q^{\lambda_i}t^{n-i})(1-cdq^{\lambda_i+1}t^{n-i})}
{(1-qq^{\lambda_i}t^{n-i})(1-abq^{\lambda_i}t^{n-i})}
=\prod_{i=1}^n
\frac{(1-abcd q^{\lambda_i}t^{n-i})(1-cdq^{\lambda_i+1}t^{n-i})}
{(1-q^{\lambda_i+1}t^{n-i})(1-abq^{\lambda_i}t^{n-i})}.
$$
Since
$$\ev_{q^\lambda t^\rho}(\Delta_s^Y(t^+)\Delta_{0,s}^Y(t^+))
=\prod_{1\le i <j\le n} \frac{(q^{\lambda_i-\lambda_j}t^{j-i+1};q)_\infty (q^{\lambda_i+\lambda_j} t^{2n-i-j+1}t_0t_n;q)_\infty}
{(q^{\lambda_i-\lambda_j}t^{j-i};q)_\infty (q^{\lambda_i+\lambda_j} t^{2n-i-j}t_0t_n;q)_\infty}
$$
and 
$$\ev_{q^{-\lambda}t^{-\rho}}(\Delta_s^{Y^{-1}}(t^+)^{-1})
=\prod_{1\le i<j\le n} \frac{ (q^{\lambda_i-\lambda_j+1} t^{j-i-1};q)_\infty (q^{\lambda_i+\lambda_j+1} t^{2n-i-j-1}t_0t_n;q)_\infty}
{ (q^{\lambda_i-\lambda_j+1} t^{j-i};q)_\infty (q^{\lambda_i+\lambda_j+1} t^{2n-i-j}t_0t_n;q)_\infty}
$$
then
$$N^s_\lambda(q,t^\pm) = \prod_{1\le i<j\le n} \frac{ \left(\begin{array}{l} 
q^{\lambda_i-\lambda_j}t^{j-i}, q^{\lambda_i+\lambda_j}t^{2n-i-j}t_0q^2t_n,\\
q^{\lambda_i-\lambda_j+1}t^{j-i},q^{\lambda_i+\lambda_j+1}t^{2n-i-j}t_0q^2t_n\end{array} ;q\right)_\infty}
{\left(\begin{array}{l} q^{\lambda_i-\lambda_j}t^{j-i+1}, q^{\lambda_i+\lambda_j}t^{2n-i-j+1}t_0q^2t_n, \\
q^{\lambda_i-\lambda_j+1}t^{j-i-1},q^{\lambda_i+\lambda_j+1}t^{2n-i-j-1}t_0 q^2 t_n\end{array} ;q\right)_\infty}
$$
and
$$N^s_{\lambda+\omega}(q,t^+) = \prod_{1\le i<j\le n} \frac{ \left( \begin{array}{l}
q^{\lambda_i-\lambda_j}t^{j-i}, q^{\lambda_i+\lambda_j+2}t^{2n-i-j}t_0t_n,\\
q^{\lambda_i-\lambda_j+1}t^{j-i},q^{\lambda_i+\lambda_j+3}t^{2n-i-j}t_0t_n\end{array} ; q\right)_\infty}
{ \left(\begin{array}{l} q^{\lambda_i-\lambda_j}t^{j-i+1}, q^{\lambda_i+\lambda_j+2}t^{2n-i-j+1}t_0t_n, \\
q^{\lambda_i-\lambda_j+1}t^{j-i-1},q^{\lambda_i+\lambda_j+3}t^{2n-i-j-1}t_0t_n\end{array};q\right)_\infty},
$$
since $\omega_i = 1$.
So
$$\frac{N_\lambda^s(q,t^\pm)}{N_{\lambda+\omega}^s(q,t^+)} = 1.
$$
Since
\begin{align*}
t_n^n\cdot \ev^t_{\lambda+\omega} \left(\frac{\kappa^\pm_{v_{\lambda+\omega}} }{\kappa^+_{v_{\lambda+\omega}} } \right)
&=\ev^t_{\lambda+\omega}\left( \prod_{i=1}^n \frac{\kappa^{Y^{-1}}_{\varepsilon_i} }{\kappa^Y_{\varepsilon_i}}\right)
\\
&=\prod_{i=1}^n \ev^t_{\lambda+\omega}\left( t_n
\frac{(1-t_n^{\frac12}t_0^{\frac12}Y_i^{-1})(1+t_n^{\frac12}t_0^{-\frac12}Y_i^{-1})}{(1-Y_i^{-2}) }
\frac{(1-Y_i^2) } {(1-t_n^{\frac12}t_0^{\frac12}Y_i)(1+t_n^{\frac12}t_0^{-\frac12}Y_i)}\right)
\\
&=\prod_{i=1}^n \ev^t_{\lambda+\omega}\left( 
\frac{(1-t_n^{\frac12}t_0^{\frac12}Y_i^{-1})(1+t_n^{\frac12}t_0^{-\frac12}Y_i^{-1})}
 {(1-t_n^{-\frac12}t_0^{-\frac12}Y_i^{-1})(1+t_n^{-\frac12}t_0^{\frac12}Y_i^{-1})}\right)
\\
&=\prod_{i=1}^n \left( 
\frac{(1-t_n^{\frac12}t_0^{\frac12}q^{\lambda_i+1}t^{n-i} t_0^{\frac12}t_n^{\frac12} )
(1+t_n^{\frac12}t_0^{-\frac12}q^{\lambda_i+1}t^{n-i} t_0^{\frac12}t_n^{\frac12})}
{(1-t_n^{-\frac12}t_0^{-\frac12}q^{\lambda_i+1}t^{n-i} t_0^{\frac12}t_n^{\frac12})
(1+t_n^{-\frac12}t_0^{\frac12}q^{\lambda_i+1}t^{n-i} t_0^{\frac12}t_n^{\frac12})}\right)
\\
&=\prod_{i=1}^n \left( 
\frac{(1-t_nt_0 q^{\lambda_i+1}t^{n-i} )
(1+t_n q^{\lambda_i+1}t^{n-i} )}
{(1- q^{\lambda_i+1}t^{n-i} )
(1+t_0 q^{\lambda_i+1}t^{n-i} )}\right)
\\
&=\prod_{i=1}^n \left(
\frac{(1-q^{-1}abcd q^{\lambda_i+1}t^{n-i} )
(1-cd q^{\lambda_i+1}t^{n-i} )}
{(1- q^{\lambda_i+1}t^{n-i} )
(1-q^{-1}ab q^{\lambda_i+1}t^{n-i} )}\right)
\\
&=\prod_{i=1}^n \left(
\frac{(1- abcd q^{\lambda_i}t^{n-i} )
(1-cd q^{\lambda_i+1}t^{n-i} )}
{(1- q^{\lambda_i+1}t^{n-i} )
(1- ab q^{\lambda_i}t^{n-i} )}\right)
\end{align*}
then
$$
\frac{N_\lambda(q,t^\pm)}{N_{\lambda+\omega}(q,t^+)}
=\frac{N^{gd}_\lambda(q,t^\pm)}{N^{gd}_{\lambda+\omega}(q,t^+)}\cdot
\frac{N^s_\lambda(q,t^\pm)}{N^s_{\lambda+\omega}(q,t^+)}
=\frac{N^{gd}_\lambda(q,t^\pm)}{N^{gd}_{\lambda+\omega}(q,t^+)}\cdot 1
=t_n^n\cdot \ev^t_{\lambda+\omega} \left(\frac{\kappa^\pm_{v_{\lambda+\omega}} }{\kappa^+_{v_{\lambda+\omega}} } \right).
$$

The proof of the second equality in (c) is similar as follows.
$$N^{gd}_\lambda(q,t^\mp) = \prod_{i=1}^n 
\frac{
(
q^2 abcd q^{2\lambda_i}(qt)^{2(n-i)},
q^2 q^{-1}abcd q^{2\lambda_i}(qt)^{2(n-i)};
q)_\infty
}{
\left(\begin{array}{l}
qq^{\lambda_i}(qt)^{n-i},
q^2 q^{-1}abcd q^{\lambda_i}(qt)^{n-i},
ab q^{\lambda_i}(qt)^{n-i},
q ac q^{\lambda_i}(qt)^{n-i}, \\
q ad q^{\lambda_i}(qt)^{n-i},
q bc q^{\lambda_i}(qt)^{n-i},
q bd q^{\lambda_i}(qt)^{n-i},
q^2 cd q^{\lambda_i}(qt)^{n-i}
\end{array}
; q\right)_\infty
}
$$
and
$$N^{gd}_{\lambda+\pi}(q,t^+) = \prod_{i=1}^n 
\frac{
(
abcd q^{2(\lambda_i+n-i)}t^{2(n-i)},
q^{-1}abcd q^{2(\lambda_i+n-i)}t^{2(n-i)};
q)_\infty
}{
\left(
\begin{array}{l}
qq^{\lambda_i+n-i}t^{n-i},
q^{-1}abcd q^{\lambda_i+n-i}t^{n-i},
ab q^{\lambda_i+n-i}t^{n-i},
ac q^{\lambda_i+n-i}t^{n-i}, \\
ad q^{\lambda_i+n-i}t^{n-i},
bc q^{\lambda_i+n-i}t^{n-i},
bd q^{\lambda_i+n-i}t^{n-i},
cd q^{\lambda_i+n-i}t^{n-i}
\end{array};
q\right)_\infty
}
$$
so that
$$\frac{N_\lambda^{gd}(q,t^\mp)}{N_\lambda^{gd}(q,t^+)} = 1.$$
Then
$$N^s_\lambda(q,t^\mp) = \prod_{1\le i<j\le n} \frac{ \left(\begin{array}{l} 
q^{\lambda_i-\lambda_j+j-i}t^{j-i}, q^{\lambda_i+\lambda_j+2n-i-j}t^{2n-i-j}t_0t_n,\\
q^{\lambda_i-\lambda_j+j-i+1}t^{j-i},q^{\lambda_i+\lambda_j+1+2n-i-j}t^{2n-i-j}t_0t_n\end{array} ;q\right)_\infty}
{\left(\begin{array}{l} q^{\lambda_i-\lambda_j+j-i+1}t^{j-i+1}, q^{\lambda_i+\lambda_j+2n-i-j+1}t^{2n-i-j+1}t_0t_n, \\
q^{\lambda_i-\lambda_j+j-i}t^{j-i-1},q^{\lambda_i+\lambda_j+2n-i-j}t^{2n-i-j-1}t_0t_n\end{array} ;q\right)_\infty}
$$
and
$$N^s_{\lambda+\pi}(q,t^+) = \prod_{1\le i<j\le n} \frac{ \left( \begin{array}{l}
q^{\lambda_i-\lambda_j+j-i}t^{j-i}, q^{\lambda_i+\lambda_j+2n-i-j}t^{2n-i-j}t_0t_n,\\
q^{\lambda_i-\lambda_j+j-i+1}t^{j-i},q^{\lambda_i+\lambda_j+1+2n-i-j}t^{2n-i-j}t_0t_n\end{array} ; q\right)_\infty}
{ \left(\begin{array}{l} q^{\lambda_i-\lambda_j+j-i}t^{j-i+1}, q^{\lambda_i+\lambda_j+2n-i-j}t^{2n-i-j+1}t_0t_n, \\
q^{\lambda_i-\lambda_j+j-i+1}t^{j-i-1},q^{\lambda_i+\lambda_j+1+2n-i-j}t^{2n-i-j-1}t_0t_n\end{array};q\right)_\infty},
$$
since $\pi_i = n-i$.
So
$$\frac{N^s_\lambda(q,t^\mp)}{N^s_{\lambda+\pi}(t^+)}
= \prod_{1\le i<j\le n} \frac{(1-q^{\lambda_i-\lambda_j+j-i}t^{j-i+1})(1-q^{\lambda_i+\lambda_j+2n-i-j}t^{2n-i-j+1}t_0t_n)}
{(1-q^{\lambda_i-\lambda_j+j-i}t^{j-i-1})(1-q^{\lambda_i+\lambda_j+2n-i-j}t^{2n-i-j-1}t_0t_n)}
$$
and
$$
\frac{N_\lambda(q,t^\mp)}{N_{\lambda+\pi}(t^+)}
=\frac{N^{gd}_\lambda(q,t^\mp)}{N^{gd}_{\lambda+\pi}(t^+)} \cdot \frac{N^s_\lambda(q,t^\mp)}{N^s_{\lambda+\pi}(t^+)}
=1\cdot \frac{N^s_\lambda(q,t^\mp)}{N^s_{\lambda+\pi}(t^+)}
=t^{n(n-1)} \ev^t_{\lambda+\pi} \left(\frac{\kappa^\mp_{v_{\lambda+\pi}} }{\kappa^+_{v_{\lambda+\pi}} } \right),
$$
since
\begin{align*}
&t^{n(n-1)} \ev^t_{\lambda+\pi} \left(\frac{\kappa^\mp_{v_{\lambda+\pi}} }{\kappa^+_{v_{\lambda+\pi}} } \right)
=\ev^t_{\lambda+\pi}\left( \prod_{1\le i<j\le n} t^2
\frac{\kappa^{Y^{-1}}_{\varepsilon_i-\varepsilon_j}\kappa^{Y^{-1}}_{\varepsilon_i+\varepsilon_j} }
{\kappa^Y_{\varepsilon_i-\varepsilon_j}\kappa^Y_{\varepsilon_i+\varepsilon_j} }
\right)
\\
&=\ev^t_{\lambda+\pi}\left( \prod_{1\le i<j\le n} t^2
\frac{(1-tY_i^{-1}Y_j)(1-tY_i^{-1}Y_j^{-1}) }{ (1-Y_i^{-1}Y_j)(1-Y_i^{-1}Y_j^{-1}) }
\frac{(1-Y_iY_j^{-1})(1-Y_iY_j) }{ (1-tY_iY_j^{-1})(1-tY_iY_j)}
\right)
\\
&=\ev^t_{\lambda+\pi}\left( \prod_{1\le i<j\le n} 
\frac{(1-tY_i^{-1}Y_j)(1-tY_i^{-1}Y_j^{-1}) }{ (1-t^{-1}Y^{-1}_iY_j)(1-t^{-1}Y^{-1}_iY^{-1}_j)}
\right)
\\
&= \prod_{1\le i<j\le n} 
\frac{(1-tq^{\lambda_i+(n-i)-\lambda_j-(n-j)}t^{-(n-j)+(n-i)} )(1-tq^{\lambda_i+(n-i)+\lambda_j+(n-j)}t^{(n-i)+(n-j)}t_0t_n) }
{(1-t^{-1}q^{\lambda_i+(n-i)-\lambda_j-(n-j)}t^{-(n-j)+(n-i)} )(1-t^{-1}q^{\lambda_i+(n-i)+\lambda_j+(n-j)}t^{(n-i)+(n-j)}t_0t_n) }
\\
&=\prod_{1\le i<j\le n} 
\frac{(1-tq^{\lambda_i-\lambda_j+j-i}t^{j-i} )(1-tq^{\lambda_i+\lambda_j+2n-i-j)}t^{2n-i-j}t_0t_n) }
{(1-t^{-1}q^{\lambda_i-\lambda_j+j-i}t^{j-i} )(1-t^{-1}q^{\lambda_i+\lambda_j+2n-i-j}t^{2n-i-j}t_0t_n) }.
\end{align*}
\end{proof}

\subsection{The constant term}

To get the constant term of $\dfrac{1}{\Delta^+_{CC}}$ specialize $\lambda=0$.
Then
$$N_0^{gd}(t^+) = 
\prod_{i=1}^n \frac{ (abcd t^{2(n-i)},q^{-1}abcd t^{2(n-i)};q)_\infty}
{(qt^{n-i}, q^{-1}abcd t^{n-i}, abt^{n-i}, ac t^{n-i}, ad t^{n-i}, bc t^{n-i}, bd t^{n-i}, cd t^{n-i};q)_\infty}
$$
and
$$N_0^s(t^+) = \prod_{1\le i<j\le n}
\frac{(t^{j-i}, t^{2n-i-j}t_0t_n, qt^{j-i},qt^{2n-i-j}t_0t_n;q)_\infty}
{(t t^{j-i},t t^{2n-i-j}t_0t_n, qt^{-1} t^{j-i}, qt^{-1} t^{2n-i-j}t_0t_n;q)_\infty},
$$
so that
$$W_0(t,t_n) \mathrm{ct}\left(\frac{1}{\Delta^+_{CC}}\right) 
= W_0(t,t_n)(1,1)_+
= \langle 1, 1\rangle_+ = N_0(t^+) = N_0^{gd}(t^+)N_0^s(t^+).$$

\newpage

\section{Appendix: Examples for $n=2$}

\subsection{The finite Weyl group and the roots in $S_0^+$}

The finite Weyl group is
$W_{\mathrm{fin}} = \{ 1, s_1, s_2, s_1s_2, s_2s_1, s_1s_2s_1, s_2s_1s_2, s_1s_2s_1s_2\}$
where
$$s_1^2=s_2^2=1, \quad s_1s_2s_1s_2=s_2s_1s_2s_1.$$
Then
$$S_{0,d}^+ = \{ \varepsilon_1, \varepsilon_2\}
\qquad\hbox{and}\qquad
S_{0,s}^+ = \{ \varepsilon_1-\varepsilon_2, \varepsilon_1+\varepsilon_2\}$$
with $s_{\varepsilon_1} = s_1s_2s_1$, $s_{\varepsilon_2} = s_2$, $s_{\varepsilon_1-\varepsilon_2} = s_1$,
$s_{\varepsilon_1+\varepsilon_2} = s_2s_1s_2$.

\subsection{Hecke symmetrizers}

The Hecke symmetrizers are
\begin{align*}
\varepsilon_+ &= T_1T_2T_1T_2 + t^{-\frac12}T_2T_1T_2+t_n^{-\frac12}T_1T_2T_1
\\
&\qquad
+t^{-\frac12}t_n^{-\frac12}T_1T_2
+t^{-\frac12}t_n^{-\frac12}T_2T_1
+t^{-1} t_n^{-\frac12}T_2+t^{-\frac12} t_n^{-1} T_1+ t^{-1} t_n^{-1}, \\
\varepsilon_\pm
&= T_1T_2T_1T_2 + t^{-\frac12}T_2T_1T_2-t_n^{\frac12}T_1T_2T_1
\\
&\qquad
-t^{-\frac12}t_n^{\frac12}T_1T_2
-t^{-\frac12}t_n^{\frac12}T_2T_1
-t^{-1} t_n^{\frac12}T_2+t^{-\frac12} t_n T_1+ t^{-1} t_n, \\
\varepsilon_\mp
&= T_1T_2T_1T_2 - t^{\frac12}T_2T_1T_2 + t_n^{-\frac12}T_1T_2T_1
\\
&\qquad
-t^{\frac12}t_n^{-\frac12}T_1T_2
-t^{\frac12}t_n^{-\frac12}T_2T_1
+t t_n^{-\frac12}T_2-t^{\frac12} t_n^{-1} T_1+ t t_n^{-1}, \\
\varepsilon_- 
&= T_1T_2T_1T_2 - t^{\frac12}T_2T_1T_2-t_n^{\frac12}T_1T_2T_1
+t^{\frac12}t_n^{\frac12}T_1T_2
+t^{\frac12}t_n^{\frac12}T_2T_1
-t t_n^{\frac12}T_2-t^{\frac12} t_n T_1+ t t_n.
\end{align*}
Since
$$\varepsilon_{\pm} 
= (T_1+t^{-\frac12})(T_2T_1T_2-t_n^{\frac12}T_2T_1-t^{-\frac12}t_n^{\frac12}T_2+t^{-\frac12}t_n)
= (T_1+t^{-\frac12})(T_2T_1-t^{-\frac12}t_n^{\frac12})(T_2-t_n^{\frac12})
$$
and $(T_1-t^{\frac12})(T_1+t^{-\frac12})=0$ then $T_1\varepsilon_{\pm} = t^{\frac12}\varepsilon_\pm$.
Similarly, since 
$$\varepsilon_{\pm} 
= (T_2-t_n^{\frac12})(T_1T_2T_1+t^{-\frac12}T_1T_2-t^{-\frac12}t_n^{\frac12}T_1-t^{-1}t_n^{\frac12})
= (T_2-t_n^{\frac12})(T_1T_2 - t^{-\frac12}t_n^{\frac12})(T_1+t^{-\frac12})
$$
and $(T_2+t_n^{-\frac12})(T_2-t_n^{\frac12}) = 0$ then $T_2\varepsilon_{\pm}=-t_n^{-\frac12}\varepsilon_\pm$.
Then
\begin{align*}
\varepsilon_\pm^2 
&= (tt_n^{-1}+t_n^{-1} + t + 1+1 + t^{-1} +t_n +t^{-1}t_n)\varepsilon_\pm
\\
&= t^{-1}t_2 (t^2t_2^{-2}+tt_2^{-2}+t^2t_2^{-1} + t t_2^{-1}+ t t_2^{-1}+t_2^{-1} +t +1)\varepsilon_\pm
\\
&= t^{-1}t_2W_0(t,t_2^{-1})\varepsilon_\pm.
\end{align*}
since
$$W_0(t,t_n) = 1+t_2+t+tt_2+tt_2+tt_2^2+t^2t_2+t^2t_2^2 = (1+t_2)(1+t)(1+tt_2).
$$

\subsection{$c$-functions}

The $c$-functions in $X$ for roots in $S_{0,s}^+$ and $S_{0,d}^+$ are
$$c_{\varepsilon_1-\varepsilon_2}^X = \frac{t^{-\frac12}-t^{\frac12}X_1X_2^{-1}}{1-X_1X_2^{-1}},
\qquad
c_{\varepsilon_1+\varepsilon_2}^X = \frac{t^{-\frac12}-t^{\frac12}X_1X_2 }{ 1-X_1 X_2 },
$$
$$c_{\varepsilon_1}^X = t_2^{-\frac12} \frac{(1-t_2^{\frac12}u_2^{\frac12}X_1)(1+t_2^{\frac12}u_2^{-\frac12}X_1)}{1-X_1^2},
\qquad
c_{\varepsilon_2}^X = t_2^{-\frac12} \frac{(1-t_2^{\frac12}u_2^{\frac12}X_2)(1+t_2^{\frac12}u_2^{-\frac12}X_2)}{1-X_2^2}.
$$
Then
$$c_\mathrm{dr}^{X^{-1}} = c_{\varepsilon_1}^{X^{-1}}c_{\varepsilon_2}^{X^{-1}}
= t_2^{-\frac{2}{2}} 
\frac{(1-t_2^{\frac12}u_2^{\frac12}X^{-1}_1)(1+t_2^{\frac12}u_2^{-\frac12}X^{-1}_1) }{ (1-X_1^{-2})}
\frac{(1-t_2^{\frac12}u_2^{\frac12}X^{-1}_2)(1+t_2^{\frac12}u_2^{-\frac12}X^{-1}_2) }{ (1-X_2^{-2})}
$$
and
$$c_\mathrm{st}^{X^{-1}} = c_{\varepsilon_1-\varepsilon_2}^{X^{-1}}c_{\varepsilon_1+\varepsilon_2}^{X^{-1}}
= \frac{(t^{-\frac12}-t^{\frac12}X_1^{-1}X_2)}{(1-X_1^{-1}X_2)}
\frac{(t^{-\frac12}-t^{\frac12}X_1^{-1}X_2^{-1})}{(1-X_1^{-1}X_2^{-1})}.
$$

\subsection{Weyl denominators}

Since $\omega = (1,1)$, $\pi = (1,0)$ and $\rho = (2,1)$ then
$$x^\omega = x_1x_2, \qquad x^\pi = x_1, \qquad x^\rho = x_1^2x_2.$$
Since
\begin{align*}
e_\pm &= s_1s_2s_1s_2-s_1s_2s_1+s_2s_1s_2-s_1s_2-s_2s_1+s_1-s_2+1, \\
e_\mp &= s_1s_2s_1s_2+s_1s_2s_1-s_2s_1s_2-s_1s_2-s_2s_1-s_1+s_2+1, \\
e_- &= s_1s_2s_1s_2-s_1s_2s_1-s_2s_1s_2+s_1s_2+s_2s_1-s_1-s_2+1,
\end{align*}
then the Weyl denominators for $n=2$ are
\begin{align*}
a^\pm_\omega &= \hbox{$\frac12$} e_\pm x^\omega = e_\pm x_1x_2 
= x_1x_2-x_1x_2^{-1}+x_1^{-1}x_2-x_1^{-1}x_2^{-1} 
= x_1^{-1}x_2^{-1}(1-x_1^{2})(1-x_2^{2}), \\
a^\mp_\pi &= \hbox{$\frac12$}e_\mp x^\delta = e_\mp x_1 
=  x_1-x_2-x_2^{-1}+x_1^{-1} =  x_1^{-1}(1-x_1x_2)(1-x_1x_2^{-1}), \\
a_\rho &=e_+  x^\rho = e_+x_1^2x_2 
= x_1^2x_2-x_1x_2^2-x_1^2x_2^{-1}+x_1x_2^{-2}+x_1^{-2}x_2-x_1^{-1}x_2^2 - x_1^{-1}x_2^{-2} + x_1^{-2}x_2^{-1}
\\
&=x_1^{-2}x_2^{-1}(1-x_1x_2)(1-x_1x_2^{-1})(1-x_1^{2})(1-x_2^{2}) = a^\pm_\omega a^\mp_\pi, 
\end{align*}
and
\begin{align*}
P_0 &= 1,
\\
A^{\pm}_\omega 
&= A^\pm_{-\varepsilon_1-\varepsilon_2}
= x_1^{-1} x_2^{-1} (1-t_2^{\frac12}u_2^{\frac12}x_1)(1+t_2^{\frac12}u_2^{-\frac12}x_1)
(1-t_2^{\frac12}u_2^{\frac12}x_2)(1+t_2^{\frac12}u_2^{-\frac12}x_2)
\\
A^{\mp}_\pi&= A^\mp_{-\varepsilon_1} =  x_1^{-1} (1-tx_1x_2^{-1})(1-tx_1x_2)
= x_1^{-1}-tx_2^{-1}-tx_2+t^2x_1,
\\
A_\rho &= x_1^{-2}x_2^{-1}
(1-tx_1x_2^{-1})(1-tx_1x_2)
(1-t_2^{\frac12}u_2^{\frac12}x_1)(1+t_2^{\frac12}u_2^{-\frac12}x_1)
(1-t_2^{\frac12}u_2^{\frac12}x_2)(1+t_2^{\frac12}u_2^{-\frac12}x_2).
\end{align*}

\subsection{Examples of Proposition~\ref{cfcnsymm},~\ref{symmexp} and~\ref{Eexps}}

Since $\chi^+(T_{s_1s_2s_1s_2}) = tt_2$ then
\begin{align*}
t t_n \varepsilon_+
&=
(\eta_{s_1s_2s_1s_2}+\eta_{s_1s_2s_1}+\eta_{s_2s_1s_2}+\eta_{s_1s_2}+\eta_{s_2s_1}+\eta_{s_1}+\eta_{s_2}+1)
\kappa_{\varepsilon_1}^Y\kappa_{\varepsilon_2}^Y 
\kappa_{\varepsilon_1-\varepsilon_2}^Y \kappa_{\varepsilon_1+\varepsilon_2}^Y 
\\
&=
(\eta_{s_2s_1s_2}+\eta_{s_2s_1s_2}+\eta_{s_1s_2}+\eta_{s_2}+1)
\kappa_{\varepsilon_1}^Y\kappa_{\varepsilon_2}^Y \kappa_{\varepsilon_1+\varepsilon_2}^Y 
(1+\eta_{s_1}) \kappa_{\varepsilon_1-\varepsilon_2}^Y  
\\
&=
\left(\begin{array}{l}
\eta_{s_2s_1s_2}\kappa_{\varepsilon_1}^Y\kappa_{\varepsilon_2}^Y \kappa_{\varepsilon_1+\varepsilon_2}^Y 
+\kappa_{\varepsilon_2}^Y \eta_{s_1s_2}\kappa_{\varepsilon_2}^Y \kappa_{\varepsilon_1+\varepsilon_2}^Y \\
+ \kappa_{\varepsilon_1}^Y \kappa_{\varepsilon_1+\varepsilon_2}^Y \eta_{s_2}\kappa_{\varepsilon_2}^{Y^{-1}}
+\kappa_{\varepsilon_1}^Y\kappa_{\varepsilon_2}^Y \kappa_{\varepsilon_1+\varepsilon_2}^Y \end{array}\right)
\chi^+(T_{s_1}) \varepsilon^+_\omega,
\\
t t_n \varepsilon_\pm 
&=
\kappa_{\varepsilon_1}^{Y^{-1}}\kappa_{\varepsilon_2}^{Y^{-1}} 
(\eta_{s_1s_2s_1s_2}-\eta_{s_1s_2s_1}+\eta_{s_2s_1s_2}-\eta_{s_1s_2}-\eta_{s_2s_1}+\eta_{s_1}-\eta_{s_2}+1)
\kappa_{\varepsilon_1-\varepsilon_2}^Y \kappa_{\varepsilon_1+\varepsilon_2}^Y
\\
&=
\kappa_{\varepsilon_1}^{Y^{-1}} \kappa_{\varepsilon_2}^{Y^{-1}} 
(\eta_{s_2s_1s_2}-\eta_{s_1s_2}-\eta_{s_2}+1)\kappa_{\varepsilon_1+\varepsilon_2}^Y 
(1+\eta_{s_1}) \kappa_{\varepsilon_1-\varepsilon_2}^Y 
\\
&=
\left(\begin{array}{l}
\eta_{s_2s_1s_2}\kappa_{\varepsilon_1}^Y \kappa_{\varepsilon_2}^Y\kappa_{\varepsilon_1+\varepsilon_2}^Y
-\kappa_{\varepsilon_2}^{Y^{-1}} \eta_{s_1s_2}\kappa_{\varepsilon_2}^Y \kappa_{\varepsilon_1+\varepsilon_2}^Y \\
-\kappa_{\varepsilon_1}^{Y^{-1}}  \kappa_{\varepsilon_1-\varepsilon_2}^Y \eta_{s_2}\kappa_{\varepsilon_2}^Y 
+\kappa_{\varepsilon_1}^{Y^{-1}} \kappa_{\varepsilon_2}^{Y^{-1}} \kappa_{\varepsilon_1+\varepsilon_2}^Y\end{array}\right) 
\chi^+(T_{s_1}) \varepsilon^+_{\omega},
\\
t t_n \varepsilon_\mp 
&=
\kappa_{\varepsilon^\vee_1-\varepsilon^\vee_2}^{Y^{-1}} \kappa_{\varepsilon^\vee_1+\varepsilon^\vee_2}^{Y^{-1}}
(\eta_{s_1s_2s_1s_2}+\eta_{s_1s_2s_1}-\eta_{s_2s_1s_2}-\eta_{s_1s_2}-\eta_{s_2s_1}-\eta_{s_1}+\eta_{s_2}+1)
\kappa_{\varepsilon_1^\vee}^Y \kappa_{\varepsilon_2^\vee}^Y 
\\
&=
\kappa_{\varepsilon^\vee_1-\varepsilon^\vee_2}^{Y^{-1}} \kappa_{\varepsilon^\vee_1+\varepsilon^\vee_2}^{Y^{-1}} 
(\eta_{s_1s_2s_1}-\eta_{s_2s_1}-\eta_{s_1}+1)\kappa_{\varepsilon_1^\vee}^Y 
(1+\eta_{s_2}) \kappa_{\varepsilon_2^\vee}^Y 
\\
&=
\left(\begin{array}{l}
\eta_{s_1s_2s_1}\kappa_{\varepsilon^\vee_1}^Y \kappa_{\varepsilon^\vee_1-\varepsilon^\vee_2}^Y
\kappa_{\varepsilon^\vee_1+\varepsilon^\vee_2}^Y
-\kappa_{\varepsilon^\vee_1-\varepsilon^\vee_2}^{Y^{-1}} \eta_{s_2s_1}
\kappa_{\varepsilon^\vee_1}^Y \kappa_{\varepsilon^\vee_1-\varepsilon^\vee_2}^Y
\\
\quad
-\kappa_{\varepsilon^\vee_1+\varepsilon^\vee_2}^{Y^{-1}}  \kappa_{\varepsilon^\vee_2}^{Y} 
\eta_{s_1}\kappa_{\varepsilon^\vee_1-\varepsilon_2^\vee}^Y 
+ \kappa_{\varepsilon^\vee_1-\varepsilon^\vee_2}^{Y^{-1}} 
\kappa_{\varepsilon^\vee_1+\varepsilon^\vee_2}^{Y^{-1}}\kappa_{\varepsilon^\vee_1}^Y
\end{array}\right)
\cdot
\chi^+(T_{s_2}) \varepsilon^+_{\pi},
\\
t t_n \varepsilon_-
&=
\kappa_{\varepsilon_1}^{Y^{-1}}\kappa_{\varepsilon_2}^{Y^{-1}}
\kappa_{\varepsilon_1-\varepsilon_2}^{Y^{-1}} \kappa_{\varepsilon_1+\varepsilon_2}^{Y^{-1}} 
(\eta_{s_1s_2s_1s_2}-\eta_{s_1s_2s_1}-\eta_{s_2s_1s_2}+\eta_{s_1s_2}+\eta_{s_2s_1}-\eta_{s_1}-\eta_{s_2}+1)
\\
&=
\eta_{s_1s_2s_1s_2}
\kappa_{\varepsilon_1}^Y\kappa_{\varepsilon_2}^Y
\kappa_{\varepsilon_1-\varepsilon_2}^Y \kappa_{\varepsilon_1+\varepsilon_2}^Y 
-\kappa_{\varepsilon_2}^{Y^{-1}} \eta_{s_1s_2s_1}
\kappa_{\varepsilon_1}^Y
\kappa_{\varepsilon_1-\varepsilon_2}^Y \kappa_{\varepsilon_1+\varepsilon_2}^Y 
-\kappa_{\varepsilon_1-\varepsilon_2}^{Y^{-1}} \eta_{s_2s_1s_2}
\kappa_{\varepsilon_1}^Y\kappa_{\varepsilon_2}^Y
\kappa_{\varepsilon_1+\varepsilon_2}^Y 
\\
&\qquad
+\kappa_{\varepsilon_1}^{Y^{-1}}
\kappa_{\varepsilon_1-\varepsilon_2}^{Y^{-1}}
\eta_{s_1s_2}
\kappa_{\varepsilon_2}^Y \kappa_{\varepsilon_1+\varepsilon_2}^Y
+\kappa_{\varepsilon_2}^{Y^{-1}}
\kappa_{\varepsilon_1+\varepsilon_2}^{Y^{-1}}
\eta_{s_2s_1}
\kappa_{\varepsilon_1}^Y \kappa_{\varepsilon_1-\varepsilon_2}^Y
\\
&\qquad
-\kappa_{\varepsilon_2}^{Y^{-1}}
\kappa_{\varepsilon_1-\varepsilon_2}^{Y^{-1}} \kappa_{\varepsilon_1+\varepsilon_2}^{Y^{-1}}
\eta_{s_1}
\kappa_{\varepsilon_1-\varepsilon_2}^Y
-\kappa_{\varepsilon_1}^{Y^{-1}}
\kappa_{\varepsilon_1-\varepsilon_2}^{Y^{-1}} \kappa_{\varepsilon_1+\varepsilon_2}^{Y^{-1}} 
\eta_{s_2}\kappa_{\varepsilon_2}^Y
+\kappa_{\varepsilon_1}^{Y^{-1}}\kappa_{\varepsilon_2}^{Y^{-1}}
\kappa_{\varepsilon_1-\varepsilon_2}^{Y^{-1}} \kappa_{\varepsilon_1+\varepsilon_2}^{Y^{-1}}.
\end{align*}
Since
\begin{align*}
\kappa_{s_1}^Y &= \kappa_{\varepsilon_1^\vee-\varepsilon^\vee_2}^Y,
&\kappa_{s_1s_2}^Y &= \kappa_{\varepsilon_2^\vee}^Y \kappa_{\varepsilon_1^\vee+\varepsilon_2^\vee}^Y,
&\kappa_{s_1s_2s_1}^Y &= \kappa_{\varepsilon_1^\vee-\varepsilon_2^\vee}^Y \kappa_{\varepsilon_1^\vee}^Y
\kappa_{\varepsilon_1^\vee+\varepsilon_2^\vee}^Y.
\\
\kappa_{s_2}^Y &= \kappa_{\varepsilon_2^\vee}^Y, 
&\kappa_{s_2s_1}^Y &= \kappa_{\varepsilon_1^\vee-\varepsilon_2^\vee}^Y \kappa_{\varepsilon_1^\vee}^Y,
&\kappa_{s_2s_1s_2}^Y &= \kappa_{\varepsilon_2^\vee}^Y \kappa_{\varepsilon_1^\vee+\varepsilon_2^\vee}^Y
\kappa_{\varepsilon_1^\vee}^Y,
\end{align*}
then
\begin{align*}
A^\pm_\omega 
&= E_{-\varepsilon_1-\varepsilon_2}
-\ev^t_{-\varepsilon_1+\varepsilon_2}(\kappa^{Y^{-1}}_{\varepsilon_2^\vee})E_{-\varepsilon_1+\varepsilon_2}
-\ev^t_{\varepsilon_1-\varepsilon_2}(\kappa^{Y^{-1}}_{\varepsilon_2^\vee}
\kappa^Y_{\varepsilon_1^\vee-\varepsilon_2^\vee})E_{\varepsilon_1-\varepsilon_2}
\\\
&\qquad
+\ev^t_{\varepsilon_1+\varepsilon_2}(\kappa^{Y^{-1}}_{\varepsilon_1^\vee}\kappa^{Y^{-1}}_{\varepsilon_2^\vee}
\kappa^Y_{\varepsilon_1^\vee+\varepsilon_2^\vee})E_{\varepsilon_1+\varepsilon_2}, \qquad\hbox{and} \\
A^\mp_\pi
&= E_{-\varepsilon_1}-\ev^t_{-\varepsilon_2}(\kappa^{Y^{-1}}_{\varepsilon_1^\vee-\varepsilon_2^\vee})E_{-\varepsilon_2}
-\ev^t_{\varepsilon_2}(\kappa^{Y^{-1}}_{\varepsilon_1^\vee+\varepsilon_2^\vee}\kappa^Y_{\varepsilon_2^\vee}) E_{\varepsilon_2}
+\ev^t_{\varepsilon_1}(\kappa_{\varepsilon^\vee_1-\varepsilon^\vee_2}^{Y^{-1}} 
\kappa_{\varepsilon^\vee_1+\varepsilon^\vee_2}^{Y^{-1}}\kappa_{\varepsilon^\vee_1}^Y)E_{\varepsilon_1}.
\end{align*}


\newpage

\begin{thebibliography}{20}




\bibitem[BT72]{BT72} F.\ Bruhat and J.\ Tits, \emph{Groupes r\'eductifs sur un corps local: I. Donn\'ees
radicielles valu\'ees}, Publications Math\'ematiques de l'Institut des Hautes \'Etudes Scientifiques,
no. \textbf{41} (1972) 5-251.



\bibitem[CGdGW16]{CGdGW16} L.\ Cantini, A.\ Garbali, J.\ de Gier and M.\ Wheeler,
\emph{Koornwinder polynomials and the stationary multi-species asymmetric exclusion process with open boundaries},
J.\ Phys.\ A49(2016), no.44, 444002, 23 pp.

\bibitem[C03]{C03} I.V.\ Cherednik, \emph{Double affine Hecke algebras and difference Fourier transforms}, 
Invent.\ Math.\ \textbf{152} (2003), no. 2, 213--303. 


\bibitem[CR22]{CR22} L.\ Colmenarejo and A.\ Ram, \emph{$c$-functions and Macdonald polynomials}, 
J.\  Algebra \textbf{655} (2024) 163-222,
https://doi.org/10.1016/j.jalgebra.2024.01.024.


\bibitem[CMW23]{CMW23} S.\ Corteel, O.\ Mandelshtam and L.\ Williams,
\emph{Rhombic staircase tableaux and Koornwinder polynomials},  arXiv:2312.17469. 


\bibitem[HKO24]{HKO24} G. Heckman, T. Koornwinder and E.\ Opdam, 
\emph{Memories of Ian G. Macdonald},  arxiv:2410.07882.

\bibitem[Hel94]{Hel94} S.\ Helgason, \emph{Harish-Chandra's $c$-function. A mathematical jewel}
in \textsl{Noncompact Lie groups and some of their applications (San Antonio, TX, 1993)}, 55-67, NATO Adv.\ Sci.\ Inst.\ Ser.\ C: Math.\ Phys.\ Sci., \textbf{429}, Kluwer Acad.\ Publ., Dordrecht, 1994.
Reprinted in Proc.\ Sympos.\ Pure Math.\ \textbf{68} Amer.\ Math.\ Soc., Providence, RI, 2000, 273-283. ,


\bibitem[Kac]{Kac} V. Kac, {\em Infinite dimensional Lie algebras}, Third edition. Cambridge University Press, Cambridge, 1990. xxii+400 pp. 

\bibitem[K92]{K92} T.H.\ Koornwinder, \emph{Askey-Wilson polynomials for root systems of type $BC$}, in: ``Hypergeometric functions on domains of positivity, Jack polynomials, and applications'' (Tampa, FL, 1991), 189–204, Contemp. Math.\ \textbf{138} Amer. Math. Soc., Providence, RI, 1992.


\bibitem[Mac]{Mac} I.G.\ Macdonald, 
\textsl{Symmetric functions and Hall polynomials}, 
Second edition, Oxford Mathematical Monographs, Oxford University Press, New York, 1995. 


\bibitem[Mac72]{Mac72} I.G.\ Macdonald,
\emph{Affine root systems and Dedekind's $\eta$-function},
Invent.\ Math\ \textbf{15} (1972) 91-143.


\bibitem[Mac87]{Mac87} I.G.\ Macdonald,
\emph{Orthogonal polynomials associated with root systems}, preprint 1987,
S\'em.\ Lothar.\ Combin.\ \textbf{45} (2000/01), Art\. B45a, 40 pp.




\bibitem[Mac03]{Mac03} I.G.\ Macdonald,  {\sl Affine Hecke Algebras and Orthogonal Polynomials}, 
Cambridge Tracts in Mathematics, vol. \textbf{157}, Cambridge University Press, Cambridge, 2003. 



\bibitem[Nou95]{Nou95} M.\ Noumi, \emph{Macdonald-Koornwinder polynomials and affine Hecke rings},
S\=uriseisekikenky\=usho K\=oky\=uroku \textbf{919}(1995) 44--55 (in Japanese).

\bibitem[Ra01]{Ra01} E.\ Rains, \emph{$BC_n$-symmetric polynomials},
Transform.\ Groups \textbf{10} (2005), no.\ 1, 63-132.

\bibitem[Ra17]{Ra17} E.\ Rains, \emph{Elliptic double affine Hecke algebras},
SIGMA Symmetry Integrability Geom.\ Methods Appl.\ \textbf{16} (2020), Paper No.\ 111, 133 pp.

\bibitem[RW15]{RW15} 
E.\ Rains and S.O.\ Warnaar, 
\emph{Bounded Littlewood identities},
Mem.\ Amer.\ Math.\ Soc.\ \textbf{270} (2021), no.\ 1317, vii+115 pp.

\bibitem[Sah99]{Sah99} S.\ Sahi, \emph{Nonsymmetric Koornwinder polynomials and duality}, Ann.\ 
Math.\ \textbf{150} (1999) no.\ 1, 267--282.

\bibitem[Sah00]{Sah00} S.\ Sahi, \emph{Some properties of Koornwinder polynomials},
Contemp.\ Math.\ \textbf{254} (2000) 395-411.

\bibitem[Sto00]{Sto00} J.V.\ Stokman, \emph{Koornwinder polynomials and affine Hecke algebras},
Internat.\ Math.\ Res.\ Notices (2000) 1005–1042.


\bibitem[Sto04]{Sto04} J.V.\ Stokman, \emph{Lecture notes on Koornwinder polynomials}
Laredo Lectures on Orthogonal Polynomials and Special Functions, 145–207, 
Adv.\ Theory Spec.\ Funct.\ Orthogonal Polynomials, Nova Science 2004.

\bibitem[Sto11]{Sto11} J.V.\ Stokman, \emph{The c-function expansion of a basic hypergeometric function associated to root systems},
Ann.\ of Math.\ (2) \textbf{179} (2014)  253-299.

\bibitem[Sto21]{Sto21}
J.V.\ Stokman, \emph{Macdonald-Koornwinder polynomials},
Encyclopedia of special functions: the Askey-Bateman project. Vol. 2. Multivariable special functions, 258-313,
Cambridge University Press 2021.


\bibitem[vD96]{vD96} J. F.\ van Diejen, \emph{Self-dual Koornwinder-Macdonald polynomials}, Inventiones Mathematicae, \textbf{126} (1996).


\bibitem[vD95]{vD95} J. F.\ van Diejen, \emph{Commuting difference operators with polynomial eigenfunctions}, 
Compositio Mathematica \textbf{95} (1995).



\bibitem[Yam20]{Yam20} K.\ Yamaguchi, \emph{Littlewood-Richardson rule for Koornwinder
polynomials},  J.\ Algebraic Combin.\ \textbf{56} (2022) 335–381.


\bibitem[YY21]{YY21} K. Yamaguchi and S. Yanagida, \emph{Specializing Koornwinder polynomials to Macdonald polynomials of type B,C,D and BC}, J.\ Algebraic Combin.\ \textbf{57} (2023) 171-226.



\end{thebibliography}
\end{document}